%% file: teth_e10.tex
\documentclass[a4paper,12pt]{article}
\usepackage{latexsym,amsthm,hhline}
\begin{document}

\newtheorem{define}{Definition}
\newtheorem{Lemma}{Lemma}
\newtheorem{note}{Remark}
\newtheorem{Theorem}{Theorem}
\newtheorem{coll}{Corollary}[Theorem]

\def \H {{I\!\!H}}
\def \E {{I\!\!E}}
\def \Z {{Z\!\!\!Z}}

\begin{center}
{\Large\bf
Coxeter Decompositions of Hyperbolic Tetrahedra.
}

\bigskip
{\large
A.~Felikson
%\footnote{\scriptsize
%Supported in part by  RFFI grant n 96-15-96050.}
}

\vspace{35pt}
\parbox{10.5cm}
{\scriptsize
{\bf Abstract}.
In this paper, we classify Coxeter decompositions of hyperbolic
tetrahedra, i.e. simplices in the hyperbolic space $\H^3$.
The paper together with \cite{simp4-9} completes the classification of
Coxeter decompositions of hyperbolic simplices.

}

\end{center}

\section{Introduction}

\begin{define}
A convex polyhedron in a space of constant curvature
is called a {\bf Coxeter polyhedron}
if all  dihedral angles of this polyhedron are the integer parts of
$\pi$.

\end{define}

\begin{define}
\label{def1}
A {\bf Coxeter decomposition} of a convex polyhedron $P$
is a decomposition of $P$ into finitely many tiles
such that each tile is a
Coxeter polyhedron and
any two tiles having a common facet are symmetric
 with respect to this facet.
\end{define}

Coxeter decompositions of hyperbolic
 triangles were studied
in~\cite{Pink}, \cite{K}, \cite{Mat}, \cite{M} and \cite{Deza}.
Coxeter decompositions of hyperbolic simplices of the dimension
greater than three were classified in~\cite{simp4-9}.
In this paper, we classify Coxeter decompositions of hyperbolic
tetrahedra, i.e. simplices in the hyperbolic space $\H^3$.
The paper completes the classification of the Coxeter decompositions of
hyperbolic simplices.

%{\bf Acknowledgments}\\

The work was partially written during the stay at the University of Bielefeld.
The author is grateful to the hospitality of this
University.
The author would like to thank O.~V.~Schwarzman and
E.~B.~Vinberg for their interest to the work.

\subsection{Definitions }

The tiles in  Definition~\ref{def1} are called
{\bf fundamental polyhedra}.
Clearly, any two fundamental polyhedra are congruent to each other.
A hyperplane $\alpha$ containing a facet of a fundamental polyhedron
is called a {\bf mirror} if $\alpha$ contains no facet of $P$.

\begin{define}
Given a Coxeter decomposition of a polyhedron $P$,
a {\bf dihedral angle} of $P$ formed up by facets
$\alpha$ and $\beta$ is called {\bf fundamental}
if no mirror contains $\alpha \cap \beta$.
A vertex $A$ is called {\bf fundamental}
if no mirror contains $A$.
\end{define}

\pagebreak

\noindent
{\bf Notation.}

\noindent
$P$ is a tetrahedron equipped with a Coxeter decomposition; \\
$F$ is a fundamental polyhedron
considered up to an isometry of $\H^3$; \\
$\Sigma(T)$ is a Coxeter diagram of a Coxeter tetrahedron $T$; \\
$N$ is the number of the fundamental polyhedra inside $P$.\\

%\noindent
%A decomposition is called {\bf non-trivial} if $N>1$.

%$k$-edge is an edge of the fundamental tetrahedron with dihedral
%angle $\frac{\pi}{k}$.\\
%$k$-$l$-$m$-vertex is a vertex of the fundamental tetrahedron,
%such that the dihedral angles at the edges intersecting in this
%vertex are equal to  $\frac{\pi}{k}$, $\frac{\pi}{l}$ and
%$\frac{\pi}{m}$.\\

A Coxeter tetrahedron $T$ can be
represented by its Coxeter diagram $\Sigma(T)$:
the nodes $v_i$ of  $\Sigma(T)$ correspond to the facets $f_i$ of $T$,
two nodes $v_i$ and $v_j$ are connected by a $k$-fold edge if
the dihedral angle formed up by $f_i$ and $f_j$ equals $\frac{\pi}{k}$.

There are finitely many Coxeter hyperbolic tetrahedra.
We denote a Coxeter hyperbolic tetrahedron by $H_i$, where $i$
is a number of this tetrahedron in
Table~\ref{h_vol}.

\section{Types of Coxeter decompositions}

It is shown in~\cite{simp4-9} that the large class of decompositions
of simplices can be listed using the {\bf inductive algorithm}
explained in section~1 of~\cite{simp4-9}.

\begin{define}
Let $\Theta(P)$ be a Coxeter decomposition of a simplex $P$.

$\Theta(P)$ is called a decomposition of the
{\bf first type} if the decomposition can be obtained by the inductive
algorithm.

$\Theta(P)$ is called a decomposition of the
{\bf second type} if all the dihedral angles of $P$ are fundamental.

$\Theta(P)$ is called a decomposition of the
{\bf third type} if $\Theta(P)$ is neither of the first type nor of the
second.

\end{define}

The decompositions of the first type
where obtained using a computer.
See Table~\ref{pic_3_h1}--\ref{pic_3_h2}
for the  list of the decompositions.
%see \cite{arhiv-tetr}.
%
%
The decompositions of the second type are studied in section~\ref{type2}.
It was proved in section~2.3 of~\cite{simp4-9} that any decomposition of
the third type is a superposition of some decompositions.
In Lemma~\ref{l-type3} we prove that no hyperbolic tetrahedron admits
a decomposition of the third type.

\section{Decompositions without non-fundamental dihedral angles}
\label{type2}

Suppose that $P$ admits a Coxeter decomposition such that any dihedral
angle of $P$ is fundamental.
In this case $P$ is a Coxeter tetrahedron.

The following properties are proved in~\cite{simp4-9}:

\begin{Lemma}
\label{properties}
Let $P$ be a simplex in $\H^n$ admitting a Coxeter decomposition of the
second type with fundamental polyhedron $F$.
Then the following properties are hold:

\begin{itemize}

\item[1.]
$F$ is a simplex.

\item[2.]
If all the vertices of $P$ are fundamental,
then  $N\ge 2^n$.

\item[3.]
{\bf Volume property}.\\
$\frac{Vol(P)}{Vol(F)}\in \mathbf Z$, where $Vol(T)$ is a volume
of a simplex $T$.

\item[4.]
{\bf Subdiagram property}. \\
Let $\Sigma(P)$ be a Coxeter diagram of $P$ and $\Sigma(F)$ be a Coxeter
diagram of $F$.  Let $v$ be a node of $\Sigma(P)$.
Then there exists a node $w$ of $\Sigma(F)$ such that
either
$\Sigma(P)\setminus v = \Sigma(F) \setminus w$
or the simplex  $p$ determined by the subdiagram $\Sigma(P) \setminus v$
admits a Coxeter decomposition of the second type with fundamental simplex
$f$ determined by $\Sigma(F) \setminus w$.

\end{itemize}

\end{Lemma}

All the volumes of hyperbolic Coxeter simplices are
calculated in~\cite{Ruth1}.
We reprint the list of volumes of hyperbolic Coxeter tetrahedra
in Table~\ref{h_vol}.
In the dimensions greater than three the Volume property together with
the Subdiagram property are nearby a sufficient condition of the existence
of a decomposition with given  $F$ and $P$.
In particular, all but two pairs $(F,P)$ satisfying these properties
correspond to some decompositions of the second type.
In the dimension three the properties are satisfied by a large number
of pairs $(F,P)$. Thus, we need some additional considerations in this
case.

There are finitely many bounded Coxeter tetrahedra
and finitely many unbounded Coxeter tetrahedra
in $\H^3$.
It is evident that $F$ and $P$ should be bounded or
unbounded simultaneously.

\begin{Lemma}\label{h_2_b}
No bounded hyperbolic tetrahedron admits a Coxeter decomposition
of the second type.

\end{Lemma}

\begin{proof}
There is a unique pair of bounded
tetrahedra $(F,P)$ satisfying the Volume property.
These tetrahedra are $F=H_1$ and $P=H_3$.
There exists a unique decomposition with ($F$,$P$)=($H_1$,$H_3$),
and in this decomposition $N=2$. Any decomposition with $N=2$ has a
non-fundamental angle, and the lemma is proved.

\end{proof}

Consider the unbounded tetrahedra.

\begin{Lemma}  \label{unb_f}
Let $P$ be an unbounded tetrahedron.
Then $P$ has a non-funda\-mental vertex.

\end{Lemma}

\begin{proof}
Suppose that any vertex of $P$ belongs to a unique fundamental
tetrahedron. Since each fundamental tetrahedron has
an ideal vertex, $N\le 4$.
This contradicts to the second part of Lemma~\ref{properties}.

\end{proof}

To prove the following
lemma we need a list of Coxeter decompositions
of Euclidean and spherical triangles $p$ such that all the
angles of $p$ are fundamental and $N\le 24$.
See Fig.~\ref{fund} for the list.

\begin{figure}[!h]
\begin{center}
\caption{
Decompositions of spherical and Euclidean triangles}
without non-fundamental angles
 ($N\le 24$ in Euclidean case)
\input pic/fund.tex
\label{fund}
\vspace{15pt}

{\scriptsize
 The decompositions \ref{fund}a--\ref{fund}e are Euclidean, the
 decomposition \ref{fund}f is spherical.  \\
 The numbers of
 fundamental triangles are written under the decompositions.

}
\end{center}
\end{figure}

\begin{Lemma}  \label{h_2_f}
Suppose that any dihedral angle of $P$
is fundamental.
Then the decomposition is one of two
decompositions shown in Fig.~\ref{f3}.

\end{Lemma}

\begin{proof}
 By Lemma~\ref{h_2_b} the tetrahedron $P$ is  unbounded.
By Lemma~\ref{unb_f} $P$ has a non-fundamental
vertex.  Let $A$ be a non-fundamental vertex of $P$.

Suppose that $A$ is not an ideal vertex.
Consider a small sphere $s$ centered in $A$.
The section of the decomposition by the sphere $s$ is
a Coxeter decomposition of a spherical triangle $p=s\bigcap P$.
Clearly, all angles of $p$ are fundamental.
The only such a decomposition consists of 15 triangles
(see Fig.~\ref{fund}f). Hence, the number of tiles
in $P$ cannot be less than 15.
The fundamental triangle in the decomposition shown in
Fig.~\ref{fund}f has an angle $\frac{\pi}{5}$.
Hence, the fundamental tetrahedron has a dihedral angle
equal to
$\frac{\pi}{5}$.  But there is no pair of unbounded hyperbolic
tetrahedra $F$ and $P$ such that $\frac{Vol(P)}{Vol(F)}\ge 15$
and $F$ has a dihedral angle  equal to $\frac{\pi}{5}$.
Therefore, any non-fundamental vertex should be ideal.

Let $A$ be an ideal non-fundamental vertex of  $P$.
Consider a small horosphere centered in $A$.
The section of the decomposition by this horosphere
is a Coxeter decomposition of a
Euclidean triangle $p$
without non-fundamental angles.
 Since the maximal ratio of volumes of unbounded
hyperbolic tetrahedra is 24, the number of triangular
tiles $f$ in  $p$ cannot be greater than
24.  Therefore, the decomposition of the triangle $p$ is one of
the decompositions
shown in Fig.~\ref{fund}a--\ref{fund}e.  Consider these decompositions of
$p$.

\begin{itemize}

\item[1)]

Suppose that $p$ is decomposed as shown in
Fig.~\ref{fund}a--\ref{fund}c.
Then all the angles of the triangles $f$ and $p$ are equal to
$\frac{\pi}{3}$. Therefore, each of the Coxeter diagrams
$\Sigma (F)$ and $\Sigma (P)$ has a subdiagram corresponding to
a triangle with angles
($\frac{\pi}{3}$,$\frac{\pi}{3}$,$\frac{\pi}{3}$).
All hyperbolic Coxeter tetrahedra
satisfying this condition are
 shown in Fig.~\ref{3.3.}.
The tetrahedra $H_{24}$ and $H_{32}$ cannot be fundamental
for a Coxeter decomposition of the second type. Indeed,
for any Coxeter tetrahedron $T$ we have
$\frac{Vol(T)}{Vol(H_{24})}<2$ and
$\frac{Vol(T)}{Vol(H_{32})}<2$.
The maximal ratio of volumes of the tetrahedra shown in
Fig.~\ref{3.3.}a--\ref{3.3.}c is 12.
Therefore, $p$ cannot be decomposed as the triangle shown in
Fig.\ref{fund}c.
Thus, $F$ is one of the tetrahedra shown in Fig.~\ref{3.3.}a,
$P$ is one of the tetrahedra shown in
Fig.~\ref{3.3.}a--\ref{3.3.}c and $p$ is decomposed as shown in
Fig.\ref{fund}a or \ref{fund}b.

\begin{figure}[!h]
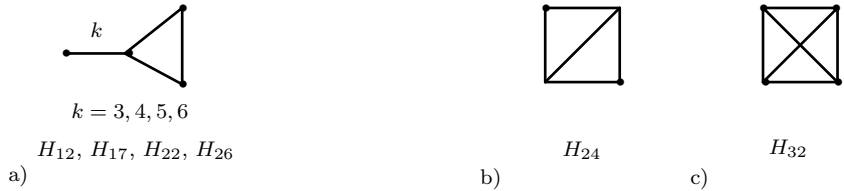

\begin{center}
{\scriptsize
\input pic/3_3.tex
}

\end{center}
\caption{Diagrams of tetrahedra with subdiagram
corresponding to
the triangle
($\frac{\pi}{3}$,$\frac{\pi}{3}$,$\frac{\pi}{3}$).
}
\label{3.3.}
\end{figure}

Now we are aimed to show that no of $H_{17}$, $H_{22}$ and $H_{26}$ can
be a fundamental polyhedron for a decomposition of the second type.
Indeed, suppose the contrary,
i.e. $\Sigma (F)$ is shown in Fig.~\ref{3.3.}a and $k>3$.
Then $P$ has no dihedral angle which is equal to $\frac{\pi}{k}$,
since $\Sigma (P)\ne \Sigma (F)$.
Consider the set of fundamental tetrahedra having a common
ideal vertex whose neighborhood is decomposed as shown in
Fig.~\ref{fund}a--\ref{fund}c.
Each of these tetrahedra has its own edge with dihedral angle
$\frac{\pi}{k}$. Since $P$ has no such a dihedral angle,
each of these edges belongs to $k$ distinct tetrahedra.
Therefore, $N\ge n\cdot k$, where $n$ is the number of the tiles
$f$ in $p$, i.e. $n=4$ or $9$.\\
If $n=9$, then $n\cdot k\ge 9\cdot 4>24$. This is impossible. \\
If $n=4$, then $n\cdot k\ge 4\cdot 4=16$; but there is no Coxeter
tetrahedron $P$ having the volume big enough.
Therefore, $k=3$.

The tetrahedron $H_{12}$
has no dihedral angle which is equal to $\frac{\pi}{l}$,
where $l=4$,$5$,$6$.
Thus, $H_{12}$
can tessellate only the
tetrahedra $H_{24}$ and $H_{32}$.
These tessellations are the decompositions shown
in Fig.~\ref{f3}.

\item[2)]
Suppose that $p$ is decomposed as shown in Fig.~\ref{fund}d.
Then the number of fundamental tetrahedra in the decomposition
cannot be smaller than nine. Further, each of the
Coxeter diagrams $\Sigma (F)$ and $\Sigma (P)$ has a
subdiagram corresponding to a triangle  with angles
($\frac{\pi}{4}$,$\frac{\pi}{4}$,$\frac{\pi}{2}$).
A unique pair of tetrahedra ($F$,$P$) satisfying these
conditions is ($H_{11}$,$H_{31}$). To prove that there is no such
a decomposition, consider
an ideal vertex $A$ of $P$. Let $F_1,...,F_9$ be the
fundamental tetrahedra
incident to $A$. Let
$f_i$ ($i=1,...,9$) be a face of $F_i$ opposite to $A$.
Since $F_i$ has a
face which is not orthogonal to $f_i$, the faces $f_i$ ($i=1$,...,9)
belong to three different
planes (see Fig.~\ref{3_planes}). Any of these planes contains at most
four faces $f_i$. Therefore, the decomposition contains at least five
fundamental tetrahedra additional to $F_1$,...,$F_9$.
This is impossible, since $\frac{Vol(P)}{Vol(F)}=12$ and $9+5=14> 12$.

\begin{figure}[!h]
{\scriptsize
\begin{center}
\input pic/3_planes.tex

\caption{
The triangles shaded by distinct ways belong to
distinct planes.}
{\normalsize
These planes intersect each other and form up
dihedral angle $\frac{2\pi}{3}$.}
\label{3_planes}
\end{center}

}
\end{figure}

\item[3)]
Suppose that $p$ is decomposed as shown in
Fig.~\ref{fund}e. Then the number of the fundamental tetrahedra
in the decomposition of $P$ cannot be smaller than 18.  Further,
$\Sigma (F)$  has a subdiagram corresponding to
a triangle  with angles
($\frac{\pi}{3}$,$\frac{\pi}{3}$,$\frac{\pi}{3}$)
and  $\Sigma (P)$ has a subdiagram
corresponding to a triangle  with angles
($\frac{\pi}{2}$,$\frac{\pi}{3}$,$\frac{\pi}{6}$).
A unique pair of unbounded hyperbolic tetrahedra ($F$,$P$)
satisfying these conditions is ($H_{10}$,$H_{32}$).
Since any dihedral angle of $P$ equals
$\frac{\pi}{3}$,  the neighborhood of any vertex of $P$
is decomposed as shown in Fig.~\ref{fund}e. Then
$N\ge 18\cdot 4>24$ ($F$ has a unique ideal vertex).
Since $\frac{Vol(P)}{Vol(F)}=24$, such a decomposition does not
exist.

\end{itemize}

All the cases are considered, the lemma is proved.

\end{proof}

\begin{Lemma}
\label{l-type3}
There is no Coxeter decomposition of hyperbolic tetrahedron
of the third type.

\end{Lemma}

\begin{proof}

It was proved in section~2.3 of~\cite{simp4-9} that any decomposition of
the third type is a superposition of some decompositions.
Let $P$ be a hyperbolic tetrahedron and $\Theta(P)$ be a decomposition of
the third type.  Suppose that $P$ is a smallest tetrahedron with the
decomposition of the third type, that is for any tetrahedron $P'$ inside
$P$ the restriction of $\Theta(P)$ onto $P'$ is a decomposition either of
the first or of the second type.

Let $\Theta_1$ and $\Theta_2$ be
decompositions of the first type and $\Theta(P)$ be a superposition of
$\Theta_1$ and $\Theta_2$.  It follows from the construction of the
inductive algorithm that $\Theta(P)$ is a decomposition of the first type.
The contradiction shows that at least one of $\Theta_1$ and $\Theta_2$
is a decomposition of the second type, i.e. at least one of two
decompositions shown in Fig.~\ref{f3}.  A superposition of any of these
two decompositions with any other decomposition is a decomposition of the
first type.

\end{proof}

%\pagebreak
%\clearpage

\section*{Tables}

%\subsubsection*{Notation}

The tetrahedra with dihedral angles
$\frac{k_i\pi}{q_i}$ $i=1$,...,6 are represented by the
following diagrams:
take a Coxeter diagram of the tetrahedron whose dihedral angles
equal $\frac{\pi}{q_i}$ $i=1$,...,6 and subdivide the edges
corresponding to the dihedral angles  $\frac{k_i\pi}{q_i}$ into $k_i$
parts. \\

The list of decompositions of the first type was obtained
by the inductive algorithm
(see section~1 of \cite{simp4-9}).
All the non-trivial decompositions of the first
type are listed here for each of the Coxeter tetrahedra.
The Coxeter diagram of the fundamental tetrahedron
is presented before the list of  decompositions with
this fundamental tetrahedron. \\
Numbers ($k$,$l$ ; $m$,$n$,$p$,$q$) under the diagrams
describe the decompositions.
\begin{itemize}
\item[]
$\bullet$
$k$ is a number of fundamental tetrahedra in the
decomposition;\\
$\bullet$
$l$ is a number of gluing necessary for construction of the
decomposition;
(gluing the tetrahedra with $l=l_1$ and $l=l_2$
we obtain a tetrahedron with $l=1+max \{l_1,l_2\}$);\\
$\bullet$
$m$ and $n$ are the numbers of the tetrahedra which
should be glued together to obtain the decomposition; \\
$\bullet$
$p$ and $q$ are the numbers of the faces glued together of the tetrahedra
$m$ and $n$ respectively,
the nodes of the diagrams are
numbered from the left to the right by the numbers 0,1,2,3.

\end{itemize}

\pagebreak
\clearpage

\begin{table}[!h]
\begin{center}
\caption{Hyperbolic Coxeter  tetrahedra.}
{ The volumes are reprinted from \cite{Ruth1}.}
%{\scriptsize
%\input pic/h_vol.tex

\end{center}
\label{h_vol}
\end{table}

\begin{center}
\input pic/vol_3_og.tex

\pagebreak

\input pic/vol_3ne1.tex

\pagebreak

\input pic/vol_3ne2.tex
\end{center}

\pagebreak
\clearpage

\begin{table}[!h]
\begin{center}
\caption{
Decompositions of bounded tetrahedra of the first type.}
\label{pic_3_h1}
\end{center}
\vspace{-45pt}
\input{pic/pic_3h1.tex}
\end{table}
\vspace{20pt}

\begin{table}[!h]
\begin{center}
\caption{
Decompositions of unbounded tetrahedra of the first type.}
\label{pic_3_h2}
\end{center}
\vspace{-45pt}
\end{table}

\noindent
\input{pic/pic_3h2.tex}

%\pagebreak
%\clearpage

\vspace{60pt}

\begin{figure}[!h]
\begin{center}
\caption{Hyperbolic tetrahedra of the second type.}
{(Ideal vertices are marked by small circles.)}
\vspace{14pt}
\input pic/f3.tex

\label{f3}
\vspace{6pt}
\end{center}
\end{figure}

%\pagebreak
%\clearpage

\end{document}

%% file: pic/fund.tex
%TeX representation of Charton picture (version 1.96)
%Copyright (c) S.Trifonov, D.Tejblum, 1993-4
%FILE: C:\XBOCT\SIMPLEX\PIC\FUND.tex
%DATE: Wed May 20 15:34:09 1998
%------------
\setlength{\unitlength}{0.075in}
\begin{picture}(69.35,14.00)
\put(57.35,0.16){{\setbox0=\hbox{f)}\lower\ht0\box0}}
\put(33.88,2.83){{\setbox0=\hbox{9}\lower\ht0\box0}}
\special{em:linewidth 0.014in}
\put(34.90,4.16){\special{em:moveto}}
\put(34.90,8.65){\special{em:lineto}}
\put(37.05,6.40){\special{em:moveto}}
\put(34.90,4.16){\special{em:lineto}}
\put(32.55,6.40){\special{em:lineto}}
\put(37.05,4.16){\special{em:moveto}}
\put(37.05,6.40){\special{em:lineto}}
\put(32.36,6.40){\special{em:lineto}}
\put(34.90,8.65){\special{em:lineto}}
\put(32.55,8.65){\special{em:lineto}}
\put(32.55,10.90){\special{em:moveto}}
\put(32.55,4.16){\special{em:lineto}}
\put(39.40,4.16){\special{em:lineto}}
\put(32.55,10.90){\special{em:lineto}}
\put(62.81,1.83){{\setbox0=\hbox{15}\lower\ht0\box0}}
\put(42.41,0.00){{\setbox0=\hbox{e)}\lower\ht0\box0}}
\put(68.95,6.66){\special{em:moveto}}
\put(67.75,5.46){\special{em:lineto}}
\put(66.23,4.26){\special{em:lineto}}
\put(65.61,3.83){\special{em:lineto}}
\put(59.48,8.00){\special{em:moveto}}
\put(68.95,6.66){\special{em:lineto}}
\put(61.75,12.00){\special{em:moveto}}
\put(64.41,11.66){\special{em:lineto}}
\put(66.55,11.16){\special{em:lineto}}
\put(65.61,3.83){\special{em:moveto}}
\put(61.75,12.00){\special{em:lineto}}
\put(61.48,3.83){\special{em:moveto}}
\put(60.76,4.83){\special{em:lineto}}
\put(60.21,5.83){\special{em:lineto}}
\put(59.81,6.83){\special{em:lineto}}
\put(59.48,8.00){\special{em:lineto}}
\put(66.55,11.16){\special{em:moveto}}
\put(61.48,3.83){\special{em:lineto}}
\put(66.55,11.33){\special{em:moveto}}
\put(66.75,8.80){\special{em:lineto}}
\put(66.65,7.21){\special{em:lineto}}
\put(66.35,5.63){\special{em:lineto}}
\put(65.98,4.43){\special{em:lineto}}
\put(65.75,3.83){\special{em:lineto}}
\put(65.75,3.83){\special{em:moveto}}
\put(63.91,4.58){\special{em:lineto}}
\put(62.18,5.58){\special{em:lineto}}
\put(60.53,6.91){\special{em:lineto}}
\put(59.48,8.00){\special{em:lineto}}
\put(59.48,8.00){\special{em:moveto}}
\put(61.63,9.41){\special{em:lineto}}
\put(63.88,10.48){\special{em:lineto}}
\put(66.55,11.33){\special{em:lineto}}
\put(63.75,14.00){\special{em:moveto}}
\put(64.93,13.05){\special{em:lineto}}
\put(66.16,11.78){\special{em:lineto}}
\put(67.13,10.51){\special{em:lineto}}
\put(67.90,9.25){\special{em:lineto}}
\put(68.63,7.56){\special{em:lineto}}
\put(69.13,5.88){\special{em:lineto}}
\put(69.35,4.66){\special{em:lineto}}
\put(69.35,4.66){\special{em:moveto}}
\put(67.45,4.11){\special{em:lineto}}
\put(65.51,3.80){\special{em:lineto}}
\put(63.88,3.71){\special{em:lineto}}
\put(61.35,3.91){\special{em:lineto}}
\put(59.28,4.38){\special{em:lineto}}
\put(58.41,4.66){\special{em:lineto}}
\put(58.41,4.66){\special{em:moveto}}
\put(58.96,6.90){\special{em:lineto}}
\put(59.88,9.13){\special{em:lineto}}
\put(60.86,10.81){\special{em:lineto}}
\put(62.20,12.50){\special{em:lineto}}
\put(63.75,14.00){\special{em:lineto}}
\put(31.61,0.00){{\setbox0=\hbox{d)}\lower\ht0\box0}}
\put(19.61,0.00){{\setbox0=\hbox{c)}\lower\ht0\box0}}
\put(9.08,0.16){{\setbox0=\hbox{b)}\lower\ht0\box0}}
\put(0.01,0.16){{\setbox0=\hbox{a)}\lower\ht0\box0}}
\put(47.08,2.50){{\setbox0=\hbox{18}\lower\ht0\box0}}
\put(24.55,2.16){{\setbox0=\hbox{16}\lower\ht0\box0}}
\put(12.95,2.16){{\setbox0=\hbox{9}\lower\ht0\box0}}
\put(4.55,2.00){{\setbox0=\hbox{4}\lower\ht0\box0}}
\put(48.68,8.83){\special{em:moveto}}
\put(46.68,10.00){\special{em:lineto}}
\put(49.35,4.00){\special{em:moveto}}
\put(51.08,5.00){\special{em:lineto}}
\put(44.68,4.00){\special{em:moveto}}
\put(44.68,6.16){\special{em:lineto}}
\put(46.68,4.00){\special{em:lineto}}
\put(48.81,8.83){\special{em:lineto}}
\put(46.15,8.83){\special{em:lineto}}
\put(49.35,4.00){\special{em:lineto}}
\put(50.28,6.33){\special{em:lineto}}
\put(44.68,6.16){\special{em:lineto}}
\put(48.68,8.83){\special{em:moveto}}
\put(44.68,6.16){\special{em:lineto}}
\put(49.35,4.00){\special{em:lineto}}
\put(48.68,8.83){\special{em:lineto}}
\put(47.35,11.00){\special{em:moveto}}
\put(43.35,4.00){\special{em:lineto}}
\put(51.88,4.00){\special{em:lineto}}
\put(47.35,11.00){\special{em:lineto}}
\put(23.21,10.33){\special{em:moveto}}
\put(24.28,12.50){\special{em:lineto}}
\put(29.35,3.66){\special{em:lineto}}
\put(26.95,3.66){\special{em:lineto}}
\put(28.01,6.00){\special{em:lineto}}
\put(25.75,6.00){\special{em:lineto}}
\put(26.81,8.33){\special{em:lineto}}
\put(24.41,8.33){\special{em:lineto}}
\put(25.48,10.50){\special{em:lineto}}
\put(23.21,10.50){\special{em:lineto}}
\put(22.15,8.33){\special{em:moveto}}
\put(23.21,10.33){\special{em:lineto}}
\put(27.08,3.66){\special{em:lineto}}
\put(24.68,3.66){\special{em:lineto}}
\put(25.75,6.00){\special{em:lineto}}
\put(23.35,6.00){\special{em:lineto}}
\put(24.28,8.33){\special{em:lineto}}
\put(22.15,8.33){\special{em:lineto}}
\put(21.08,6.00){\special{em:moveto}}
\put(23.35,6.00){\special{em:lineto}}
\put(22.28,3.66){\special{em:lineto}}
\put(21.08,6.00){\special{em:lineto}}
\put(22.15,8.33){\special{em:moveto}}
\put(19.88,3.66){\special{em:lineto}}
\put(24.68,3.66){\special{em:lineto}}
\put(22.15,8.33){\special{em:lineto}}
\put(11.75,8.33){\special{em:moveto}}
\put(12.81,10.33){\special{em:lineto}}
\put(16.68,3.66){\special{em:lineto}}
\put(14.28,3.66){\special{em:lineto}}
\put(15.35,6.00){\special{em:lineto}}
\put(12.95,6.00){\special{em:lineto}}
\put(13.88,8.33){\special{em:lineto}}
\put(11.75,8.33){\special{em:lineto}}
\put(10.68,6.00){\special{em:moveto}}
\put(12.95,6.00){\special{em:lineto}}
\put(11.88,3.66){\special{em:lineto}}
\put(10.68,6.00){\special{em:lineto}}
\put(11.75,8.33){\special{em:moveto}}
\put(9.48,3.66){\special{em:lineto}}
\put(14.28,3.66){\special{em:lineto}}
\put(11.75,8.33){\special{em:lineto}}
\put(3.61,6.00){\special{em:moveto}}
\put(5.88,6.00){\special{em:lineto}}
\put(4.81,3.66){\special{em:lineto}}
\put(3.61,6.00){\special{em:lineto}}
\put(4.68,8.33){\special{em:moveto}}
\put(2.41,3.66){\special{em:lineto}}
\put(7.21,3.66){\special{em:lineto}}
\put(4.68,8.33){\special{em:lineto}}
\end{picture}

%% file: pic/3_3.tex
%TeX representation of Charton picture (version 1.96)
%Copyright (c) S.Trifonov, D.Tejblum, 1993-4
%FILE: 3_3.tex
%DATE: Thu Oct 01 16:12:56 1998
%------------
\setlength{\unitlength}{0.075in}
\begin{picture}(57.88,11.16)
\put(2.01,1.83){{\setbox0=\hbox{$H_{12}$, $H_{17}$, $H_{22}$, $H_{26}$}\lower\ht0\box0}}
\put(38.68,1.83){{\setbox0=\hbox{$H_{24}$}\lower\ht0\box0}}
\put(53.08,2.00){{\setbox0=\hbox{$H_{32}$}\lower\ht0\box0}}
\put(0.01,0.16){{\setbox0=\hbox{a)}\lower\ht0\box0}}
\put(32.95,0.00){{\setbox0=\hbox{b)}\lower\ht0\box0}}
\special{em:linewidth 0.014in}
\put(52.68,11.16){\special{em:moveto}}
\put(57.88,6.00){\special{em:lineto}}
\put(52.81,6.00){{\setbox0=\hbox{$\scriptstyle\bullet$}\kern-.4\wd0\lower.5\ht0\box0}}
\put(57.75,11.16){{\setbox0=\hbox{$\scriptstyle\bullet$}\kern-.4\wd0\lower.5\ht0\box0}}
\put(57.88,6.00){{\setbox0=\hbox{$\scriptstyle\bullet$}\kern-.4\wd0\lower.5\ht0\box0}}
\put(52.68,11.16){{\setbox0=\hbox{$\scriptstyle\bullet$}\kern-.4\wd0\lower.5\ht0\box0}}
\put(57.88,11.16){\special{em:moveto}}
\put(52.68,6.00){\special{em:lineto}}
\put(52.68,11.16){\special{em:moveto}}
\put(52.68,6.00){\special{em:lineto}}
\put(57.88,6.00){\special{em:lineto}}
\put(57.88,11.16){\special{em:lineto}}
\put(52.68,11.16){\special{em:lineto}}
\put(42.68,6.00){{\setbox0=\hbox{$\scriptstyle\bullet$}\kern-.4\wd0\lower.5\ht0\box0}}
\put(37.48,11.16){{\setbox0=\hbox{$\scriptstyle\bullet$}\kern-.4\wd0\lower.5\ht0\box0}}
\put(42.68,11.16){\special{em:moveto}}
\put(37.48,6.00){\special{em:lineto}}
\put(37.48,11.16){\special{em:moveto}}
\put(37.48,6.00){\special{em:lineto}}
\put(42.68,6.00){\special{em:lineto}}
\put(42.68,11.16){\special{em:lineto}}
\put(37.48,11.16){\special{em:lineto}}
\put(47.61,0.00){{\setbox0=\hbox{c)}\lower\ht0\box0}}
\put(4.41,4.50){{\setbox0=\hbox{$k=3,4,5,6$}\lower\ht0\box0}}
\put(8.41,8.00){{\setbox0=\hbox{$\scriptstyle\bullet$}\kern-.4\wd0\lower.5\ht0\box0}}
\put(5.75,10.16){{\setbox0=\hbox{$k$}\lower\ht0\box0}}
\put(4.01,8.00){{\setbox0=\hbox{$\scriptstyle\bullet$}\kern-.4\wd0\lower.5\ht0\box0}}
\put(12.15,11.16){{\setbox0=\hbox{$\scriptstyle\bullet$}\kern-.4\wd0\lower.5\ht0\box0}}
\put(12.15,5.83){{\setbox0=\hbox{$\scriptstyle\bullet$}\kern-.4\wd0\lower.5\ht0\box0}}
\put(4.01,8.00){\special{em:moveto}}
\put(8.15,8.00){\special{em:lineto}}
\put(12.15,11.16){\special{em:lineto}}
\put(12.15,5.83){\special{em:lineto}}
\put(8.15,8.00){\special{em:lineto}}
\end{picture}

%% file: pic/3_planes.tex
%TeX representation of Charton picture (version 1.96)
%Copyright (c) S.Trifonov, D.Tejblum, 1993-4
%FILE: C:\XBOCT\SIMPLEX\PIC\3_PLANES.tex
%DATE: Thu May 21 14:08:47 1998
%------------
\setlength{\unitlength}{0.075in}
\begin{picture}(11.56,11.35)
\special{em:linewidth 0.014in}
\put(5.08,6.40){\special{em:moveto}}
\put(4.88,6.15){\special{em:lineto}}
\put(4.68,5.90){\special{em:moveto}}
\put(4.48,5.65){\special{em:lineto}}
\put(4.28,5.40){\special{em:moveto}}
\put(4.15,5.23){\special{em:lineto}}
\put(10.61,0.93){\special{em:moveto}}
\put(10.36,0.73){\special{em:lineto}}
\put(10.11,0.53){\special{em:moveto}}
\put(9.86,0.33){\special{em:lineto}}
\put(9.60,0.13){\special{em:moveto}}
\put(9.46,0.01){\special{em:lineto}}
\put(9.88,1.71){\special{em:moveto}}
\put(9.63,1.50){\special{em:lineto}}
\put(9.38,1.28){\special{em:moveto}}
\put(9.13,1.06){\special{em:lineto}}
\put(8.88,0.85){\special{em:moveto}}
\put(8.63,0.63){\special{em:lineto}}
\put(8.38,0.41){\special{em:moveto}}
\put(8.13,0.20){\special{em:lineto}}
\put(9.16,2.50){\special{em:moveto}}
\put(8.93,2.28){\special{em:lineto}}
\put(8.70,2.06){\special{em:moveto}}
\put(8.46,1.85){\special{em:lineto}}
\put(8.23,1.63){\special{em:moveto}}
\put(8.00,1.41){\special{em:lineto}}
\put(7.76,1.18){\special{em:moveto}}
\put(7.53,0.96){\special{em:lineto}}
\put(7.30,0.75){\special{em:moveto}}
\put(7.06,0.51){\special{em:lineto}}
\put(6.83,0.30){\special{em:moveto}}
\put(6.60,0.06){\special{em:lineto}}
\put(8.33,3.28){\special{em:moveto}}
\put(8.10,3.05){\special{em:lineto}}
\put(7.86,2.81){\special{em:moveto}}
\put(7.63,2.60){\special{em:lineto}}
\put(7.40,2.38){\special{em:moveto}}
\put(7.16,2.15){\special{em:lineto}}
\put(6.93,1.93){\special{em:moveto}}
\put(6.70,1.70){\special{em:lineto}}
\put(6.46,1.46){\special{em:moveto}}
\put(6.23,1.23){\special{em:lineto}}
\put(6.00,1.00){\special{em:moveto}}
\put(5.76,0.76){\special{em:lineto}}
\put(5.53,0.55){\special{em:moveto}}
\put(5.30,0.31){\special{em:lineto}}
\put(5.06,0.10){\special{em:moveto}}
\put(4.98,0.01){\special{em:lineto}}
\put(4.15,0.80){\special{em:moveto}}
\put(4.38,1.03){\special{em:lineto}}
\put(4.61,1.26){\special{em:moveto}}
\put(4.85,1.48){\special{em:lineto}}
\put(5.08,1.70){\special{em:moveto}}
\put(5.31,1.93){\special{em:lineto}}
\put(5.55,2.15){\special{em:moveto}}
\put(5.78,2.38){\special{em:lineto}}
\put(6.01,2.61){\special{em:moveto}}
\put(6.25,2.85){\special{em:lineto}}
\put(6.48,3.08){\special{em:moveto}}
\put(6.71,3.31){\special{em:lineto}}
\put(6.95,3.53){\special{em:moveto}}
\put(7.18,3.76){\special{em:lineto}}
\put(7.41,3.98){\special{em:moveto}}
\put(7.50,4.06){\special{em:lineto}}
\put(4.15,2.36){\special{em:moveto}}
\put(4.38,2.60){\special{em:lineto}}
\put(4.61,2.83){\special{em:moveto}}
\put(4.85,3.06){\special{em:lineto}}
\put(5.08,3.30){\special{em:moveto}}
\put(5.31,3.53){\special{em:lineto}}
\put(5.55,3.76){\special{em:moveto}}
\put(5.78,4.00){\special{em:lineto}}
\put(6.01,4.23){\special{em:moveto}}
\put(6.25,4.45){\special{em:lineto}}
\put(6.48,4.68){\special{em:moveto}}
\put(6.65,4.85){\special{em:lineto}}
\put(5.93,5.63){\special{em:moveto}}
\put(5.71,5.40){\special{em:lineto}}
\put(5.48,5.16){\special{em:moveto}}
\put(5.25,4.93){\special{em:lineto}}
\put(5.01,4.70){\special{em:moveto}}
\put(4.78,4.46){\special{em:lineto}}
\put(4.56,4.23){\special{em:moveto}}
\put(4.35,4.00){\special{em:lineto}}
\put(0.50,10.86){\special{em:moveto}}
\put(0.40,10.75){\special{em:lineto}}
\put(0.30,10.63){\special{em:moveto}}
\put(0.20,10.51){\special{em:lineto}}
\put(0.10,10.38){\special{em:moveto}}
\put(0.00,10.25){\special{em:lineto}}
\put(3.40,8.00){\special{em:moveto}}
\put(3.28,7.90){\special{em:lineto}}
\put(3.15,7.80){\special{em:moveto}}
\put(3.01,7.70){\special{em:lineto}}
\put(2.88,7.60){\special{em:moveto}}
\put(2.78,7.51){\special{em:lineto}}
\put(3.01,8.41){\special{em:moveto}}
\put(2.90,8.31){\special{em:lineto}}
\put(2.78,8.21){\special{em:moveto}}
\put(2.66,8.11){\special{em:lineto}}
\put(2.55,8.01){\special{em:moveto}}
\put(2.43,7.91){\special{em:lineto}}
\put(2.31,7.81){\special{em:moveto}}
\put(2.20,7.71){\special{em:lineto}}
\put(2.08,7.61){\special{em:moveto}}
\put(1.96,7.51){\special{em:lineto}}
\put(2.63,8.81){\special{em:moveto}}
\put(2.51,8.71){\special{em:lineto}}
\put(2.40,8.61){\special{em:moveto}}
\put(2.28,8.50){\special{em:lineto}}
\put(2.16,8.38){\special{em:moveto}}
\put(2.05,8.26){\special{em:lineto}}
\put(1.93,8.16){\special{em:moveto}}
\put(1.81,8.05){\special{em:lineto}}
\put(1.70,7.93){\special{em:moveto}}
\put(1.58,7.83){\special{em:lineto}}
\put(1.46,7.71){\special{em:moveto}}
\put(1.35,7.60){\special{em:lineto}}
\put(2.18,9.23){\special{em:moveto}}
\put(2.06,9.11){\special{em:lineto}}
\put(1.95,9.00){\special{em:moveto}}
\put(1.83,8.88){\special{em:lineto}}
\put(1.71,8.76){\special{em:moveto}}
\put(1.60,8.65){\special{em:lineto}}
\put(1.48,8.53){\special{em:moveto}}
\put(1.36,8.41){\special{em:lineto}}
\put(1.25,8.30){\special{em:moveto}}
\put(1.13,8.18){\special{em:lineto}}
\put(1.01,8.08){\special{em:moveto}}
\put(0.90,7.96){\special{em:lineto}}
\put(0.78,7.85){\special{em:moveto}}
\put(0.66,7.75){\special{em:lineto}}
\put(0.55,7.63){\special{em:moveto}}
\put(0.43,7.51){\special{em:lineto}}
\put(0.00,7.93){\special{em:moveto}}
\put(0.11,8.03){\special{em:lineto}}
\put(0.23,8.15){\special{em:moveto}}
\put(0.35,8.26){\special{em:lineto}}
\put(0.46,8.38){\special{em:moveto}}
\put(0.58,8.50){\special{em:lineto}}
\put(0.70,8.61){\special{em:moveto}}
\put(0.81,8.73){\special{em:lineto}}
\put(0.93,8.85){\special{em:moveto}}
\put(1.05,8.96){\special{em:lineto}}
\put(1.16,9.06){\special{em:moveto}}
\put(1.28,9.18){\special{em:lineto}}
\put(1.40,9.30){\special{em:moveto}}
\put(1.51,9.40){\special{em:lineto}}
\put(1.63,9.51){\special{em:moveto}}
\put(1.75,9.63){\special{em:lineto}}
\put(0.00,8.75){\special{em:moveto}}
\put(0.11,8.86){\special{em:lineto}}
\put(0.23,8.98){\special{em:moveto}}
\put(0.35,9.10){\special{em:lineto}}
\put(0.46,9.21){\special{em:moveto}}
\put(0.58,9.33){\special{em:lineto}}
\put(0.70,9.45){\special{em:moveto}}
\put(0.81,9.56){\special{em:lineto}}
\put(0.93,9.68){\special{em:moveto}}
\put(1.05,9.80){\special{em:lineto}}
\put(1.16,9.91){\special{em:moveto}}
\put(1.28,10.01){\special{em:lineto}}
\put(0.93,10.46){\special{em:moveto}}
\put(0.81,10.35){\special{em:lineto}}
\put(0.70,10.23){\special{em:moveto}}
\put(0.58,10.11){\special{em:lineto}}
\put(0.46,10.00){\special{em:moveto}}
\put(0.36,9.88){\special{em:lineto}}
\put(0.25,9.76){\special{em:moveto}}
\put(0.13,9.65){\special{em:lineto}}
\put(0.03,9.53){\special{em:moveto}}
\put(0.00,9.50){\special{em:lineto}}
\put(0.00,7.18){\special{em:moveto}}
\put(3.61,7.18){\special{em:lineto}}
\put(3.61,0.00){\special{em:lineto}}
\put(3.96,0.00){\special{em:moveto}}
\put(3.96,7.55){\special{em:lineto}}
\put(7.60,3.76){\special{em:moveto}}
\put(3.96,0.00){\special{em:lineto}}
\put(0.00,3.76){\special{em:lineto}}
\put(7.60,0.00){\special{em:moveto}}
\put(7.60,3.76){\special{em:lineto}}
\put(0.00,3.76){\special{em:lineto}}
\put(3.96,7.55){\special{em:lineto}}
\put(0.00,7.55){\special{em:lineto}}
\put(0.00,11.35){\special{em:moveto}}
\put(0.00,0.00){\special{em:lineto}}
\put(11.56,0.00){\special{em:lineto}}
\put(0.00,11.35){\special{em:lineto}}
\end{picture}

%% file: pic/vol_3_og.tex
%TeX representation of Charton picture (version 1.96)
%Copyright (c) S.Trifonov, D.Tejblum, 1993-4
%FILE: VOL_3_OG.tex
%DATE: Thu Nov 08 21:21:07 2001
%------------
\setlength{\unitlength}{0.09in}
\begin{picture}(54.26,60.83)
\put(19.33,57.66){{\setbox0=\hbox{Notation}\lower\ht0\box0}}
\special{em:linewidth 0.014in}
\put(0.00,54.16){\special{em:moveto}}
\put(54.26,54.16){\special{em:lineto}}
\put(54.26,60.83){\special{em:moveto}}
\put(0.00,60.83){\special{em:lineto}}
\put(0.00,60.83){\special{em:moveto}}
\put(0.00,0.00){\special{em:lineto}}
\put(20.26,6.83){{\setbox0=\hbox{$H_9$}\lower\ht0\box0}}
\put(20.26,13.50){{\setbox0=\hbox{$H_8$}\lower\ht0\box0}}
\put(20.26,20.16){{\setbox0=\hbox{$H_7$}\lower\ht0\box0}}
\put(20.26,26.83){{\setbox0=\hbox{$H_6$}\lower\ht0\box0}}
\put(20.26,31.83){{\setbox0=\hbox{$H_5$}\lower\ht0\box0}}
\put(20.26,36.83){{\setbox0=\hbox{$H_4$}\lower\ht0\box0}}
\put(20.26,41.83){{\setbox0=\hbox{$H_3$}\lower\ht0\box0}}
\put(20.26,46.83){{\setbox0=\hbox{$H_2$}\lower\ht0\box0}}
\put(10.40,2.83){{\setbox0=\hbox{5}\lower\ht0\box0}}
\put(10.66,9.00){{\setbox0=\hbox{5}\lower\ht0\box0}}
\put(3.60,15.33){{\setbox0=\hbox{4}\lower\ht0\box0}}
\put(3.60,9.33){{\setbox0=\hbox{5}\lower\ht0\box0}}
\put(10.26,16.66){{\setbox0=\hbox{4}\lower\ht0\box0}}
\put(10.26,22.83){{\setbox0=\hbox{4}\lower\ht0\box0}}
\put(6.40,26.00){{\setbox0=\hbox{5}\lower\ht0\box0}}
\put(10.26,31.83){{\setbox0=\hbox{5}\lower\ht0\box0}}
\put(3.60,31.83){{\setbox0=\hbox{5}\lower\ht0\box0}}
\put(6.40,35.83){{\setbox0=\hbox{4}\lower\ht0\box0}}
\put(3.60,41.83){{\setbox0=\hbox{5}\lower\ht0\box0}}
\put(6.66,46.83){{\setbox0=\hbox{5}\lower\ht0\box0}}
\put(10.26,52.00){{\setbox0=\hbox{5}\lower\ht0\box0}}
\put(3.46,52.00){{\setbox0=\hbox{4}\lower\ht0\box0}}
\put(0.00,0.00){\special{em:moveto}}
\put(54.26,0.00){\special{em:lineto}}
\put(37.20,57.66){{\setbox0=\hbox{Volume}\lower\ht0\box0}}
\put(2.53,57.33){{\setbox0=\hbox{diagram}\lower\ht0\box0}}
\put(2.26,36.83){\special{em:moveto}}
\put(5.60,36.83){\special{em:lineto}}
\put(2.26,33.50){\special{em:moveto}}
\put(5.60,33.50){\special{em:lineto}}
\put(5.60,36.83){\special{em:moveto}}
\put(5.60,33.50){\special{em:lineto}}
\put(2.26,36.83){\special{em:moveto}}
\put(2.26,33.50){\special{em:lineto}}
\put(2.26,33.50){{\setbox0=\hbox{$\scriptstyle\bullet$}\kern-.4\wd0\lower.5\ht0\box0}}
\put(2.26,36.83){{\setbox0=\hbox{$\scriptstyle\bullet$}\kern-.4\wd0\lower.5\ht0\box0}}
\put(5.60,33.50){{\setbox0=\hbox{$\scriptstyle\bullet$}\kern-.4\wd0\lower.5\ht0\box0}}
\put(5.60,36.83){{\setbox0=\hbox{$\scriptstyle\bullet$}\kern-.4\wd0\lower.5\ht0\box0}}
\put(54.26,60.83){\special{em:moveto}}
\put(54.26,0.00){\special{em:lineto}}
\put(9.33,7.16){\special{em:moveto}}
\put(12.80,7.16){\special{em:lineto}}
\put(36.40,5.50){{\setbox0=\hbox{0.5021308905}\lower\ht0\box0}}
\put(36.40,12.16){{\setbox0=\hbox{0.3586534401}\lower\ht0\box0}}
\put(36.40,18.83){{\setbox0=\hbox{0.2222287320}\lower\ht0\box0}}
\put(36.40,25.50){{\setbox0=\hbox{0.2052887885}\lower\ht0\box0}}
\put(36.40,30.50){{\setbox0=\hbox{0.0933255395}\lower\ht0\box0}}
\put(36.40,35.50){{\setbox0=\hbox{0.0857701820}\lower\ht0\box0}}
\put(36.40,40.50){{\setbox0=\hbox{0.0717701267}\lower\ht0\box0}}
\put(36.40,45.50){{\setbox0=\hbox{0.0390502856}\lower\ht0\box0}}
\put(36.40,50.50){{\setbox0=\hbox{0.0358850633}\lower\ht0\box0}}
\put(9.33,3.83){\special{em:moveto}}
\put(12.80,3.83){\special{em:lineto}}
\put(12.80,7.16){\special{em:moveto}}
\put(12.80,3.83){\special{em:lineto}}
\put(9.46,7.16){\special{em:moveto}}
\put(9.46,3.83){\special{em:lineto}}
\put(9.46,3.83){{\setbox0=\hbox{$\scriptstyle\bullet$}\kern-.4\wd0\lower.5\ht0\box0}}
\put(9.46,7.16){{\setbox0=\hbox{$\scriptstyle\bullet$}\kern-.4\wd0\lower.5\ht0\box0}}
\put(12.80,3.83){{\setbox0=\hbox{$\scriptstyle\bullet$}\kern-.4\wd0\lower.5\ht0\box0}}
\put(12.80,7.16){{\setbox0=\hbox{$\scriptstyle\bullet$}\kern-.4\wd0\lower.5\ht0\box0}}
\put(2.13,10.16){\special{em:moveto}}
\put(5.60,10.16){\special{em:lineto}}
\put(2.26,13.50){\special{em:moveto}}
\put(5.60,13.50){\special{em:lineto}}
\put(5.60,13.50){\special{em:moveto}}
\put(5.60,10.16){\special{em:lineto}}
\put(2.26,13.50){\special{em:moveto}}
\put(2.26,10.16){\special{em:lineto}}
\put(2.26,10.16){{\setbox0=\hbox{$\scriptstyle\bullet$}\kern-.4\wd0\lower.5\ht0\box0}}
\put(2.26,13.50){{\setbox0=\hbox{$\scriptstyle\bullet$}\kern-.4\wd0\lower.5\ht0\box0}}
\put(5.60,10.16){{\setbox0=\hbox{$\scriptstyle\bullet$}\kern-.4\wd0\lower.5\ht0\box0}}
\put(5.60,13.50){{\setbox0=\hbox{$\scriptstyle\bullet$}\kern-.4\wd0\lower.5\ht0\box0}}
\put(9.06,21.00){\special{em:moveto}}
\put(12.40,21.00){\special{em:lineto}}
\put(9.06,17.66){\special{em:moveto}}
\put(12.40,17.66){\special{em:lineto}}
\put(12.40,21.00){\special{em:moveto}}
\put(12.40,17.66){\special{em:lineto}}
\put(9.06,21.00){\special{em:moveto}}
\put(9.06,17.66){\special{em:lineto}}
\put(9.06,17.66){{\setbox0=\hbox{$\scriptstyle\bullet$}\kern-.4\wd0\lower.5\ht0\box0}}
\put(9.06,21.00){{\setbox0=\hbox{$\scriptstyle\bullet$}\kern-.4\wd0\lower.5\ht0\box0}}
\put(12.40,17.66){{\setbox0=\hbox{$\scriptstyle\bullet$}\kern-.4\wd0\lower.5\ht0\box0}}
\put(12.40,21.00){{\setbox0=\hbox{$\scriptstyle\bullet$}\kern-.4\wd0\lower.5\ht0\box0}}
\put(2.26,23.50){{\setbox0=\hbox{$\scriptstyle\bullet$}\kern-.4\wd0\lower.5\ht0\box0}}
\put(2.26,26.83){{\setbox0=\hbox{$\scriptstyle\bullet$}\kern-.4\wd0\lower.5\ht0\box0}}
\put(5.60,23.50){{\setbox0=\hbox{$\scriptstyle\bullet$}\kern-.4\wd0\lower.5\ht0\box0}}
\put(5.60,26.83){{\setbox0=\hbox{$\scriptstyle\bullet$}\kern-.4\wd0\lower.5\ht0\box0}}
\put(2.26,26.83){\special{em:moveto}}
\put(5.60,26.83){\special{em:lineto}}
\put(5.60,23.50){\special{em:lineto}}
\put(2.26,23.50){\special{em:lineto}}
\put(2.26,26.83){\special{em:lineto}}
\put(8.93,30.16){\special{em:moveto}}
\put(12.40,30.16){\special{em:lineto}}
\put(12.40,30.16){{\setbox0=\hbox{$\scriptstyle\bullet$}\kern-.4\wd0\lower.5\ht0\box0}}
\put(5.06,30.16){\special{em:moveto}}
\put(8.93,30.16){\special{em:lineto}}
\put(8.93,30.16){{\setbox0=\hbox{$\scriptstyle\bullet$}\kern-.4\wd0\lower.5\ht0\box0}}
\put(2.26,30.16){\special{em:moveto}}
\put(5.73,30.16){\special{em:lineto}}
\put(5.73,30.16){{\setbox0=\hbox{$\scriptstyle\bullet$}\kern-.4\wd0\lower.5\ht0\box0}}
\put(2.26,30.16){{\setbox0=\hbox{$\scriptstyle\bullet$}\kern-.4\wd0\lower.5\ht0\box0}}
\put(8.93,38.50){{\setbox0=\hbox{$\scriptstyle\bullet$}\kern-.4\wd0\lower.5\ht0\box0}}
\put(8.93,41.83){{\setbox0=\hbox{$\scriptstyle\bullet$}\kern-.4\wd0\lower.5\ht0\box0}}
\put(8.93,41.83){\special{em:moveto}}
\put(5.60,40.16){\special{em:lineto}}
\put(8.93,38.50){\special{em:lineto}}
\put(2.26,40.16){\special{em:moveto}}
\put(5.73,40.16){\special{em:lineto}}
\put(5.73,40.16){{\setbox0=\hbox{$\scriptstyle\bullet$}\kern-.4\wd0\lower.5\ht0\box0}}
\put(2.26,40.16){{\setbox0=\hbox{$\scriptstyle\bullet$}\kern-.4\wd0\lower.5\ht0\box0}}
\put(8.40,45.16){\special{em:moveto}}
\put(12.26,45.16){\special{em:lineto}}
\put(12.26,45.16){{\setbox0=\hbox{$\scriptstyle\bullet$}\kern-.4\wd0\lower.5\ht0\box0}}
\put(5.60,45.16){\special{em:moveto}}
\put(9.06,45.16){\special{em:lineto}}
\put(9.06,45.16){{\setbox0=\hbox{$\scriptstyle\bullet$}\kern-.4\wd0\lower.5\ht0\box0}}
\put(5.60,45.16){{\setbox0=\hbox{$\scriptstyle\bullet$}\kern-.4\wd0\lower.5\ht0\box0}}
\put(2.40,45.16){{\setbox0=\hbox{$\scriptstyle\bullet$}\kern-.4\wd0\lower.5\ht0\box0}}
\put(2.40,45.16){\special{em:moveto}}
\put(5.60,45.16){\special{em:lineto}}
\put(8.80,50.16){\special{em:moveto}}
\put(12.26,50.16){\special{em:lineto}}
\put(2.53,50.16){\special{em:moveto}}
\put(5.60,50.16){\special{em:lineto}}
\put(12.26,50.16){{\setbox0=\hbox{$\scriptstyle\bullet$}\kern-.4\wd0\lower.5\ht0\box0}}
\put(2.53,50.16){{\setbox0=\hbox{$\scriptstyle\bullet$}\kern-.4\wd0\lower.5\ht0\box0}}
\put(8.80,50.16){{\setbox0=\hbox{$\scriptstyle\bullet$}\kern-.4\wd0\lower.5\ht0\box0}}
\put(5.60,50.16){{\setbox0=\hbox{$\scriptstyle\bullet$}\kern-.4\wd0\lower.5\ht0\box0}}
\put(5.60,50.16){\special{em:moveto}}
\put(8.80,50.16){\special{em:lineto}}
\put(20.26,51.83){{\setbox0=\hbox{$H_1$}\lower\ht0\box0}}
\put(14.93,60.83){\special{em:moveto}}
\put(15.06,0.00){\special{em:lineto}}
\put(32.53,60.83){\special{em:moveto}}
\put(32.53,0.00){\special{em:lineto}}
\put(2.53,59.50){{\setbox0=\hbox{Coxeter}\lower\ht0\box0}}
\end{picture}

%% file: pic/vol_3ne1.tex
%TeX representation of Charton picture (version 1.96)
%Copyright (c) S.Trifonov, D.Tejblum, 1993-4
%FILE: VOL_3NE1.tex
%DATE: Thu Nov 08 21:57:03 2001
%------------
\setlength{\unitlength}{0.1in}
\begin{picture}(51.20,70.16)
\put(22.53,67.00){{\setbox0=\hbox{Notation}\lower\ht0\box0}}
\special{em:linewidth 0.014in}
\put(0.13,63.50){\special{em:moveto}}
\put(51.20,63.50){\special{em:lineto}}
\put(37.80,67.00){{\setbox0=\hbox{Volume}\lower\ht0\box0}}
\put(4.80,65.83){{\setbox0=\hbox{diagram}\lower\ht0\box0}}
\put(4.80,68.33){{\setbox0=\hbox{Coxeter}\lower\ht0\box0}}
\put(51.20,70.16){\special{em:moveto}}
\put(0.13,70.16){\special{em:lineto}}
\put(10.26,6.66){{\setbox0=\hbox{5}\lower\ht0\box0}}
\put(4.26,8.16){\special{em:moveto}}
\put(7.60,8.16){\special{em:lineto}}
\put(7.60,11.50){\special{em:lineto}}
\put(5.46,13.16){{\setbox0=\hbox{4}\lower\ht0\box0}}
\put(8.40,10.33){{\setbox0=\hbox{4}\lower\ht0\box0}}
\put(4.26,18.16){\special{em:moveto}}
\put(14.40,18.16){\special{em:lineto}}
\put(12.26,20.00){{\setbox0=\hbox{4}\lower\ht0\box0}}
\put(8.80,20.00){{\setbox0=\hbox{4}\lower\ht0\box0}}
\put(5.20,20.00){{\setbox0=\hbox{4}\lower\ht0\box0}}
\put(5.06,30.00){{\setbox0=\hbox{4}\lower\ht0\box0}}
\put(12.26,33.50){{\setbox0=\hbox{5}\lower\ht0\box0}}
\put(8.13,39.00){{\setbox0=\hbox{4}\lower\ht0\box0}}
\put(8.13,42.33){{\setbox0=\hbox{4}\lower\ht0\box0}}
\put(12.26,46.66){{\setbox0=\hbox{4}\lower\ht0\box0}}
\put(8.93,56.66){{\setbox0=\hbox{4}\lower\ht0\box0}}
\put(5.60,56.66){{\setbox0=\hbox{4}\lower\ht0\box0}}
\put(0.00,0.00){\special{em:moveto}}
\put(51.20,0.00){\special{em:lineto}}
\put(35.20,70.16){\special{em:moveto}}
\put(35.20,0.00){\special{em:lineto}}
\put(51.20,70.16){\special{em:moveto}}
\put(51.20,0.00){\special{em:lineto}}
\put(37.33,5.50){{\setbox0=\hbox{0.3430033226}\lower\ht0\box0}}
\put(37.20,10.33){{\setbox0=\hbox{0.3053218647}\lower\ht0\box0}}
\put(25.06,6.16){{\setbox0=\hbox{$H_{22}$}\lower\ht0\box0}}
\put(9.16,5.00){\special{em:moveto}}
\put(12.63,5.00){\special{em:lineto}}
\put(9.16,5.00){{\setbox0=\hbox{$\scriptstyle\bullet$}\kern-.4\wd0\lower.5\ht0\box0}}
\put(15.83,3.33){{\setbox0=\hbox{$\scriptstyle\bullet$}\kern-.4\wd0\lower.5\ht0\box0}}
\put(15.83,6.66){{\setbox0=\hbox{$\scriptstyle\bullet$}\kern-.4\wd0\lower.5\ht0\box0}}
\put(12.50,5.00){{\setbox0=\hbox{$\scriptstyle\bullet$}\kern-.4\wd0\lower.5\ht0\box0}}
\put(12.50,5.00){\special{em:moveto}}
\put(15.83,6.66){\special{em:lineto}}
\put(15.83,3.33){\special{em:lineto}}
\put(12.50,5.00){\special{em:lineto}}
\put(4.26,8.16){\special{em:moveto}}
\put(4.26,11.50){\special{em:lineto}}
\put(4.26,11.50){\special{em:moveto}}
\put(7.60,11.50){\special{em:lineto}}
\put(4.26,8.16){{\setbox0=\hbox{$\scriptstyle\bullet$}\kern-.4\wd0\lower.5\ht0\box0}}
\put(4.26,11.50){{\setbox0=\hbox{$\scriptstyle\bullet$}\kern-.4\wd0\lower.5\ht0\box0}}
\put(7.60,8.16){{\setbox0=\hbox{$\scriptstyle\bullet$}\kern-.4\wd0\lower.5\ht0\box0}}
\put(7.60,11.50){{\setbox0=\hbox{$\scriptstyle\bullet$}\kern-.4\wd0\lower.5\ht0\box0}}
\put(25.06,11.00){{\setbox0=\hbox{$H_{21}$}\lower\ht0\box0}}
\put(18.53,70.16){\special{em:moveto}}
\put(18.53,0.00){\special{em:lineto}}
\put(37.20,15.33){{\setbox0=\hbox{0.2537354016}\lower\ht0\box0}}
\put(37.20,18.66){{\setbox0=\hbox{0.2289913985}\lower\ht0\box0}}
\put(37.20,23.66){{\setbox0=\hbox{0.2114461680}\lower\ht0\box0}}
\put(37.20,28.83){{\setbox0=\hbox{0.2114461680}\lower\ht0\box0}}
\put(37.20,32.16){{\setbox0=\hbox{0.1715016613}\lower\ht0\box0}}
\put(37.20,35.50){{\setbox0=\hbox{0.1691569344}\lower\ht0\box0}}
\put(37.20,40.50){{\setbox0=\hbox{0.1526609324}\lower\ht0\box0}}
\put(37.20,45.50){{\setbox0=\hbox{0.1057230840}\lower\ht0\box0}}
\put(37.20,50.50){{\setbox0=\hbox{0.0845784672}\lower\ht0\box0}}
\put(37.20,55.50){{\setbox0=\hbox{0.0763304662}\lower\ht0\box0}}
\put(37.20,60.33){{\setbox0=\hbox{0.0422892336}\lower\ht0\box0}}
\put(25.06,16.00){{\setbox0=\hbox{$H_{20}$}\lower\ht0\box0}}
\put(12.26,16.50){{\setbox0=\hbox{6}\lower\ht0\box0}}
\put(5.20,16.50){{\setbox0=\hbox{6}\lower\ht0\box0}}
\put(4.26,14.83){\special{em:moveto}}
\put(7.60,14.83){\special{em:lineto}}
\put(10.80,14.83){\special{em:moveto}}
\put(14.26,14.83){\special{em:lineto}}
\put(14.26,14.83){{\setbox0=\hbox{$\scriptstyle\bullet$}\kern-.4\wd0\lower.5\ht0\box0}}
\put(4.26,14.83){{\setbox0=\hbox{$\scriptstyle\bullet$}\kern-.4\wd0\lower.5\ht0\box0}}
\put(10.80,14.83){{\setbox0=\hbox{$\scriptstyle\bullet$}\kern-.4\wd0\lower.5\ht0\box0}}
\put(7.60,14.83){{\setbox0=\hbox{$\scriptstyle\bullet$}\kern-.4\wd0\lower.5\ht0\box0}}
\put(7.60,14.83){\special{em:moveto}}
\put(10.80,14.83){\special{em:lineto}}
\put(25.06,19.33){{\setbox0=\hbox{$H_{19}$}\lower\ht0\box0}}
\put(11.06,18.16){{\setbox0=\hbox{$\scriptstyle\bullet$}\kern-.4\wd0\lower.5\ht0\box0}}
\put(14.40,18.16){{\setbox0=\hbox{$\scriptstyle\bullet$}\kern-.4\wd0\lower.5\ht0\box0}}
\put(4.26,18.16){{\setbox0=\hbox{$\scriptstyle\bullet$}\kern-.4\wd0\lower.5\ht0\box0}}
\put(7.60,18.16){{\setbox0=\hbox{$\scriptstyle\bullet$}\kern-.4\wd0\lower.5\ht0\box0}}
\put(5.33,61.66){{\setbox0=\hbox{6}\lower\ht0\box0}}
\put(4.26,59.83){\special{em:moveto}}
\put(7.60,59.83){\special{em:lineto}}
\put(10.80,59.83){\special{em:moveto}}
\put(14.26,59.83){\special{em:lineto}}
\put(14.26,59.83){{\setbox0=\hbox{$\scriptstyle\bullet$}\kern-.4\wd0\lower.5\ht0\box0}}
\put(4.26,59.83){{\setbox0=\hbox{$\scriptstyle\bullet$}\kern-.4\wd0\lower.5\ht0\box0}}
\put(10.80,59.83){{\setbox0=\hbox{$\scriptstyle\bullet$}\kern-.4\wd0\lower.5\ht0\box0}}
\put(7.60,59.83){{\setbox0=\hbox{$\scriptstyle\bullet$}\kern-.4\wd0\lower.5\ht0\box0}}
\put(7.60,59.83){\special{em:moveto}}
\put(10.80,59.83){\special{em:lineto}}
\put(25.06,24.33){{\setbox0=\hbox{$H_{18}$}\lower\ht0\box0}}
\put(10.80,24.33){{\setbox0=\hbox{6}\lower\ht0\box0}}
\put(16.66,20.83){{\setbox0=\hbox{$\scriptstyle\bullet$}\kern-.4\wd0\lower.5\ht0\box0}}
\put(16.66,24.16){{\setbox0=\hbox{$\scriptstyle\bullet$}\kern-.4\wd0\lower.5\ht0\box0}}
\put(16.66,24.16){\special{em:moveto}}
\put(13.33,22.50){\special{em:lineto}}
\put(16.66,20.83){\special{em:lineto}}
\put(10.00,22.50){\special{em:moveto}}
\put(13.46,22.50){\special{em:lineto}}
\put(13.46,22.50){{\setbox0=\hbox{$\scriptstyle\bullet$}\kern-.4\wd0\lower.5\ht0\box0}}
\put(10.00,22.50){{\setbox0=\hbox{$\scriptstyle\bullet$}\kern-.4\wd0\lower.5\ht0\box0}}
\put(25.06,29.50){{\setbox0=\hbox{$H_{17}$}\lower\ht0\box0}}
\put(25.06,32.83){{\setbox0=\hbox{$H_{16}$}\lower\ht0\box0}}
\put(25.06,36.16){{\setbox0=\hbox{$H_{15}$}\lower\ht0\box0}}
\put(9.20,36.66){{\setbox0=\hbox{6}\lower\ht0\box0}}
\put(4.26,35.00){\special{em:moveto}}
\put(7.60,35.00){\special{em:lineto}}
\put(10.80,35.00){\special{em:moveto}}
\put(14.26,35.00){\special{em:lineto}}
\put(14.26,35.00){{\setbox0=\hbox{$\scriptstyle\bullet$}\kern-.4\wd0\lower.5\ht0\box0}}
\put(4.26,35.00){{\setbox0=\hbox{$\scriptstyle\bullet$}\kern-.4\wd0\lower.5\ht0\box0}}
\put(10.80,35.00){{\setbox0=\hbox{$\scriptstyle\bullet$}\kern-.4\wd0\lower.5\ht0\box0}}
\put(7.60,35.00){{\setbox0=\hbox{$\scriptstyle\bullet$}\kern-.4\wd0\lower.5\ht0\box0}}
\put(7.60,35.00){\special{em:moveto}}
\put(10.80,35.00){\special{em:lineto}}
\put(7.60,40.00){\special{em:moveto}}
\put(10.66,38.66){\special{em:lineto}}
\put(7.56,40.00){\special{em:moveto}}
\put(10.80,41.33){\special{em:lineto}}
\put(4.26,40.00){\special{em:moveto}}
\put(7.60,40.00){\special{em:lineto}}
\put(4.26,40.00){{\setbox0=\hbox{$\scriptstyle\bullet$}\kern-.4\wd0\lower.5\ht0\box0}}
\put(10.66,38.66){{\setbox0=\hbox{$\scriptstyle\bullet$}\kern-.4\wd0\lower.5\ht0\box0}}
\put(10.80,41.33){{\setbox0=\hbox{$\scriptstyle\bullet$}\kern-.4\wd0\lower.5\ht0\box0}}
\put(7.60,40.00){{\setbox0=\hbox{$\scriptstyle\bullet$}\kern-.4\wd0\lower.5\ht0\box0}}
\put(4.26,28.16){\special{em:moveto}}
\put(7.60,28.16){\special{em:lineto}}
\put(4.26,28.16){{\setbox0=\hbox{$\scriptstyle\bullet$}\kern-.4\wd0\lower.5\ht0\box0}}
\put(10.93,26.50){{\setbox0=\hbox{$\scriptstyle\bullet$}\kern-.4\wd0\lower.5\ht0\box0}}
\put(10.93,29.83){{\setbox0=\hbox{$\scriptstyle\bullet$}\kern-.4\wd0\lower.5\ht0\box0}}
\put(7.60,28.16){{\setbox0=\hbox{$\scriptstyle\bullet$}\kern-.4\wd0\lower.5\ht0\box0}}
\put(7.60,28.16){\special{em:moveto}}
\put(10.93,29.83){\special{em:lineto}}
\put(10.93,26.50){\special{em:lineto}}
\put(7.60,28.16){\special{em:lineto}}
\put(25.06,41.16){{\setbox0=\hbox{$H_{14}$}\lower\ht0\box0}}
\put(25.06,46.16){{\setbox0=\hbox{$H_{13}$}\lower\ht0\box0}}
\put(10.80,45.00){\special{em:moveto}}
\put(14.26,45.00){\special{em:lineto}}
\put(5.60,46.66){{\setbox0=\hbox{6}\lower\ht0\box0}}
\put(4.26,45.00){\special{em:moveto}}
\put(7.60,45.00){\special{em:lineto}}
\put(14.26,45.00){{\setbox0=\hbox{$\scriptstyle\bullet$}\kern-.4\wd0\lower.5\ht0\box0}}
\put(4.26,45.00){{\setbox0=\hbox{$\scriptstyle\bullet$}\kern-.4\wd0\lower.5\ht0\box0}}
\put(10.80,45.00){{\setbox0=\hbox{$\scriptstyle\bullet$}\kern-.4\wd0\lower.5\ht0\box0}}
\put(7.60,45.00){{\setbox0=\hbox{$\scriptstyle\bullet$}\kern-.4\wd0\lower.5\ht0\box0}}
\put(7.60,45.00){\special{em:moveto}}
\put(10.80,45.00){\special{em:lineto}}
\put(4.26,50.00){{\setbox0=\hbox{$\scriptstyle\bullet$}\kern-.4\wd0\lower.5\ht0\box0}}
\put(25.06,51.16){{\setbox0=\hbox{$H_{12}$}\lower\ht0\box0}}
\put(10.93,48.33){{\setbox0=\hbox{$\scriptstyle\bullet$}\kern-.4\wd0\lower.5\ht0\box0}}
\put(10.93,51.66){{\setbox0=\hbox{$\scriptstyle\bullet$}\kern-.4\wd0\lower.5\ht0\box0}}
\put(7.60,50.00){{\setbox0=\hbox{$\scriptstyle\bullet$}\kern-.4\wd0\lower.5\ht0\box0}}
\put(7.60,50.00){\special{em:moveto}}
\put(10.93,51.66){\special{em:lineto}}
\put(10.93,48.33){\special{em:lineto}}
\put(7.60,50.00){\special{em:lineto}}
\put(4.26,50.00){\special{em:moveto}}
\put(7.60,50.00){\special{em:lineto}}
\put(25.06,56.16){{\setbox0=\hbox{$H_{11}$}\lower\ht0\box0}}
\put(7.46,55.00){\special{em:moveto}}
\put(10.93,55.00){\special{em:lineto}}
\put(10.93,55.00){\special{em:moveto}}
\put(14.26,55.00){\special{em:lineto}}
\put(10.93,55.00){{\setbox0=\hbox{$\scriptstyle\bullet$}\kern-.4\wd0\lower.5\ht0\box0}}
\put(4.13,55.00){\special{em:moveto}}
\put(7.46,55.00){\special{em:lineto}}
\put(14.26,55.00){{\setbox0=\hbox{$\scriptstyle\bullet$}\kern-.4\wd0\lower.5\ht0\box0}}
\put(4.13,55.00){{\setbox0=\hbox{$\scriptstyle\bullet$}\kern-.4\wd0\lower.5\ht0\box0}}
\put(7.46,55.00){{\setbox0=\hbox{$\scriptstyle\bullet$}\kern-.4\wd0\lower.5\ht0\box0}}
\put(5.60,33.50){{\setbox0=\hbox{6}\lower\ht0\box0}}
\put(4.26,31.66){\special{em:moveto}}
\put(7.60,31.66){\special{em:lineto}}
\put(14.26,31.66){{\setbox0=\hbox{$\scriptstyle\bullet$}\kern-.4\wd0\lower.5\ht0\box0}}
\put(4.26,31.66){{\setbox0=\hbox{$\scriptstyle\bullet$}\kern-.4\wd0\lower.5\ht0\box0}}
\put(10.80,31.66){{\setbox0=\hbox{$\scriptstyle\bullet$}\kern-.4\wd0\lower.5\ht0\box0}}
\put(7.60,31.66){{\setbox0=\hbox{$\scriptstyle\bullet$}\kern-.4\wd0\lower.5\ht0\box0}}
\put(7.60,31.66){\special{em:moveto}}
\put(10.80,31.66){\special{em:lineto}}
\put(25.06,60.83){{\setbox0=\hbox{$H_{10}$}\lower\ht0\box0}}
\put(10.80,31.66){\special{em:moveto}}
\put(14.26,31.66){\special{em:lineto}}
\put(0.00,70.16){\special{em:moveto}}
\put(0.00,0.00){\special{em:lineto}}
\end{picture}

%% file: pic/vol_3ne2.tex
%TeX representation of Charton picture (version 1.96)
%Copyright (c) S.Trifonov, D.Tejblum, 1993-4
%FILE: VOL_3NE2.tex
%DATE: Thu Nov 08 21:56:20 2001
%------------
\setlength{\unitlength}{0.1in}
\begin{picture}(51.86,68.50)
\put(21.13,65.33){{\setbox0=\hbox{Notation}\lower\ht0\box0}}
\special{em:linewidth 0.014in}
\put(0.13,61.83){\special{em:moveto}}
\put(51.86,61.83){\special{em:lineto}}
\put(0.00,68.50){\special{em:moveto}}
\put(51.73,68.50){\special{em:lineto}}
\put(36.93,65.33){{\setbox0=\hbox{Volume}\lower\ht0\box0}}
\put(4.00,64.16){{\setbox0=\hbox{diagram}\lower\ht0\box0}}
\put(3.73,66.66){{\setbox0=\hbox{Coxeter}\lower\ht0\box0}}
\put(51.73,68.33){\special{em:moveto}}
\put(51.86,0.00){\special{em:lineto}}
\put(2.26,11.83){\special{em:moveto}}
\put(5.60,11.83){\special{em:lineto}}
\put(5.60,8.46){\special{em:lineto}}
\put(2.26,8.46){\special{em:lineto}}
\put(2.26,11.83){\special{em:lineto}}
\put(3.46,13.66){{\setbox0=\hbox{4}\lower\ht0\box0}}
\put(3.33,8.00){{\setbox0=\hbox{4}\lower\ht0\box0}}
\put(6.26,10.83){{\setbox0=\hbox{4}\lower\ht0\box0}}
\put(0.93,10.66){{\setbox0=\hbox{4}\lower\ht0\box0}}
\put(3.60,20.66){{\setbox0=\hbox{5}\lower\ht0\box0}}
\put(13.86,31.83){\special{em:moveto}}
\put(10.53,31.83){\special{em:lineto}}
\put(10.53,28.50){\special{em:lineto}}
\put(13.86,28.50){\special{em:lineto}}
\put(11.33,33.66){{\setbox0=\hbox{4}\lower\ht0\box0}}
\put(11.60,27.83){{\setbox0=\hbox{4}\lower\ht0\box0}}
\put(9.20,30.66){{\setbox0=\hbox{4}\lower\ht0\box0}}
\put(3.60,34.33){{\setbox0=\hbox{4}\lower\ht0\box0}}
\put(6.53,45.50){{\setbox0=\hbox{4}\lower\ht0\box0}}
\put(6.80,49.83){{\setbox0=\hbox{4}\lower\ht0\box0}}
\put(3.33,48.83){{\setbox0=\hbox{4}\lower\ht0\box0}}
\put(0.00,0.16){\special{em:moveto}}
\put(51.86,0.16){\special{em:lineto}}
\put(33.46,68.33){\special{em:moveto}}
\put(33.33,0.00){\special{em:lineto}}
\put(0.00,68.50){\special{em:moveto}}
\put(0.00,0.33){\special{em:lineto}}
\put(35.46,3.83){{\setbox0=\hbox{1.014916064}\lower\ht0\box0}}
\put(35.46,10.50){{\setbox0=\hbox{0.9159655942}\lower\ht0\box0}}
\put(35.46,17.16){{\setbox0=\hbox{0.8457846720}\lower\ht0\box0}}
\put(35.46,23.83){{\setbox0=\hbox{0.6729858045}\lower\ht0\box0}}
\put(35.46,30.50){{\setbox0=\hbox{0.5562821156}\lower\ht0\box0}}
\put(35.46,37.16){{\setbox0=\hbox{0.5258402692}\lower\ht0\box0}}
\put(35.46,42.16){{\setbox0=\hbox{0.5074708032}\lower\ht0\box0}}
\put(35.46,47.16){{\setbox0=\hbox{0.4579827971}\lower\ht0\box0}}
\put(35.46,52.16){{\setbox0=\hbox{0.4228923360}\lower\ht0\box0}}
\put(35.46,58.83){{\setbox0=\hbox{0.3641071004}\lower\ht0\box0}}
\put(14.00,1.83){\special{em:moveto}}
\put(10.66,5.16){\special{em:lineto}}
\put(14.00,5.16){\special{em:moveto}}
\put(10.66,1.83){\special{em:lineto}}
\put(10.66,5.16){\special{em:moveto}}
\put(10.66,1.83){\special{em:lineto}}
\put(14.03,1.91){\special{em:moveto}}
\put(14.00,5.21){\special{em:lineto}}
\put(10.70,5.18){\special{em:moveto}}
\put(14.00,5.21){\special{em:lineto}}
\put(10.73,1.88){\special{em:moveto}}
\put(14.03,1.91){\special{em:lineto}}
\put(14.03,1.91){{\setbox0=\hbox{$\scriptstyle\bullet$}\kern-.4\wd0\lower.5\ht0\box0}}
\put(10.73,1.88){{\setbox0=\hbox{$\scriptstyle\bullet$}\kern-.4\wd0\lower.5\ht0\box0}}
\put(14.00,5.21){{\setbox0=\hbox{$\scriptstyle\bullet$}\kern-.4\wd0\lower.5\ht0\box0}}
\put(10.70,5.18){{\setbox0=\hbox{$\scriptstyle\bullet$}\kern-.4\wd0\lower.5\ht0\box0}}
\put(23.20,4.00){{\setbox0=\hbox{$H_{32}$}\lower\ht0\box0}}
\put(23.20,10.66){{\setbox0=\hbox{$H_{31}$}\lower\ht0\box0}}
\put(2.26,8.50){{\setbox0=\hbox{$\scriptstyle\bullet$}\kern-.4\wd0\lower.5\ht0\box0}}
\put(2.26,11.83){{\setbox0=\hbox{$\scriptstyle\bullet$}\kern-.4\wd0\lower.5\ht0\box0}}
\put(5.60,8.50){{\setbox0=\hbox{$\scriptstyle\bullet$}\kern-.4\wd0\lower.5\ht0\box0}}
\put(5.60,11.83){{\setbox0=\hbox{$\scriptstyle\bullet$}\kern-.4\wd0\lower.5\ht0\box0}}
\put(23.20,17.33){{\setbox0=\hbox{$H_{30}$}\lower\ht0\box0}}
\put(11.86,14.33){{\setbox0=\hbox{6}\lower\ht0\box0}}
\put(10.66,18.50){\special{em:moveto}}
\put(10.66,15.16){\special{em:lineto}}
\put(14.03,15.25){\special{em:moveto}}
\put(14.00,18.55){\special{em:lineto}}
\put(10.70,18.51){\special{em:moveto}}
\put(14.00,18.55){\special{em:lineto}}
\put(10.73,15.21){\special{em:moveto}}
\put(14.03,15.25){\special{em:lineto}}
\put(14.03,15.25){{\setbox0=\hbox{$\scriptstyle\bullet$}\kern-.4\wd0\lower.5\ht0\box0}}
\put(10.73,15.21){{\setbox0=\hbox{$\scriptstyle\bullet$}\kern-.4\wd0\lower.5\ht0\box0}}
\put(14.00,18.55){{\setbox0=\hbox{$\scriptstyle\bullet$}\kern-.4\wd0\lower.5\ht0\box0}}
\put(10.70,18.51){{\setbox0=\hbox{$\scriptstyle\bullet$}\kern-.4\wd0\lower.5\ht0\box0}}
\put(11.73,20.33){{\setbox0=\hbox{6}\lower\ht0\box0}}
\put(23.06,24.00){{\setbox0=\hbox{$H_{29}$}\lower\ht0\box0}}
\put(3.46,27.00){{\setbox0=\hbox{6}\lower\ht0\box0}}
\put(5.60,21.83){{\setbox0=\hbox{$\scriptstyle\bullet$}\kern-.4\wd0\lower.5\ht0\box0}}
\put(2.26,21.83){{\setbox0=\hbox{$\scriptstyle\bullet$}\kern-.4\wd0\lower.5\ht0\box0}}
\put(5.56,25.16){{\setbox0=\hbox{$\scriptstyle\bullet$}\kern-.4\wd0\lower.5\ht0\box0}}
\put(2.30,25.16){{\setbox0=\hbox{$\scriptstyle\bullet$}\kern-.4\wd0\lower.5\ht0\box0}}
\put(2.26,21.83){\special{em:moveto}}
\put(2.30,25.16){\special{em:lineto}}
\put(5.56,25.16){\special{em:lineto}}
\put(5.60,21.83){\special{em:lineto}}
\put(3.98,21.83){\special{em:lineto}}
\put(2.26,21.83){\special{em:lineto}}
\put(23.20,30.66){{\setbox0=\hbox{$H_{28}$}\lower\ht0\box0}}
\put(13.86,31.83){\special{em:moveto}}
\put(13.86,28.50){\special{em:lineto}}
\put(10.53,28.50){{\setbox0=\hbox{$\scriptstyle\bullet$}\kern-.4\wd0\lower.5\ht0\box0}}
\put(10.53,31.83){{\setbox0=\hbox{$\scriptstyle\bullet$}\kern-.4\wd0\lower.5\ht0\box0}}
\put(13.86,28.50){{\setbox0=\hbox{$\scriptstyle\bullet$}\kern-.4\wd0\lower.5\ht0\box0}}
\put(13.86,31.83){{\setbox0=\hbox{$\scriptstyle\bullet$}\kern-.4\wd0\lower.5\ht0\box0}}
\put(23.20,37.33){{\setbox0=\hbox{$H_{27}$}\lower\ht0\box0}}
\put(2.30,35.23){\special{em:moveto}}
\put(2.30,38.50){\special{em:lineto}}
\put(5.63,35.25){\special{em:moveto}}
\put(5.60,38.55){\special{em:lineto}}
\put(2.30,38.53){\special{em:moveto}}
\put(5.60,38.55){\special{em:lineto}}
\put(2.30,35.23){\special{em:moveto}}
\put(5.63,35.25){\special{em:lineto}}
\put(5.63,35.25){{\setbox0=\hbox{$\scriptstyle\bullet$}\kern-.4\wd0\lower.5\ht0\box0}}
\put(2.30,35.23){{\setbox0=\hbox{$\scriptstyle\bullet$}\kern-.4\wd0\lower.5\ht0\box0}}
\put(5.60,38.55){{\setbox0=\hbox{$\scriptstyle\bullet$}\kern-.4\wd0\lower.5\ht0\box0}}
\put(2.30,38.50){{\setbox0=\hbox{$\scriptstyle\bullet$}\kern-.4\wd0\lower.5\ht0\box0}}
\put(3.33,40.33){{\setbox0=\hbox{6}\lower\ht0\box0}}
\put(23.20,42.33){{\setbox0=\hbox{$H_{26}$}\lower\ht0\box0}}
\put(8.26,43.50){{\setbox0=\hbox{6}\lower\ht0\box0}}
\put(7.33,41.66){{\setbox0=\hbox{$\scriptstyle\bullet$}\kern-.4\wd0\lower.5\ht0\box0}}
\put(14.00,40.00){{\setbox0=\hbox{$\scriptstyle\bullet$}\kern-.4\wd0\lower.5\ht0\box0}}
\put(14.00,43.33){{\setbox0=\hbox{$\scriptstyle\bullet$}\kern-.4\wd0\lower.5\ht0\box0}}
\put(10.66,41.66){{\setbox0=\hbox{$\scriptstyle\bullet$}\kern-.4\wd0\lower.5\ht0\box0}}
\put(10.66,41.66){\special{em:moveto}}
\put(14.00,43.33){\special{em:lineto}}
\put(14.00,40.00){\special{em:lineto}}
\put(10.66,41.66){\special{em:lineto}}
\put(7.33,41.66){\special{em:moveto}}
\put(10.66,41.66){\special{em:lineto}}
\put(23.20,47.16){{\setbox0=\hbox{$H_{25}$}\lower\ht0\box0}}
\put(2.40,46.83){\special{em:moveto}}
\put(5.70,46.83){\special{em:lineto}}
\put(5.73,46.83){\special{em:moveto}}
\put(9.06,45.16){\special{em:lineto}}
\put(5.73,46.83){\special{em:moveto}}
\put(9.06,48.50){\special{em:lineto}}
\put(2.40,46.83){{\setbox0=\hbox{$\scriptstyle\bullet$}\kern-.4\wd0\lower.5\ht0\box0}}
\put(9.06,45.16){{\setbox0=\hbox{$\scriptstyle\bullet$}\kern-.4\wd0\lower.5\ht0\box0}}
\put(9.06,48.50){{\setbox0=\hbox{$\scriptstyle\bullet$}\kern-.4\wd0\lower.5\ht0\box0}}
\put(5.73,46.83){{\setbox0=\hbox{$\scriptstyle\bullet$}\kern-.4\wd0\lower.5\ht0\box0}}
\put(23.20,52.33){{\setbox0=\hbox{$H_{24}$}\lower\ht0\box0}}
\put(10.53,50.83){\special{em:moveto}}
\put(13.86,54.16){\special{em:lineto}}
\put(10.53,50.83){{\setbox0=\hbox{$\scriptstyle\bullet$}\kern-.4\wd0\lower.5\ht0\box0}}
\put(10.53,54.16){{\setbox0=\hbox{$\scriptstyle\bullet$}\kern-.4\wd0\lower.5\ht0\box0}}
\put(13.86,50.83){{\setbox0=\hbox{$\scriptstyle\bullet$}\kern-.4\wd0\lower.5\ht0\box0}}
\put(13.86,54.16){{\setbox0=\hbox{$\scriptstyle\bullet$}\kern-.4\wd0\lower.5\ht0\box0}}
\put(10.53,54.16){\special{em:moveto}}
\put(13.86,54.16){\special{em:lineto}}
\put(13.86,50.83){\special{em:lineto}}
\put(10.53,50.83){\special{em:lineto}}
\put(10.53,54.16){\special{em:lineto}}
\put(23.20,59.00){{\setbox0=\hbox{$H_{23}$}\lower\ht0\box0}}
\put(6.40,59.00){{\setbox0=\hbox{6}\lower\ht0\box0}}
\put(2.26,56.83){{\setbox0=\hbox{$\scriptstyle\bullet$}\kern-.4\wd0\lower.5\ht0\box0}}
\put(2.26,60.16){{\setbox0=\hbox{$\scriptstyle\bullet$}\kern-.4\wd0\lower.5\ht0\box0}}
\put(5.60,56.83){{\setbox0=\hbox{$\scriptstyle\bullet$}\kern-.4\wd0\lower.5\ht0\box0}}
\put(5.60,60.16){{\setbox0=\hbox{$\scriptstyle\bullet$}\kern-.4\wd0\lower.5\ht0\box0}}
\put(2.26,60.16){\special{em:moveto}}
\put(5.60,60.16){\special{em:lineto}}
\put(5.60,56.83){\special{em:lineto}}
\put(2.26,56.83){\special{em:lineto}}
\put(2.26,60.16){\special{em:lineto}}
\put(16.66,68.50){\special{em:moveto}}
\put(16.53,0.16){\special{em:lineto}}
\end{picture}

%% file: pic/pic_3h1.tex
\special{em:linewidth 0.010in}
\setlength{\unitlength}{0.20in}

 \begin{picture}(5,3.500000)
{\bf \put(-2,1){1}
}\multiput(0,0)(1,0){4}{\circle*{0.3}}
\scriptsize\put(-1,0){0}
\put(0,-0.1){\line(1,0){1}}
 \put(0,0){\line(1,0){1}}
 \put(0,0.1){\line(1,0){1}}
 \put(1,0){\line(1,0){1}}
 \put(2,-0.05){\line(1,0){1}}
 \put(2,0.05){\line(1,0){1}}
 \end{picture}
 \begin{picture}(5,3.500000)
\multiput(0,0)(1,0){4}{\circle*{0.3}}
\scriptsize\put(-1,0){1}
{\scriptsize {\put(-1,-1.250000){(2,3 ; 0,0,0,0)}}}
\put(0,0){\line(1,0){1}}
 \put(1.500000,0){\oval(2.800000,1.300000)[t]}
\put(1.500000,0){\oval(3.000000,1.500000)[t]}
\put(1.500000,0){\oval(3.200000,1.700000)[t]}
\put(1.500000,0.600000){\line(0,1){0.3}}
\put(1,-0.05){\line(1,0){1}}
 \put(1,0.05){\line(1,0){1}}
 \put(2.000000,0){\oval(2.000000,1.000000)[t]}
\end{picture}
 \begin{picture}(5,3.500000)
\multiput(0,0)(1,0){4}{\circle*{0.3}}
\scriptsize\put(-1,0){2}
{\scriptsize {\put(-1,-1.250000){(2,3 ; 0,0,3,3)}}}
\put(1.000000,0){\oval(2.000000,1.000000)[t]}
\put(1.500000,0){\oval(2.900000,1.400000)[t]}
\put(1.500000,0){\oval(3.100000,1.600000)[t]}
\put(1.500000,0.600000){\line(0,1){0.3}}
\put(1,-0.1){\line(1,0){1}}
 \put(1,0){\line(1,0){1}}
 \put(1,0.1){\line(1,0){1}}
 \put(2,0){\line(1,0){1}}
 \end{picture}
 \begin{picture}(5,3.500000)
\multiput(0,0)(1,0){4}{\circle*{0.3}}
\scriptsize\put(-1,0){3}
{\scriptsize {\put(-1,-1.250000){(4,4 ; 1,1,2,2)}}}
\put(0,0){\line(1,0){1}}
 \put(1.000000,0){\oval(2.000000,1.000000)[t]}
\put(1.500000,0){\oval(2.900000,1.400000)[t]}
\put(1.500000,0){\oval(3.100000,1.600000)[t]}
\put(1.500000,0.600000){\line(0,1){0.3}}
\put(1,-0.1){\line(1,0){1}}
 \put(1,0){\line(1,0){1}}
 \put(1,0.1){\line(1,0){1}}
 \put(1.500000,-0.150000){\line(0,1){0.3}}
\put(2.000000,0){\oval(2.000000,1.000000)[b]}
\put(2,0){\line(1,0){1}}
 \end{picture}
 \begin{picture}(5,3.500000)
\multiput(0,0)(1,0){4}{\circle*{0.3}}
\scriptsize\put(-1,0){4}
{\scriptsize {\put(-1,-1.250000){(4,4 ; 2,2,0,0)}}}
\put(0,0){\line(1,0){1}}
 \put(1.000000,0){\oval(1.800000,0.800000)[t]}
\put(1.000000,0){\oval(2.000000,1.000000)[t]}
\put(1.000000,0){\oval(2.200000,1.200000)[t]}
\put(1.500000,0){\oval(3.000000,1.500000)[t]}
\put(1.500000,0.600000){\line(0,1){0.3}}
\put(2.000000,0){\oval(2.000000,1.000000)[b]}
\put(2,-0.1){\line(1,0){1}}
 \put(2,0){\line(1,0){1}}
 \put(2,0.1){\line(1,0){1}}
 \end{picture}

 \begin{picture}(5,3.500000)
\multiput(0,0)(1,0){4}{\circle*{0.3}}
\scriptsize\put(-1,0){5}
{\scriptsize {\put(-1,-1.250000){(8,5 ; 4,4,0,0)}}}
\put(0,-0.1){\line(1,0){1}}
 \put(0,0){\line(1,0){1}}
 \put(0,0.1){\line(1,0){1}}
 \put(1.500000,0){\oval(2.800000,1.300000)[t]}
\put(1.500000,0){\oval(3.000000,1.500000)[t]}
\put(1.500000,0){\oval(3.200000,1.700000)[t]}
\put(1.500000,0.600000){\line(0,1){0.3}}
\put(1,0){\line(1,0){1}}
 \put(2,-0.1){\line(1,0){1}}
 \put(2,0){\line(1,0){1}}
 \put(2,0.1){\line(1,0){1}}
 \end{picture}
\\

 \begin{picture}(5,3.500000)
{\bf \put(-2,1){2}
}\multiput(0,0)(1,0){4}{\circle*{0.3}}
\scriptsize\put(-1,0){0}
\put(0,0){\line(1,0){1}}
 \put(1,-0.1){\line(1,0){1}}
 \put(1,0){\line(1,0){1}}
 \put(1,0.1){\line(1,0){1}}
 \put(2,0){\line(1,0){1}}
 \end{picture}
 \begin{picture}(5,3.500000)
\multiput(0,0)(1,0){4}{\circle*{0.3}}
\scriptsize\put(-1,0){1}
{\scriptsize {\put(-1,-1.250000){(2,3 ; 0,0,0,0)}}}
\put(0,-0.1){\line(1,0){1}}
 \put(0,0){\line(1,0){1}}
 \put(0,0.1){\line(1,0){1}}
 \put(1.500000,0){\oval(3.000000,1.500000)[t]}
\put(1.500000,0.600000){\line(0,1){0.3}}
\put(1,0){\line(1,0){1}}
 \put(2.000000,0){\oval(1.800000,0.800000)[t]}
\put(2.000000,0){\oval(2.000000,1.000000)[t]}
\put(2.000000,0){\oval(2.200000,1.200000)[t]}
\end{picture}
 \begin{picture}(5,3.500000)
\multiput(0,0)(1,0){4}{\circle*{0.3}}
\scriptsize\put(-1,0){2}
{\scriptsize {\put(-1,-1.250000){(3,4 ; 1,0,0,1)}}}
\put(1.000000,0){\oval(2.000000,1.000000)[t]}
\put(1.500000,0){\oval(2.800000,1.300000)[t]}
\put(1.500000,0){\oval(3.000000,1.500000)[t]}
\put(1.500000,0){\oval(3.200000,1.700000)[t]}
\put(1.500000,0.600000){\line(0,1){0.3}}
\put(2.000000,0){\oval(1.800000,0.800000)[b]}
\put(2.000000,0){\oval(2.000000,1.000000)[b]}
\put(2.000000,0){\oval(2.200000,1.200000)[b]}
\put(2,0){\line(1,0){1}}
 \end{picture}
 \begin{picture}(5,3.500000)
\multiput(0,0)(1,0){4}{\circle*{0.3}}
\scriptsize\put(-1,0){3}
{\scriptsize {\put(-1,-1.250000){(4,4 ; 1,1,2,2)}}}
\put(0,-0.1){\line(1,0){1}}
 \put(0,0){\line(1,0){1}}
 \put(0,0.1){\line(1,0){1}}
 \put(1.000000,0){\oval(1.800000,0.800000)[t]}
\put(1.000000,0){\oval(2.000000,1.000000)[t]}
\put(1.000000,0){\oval(2.200000,1.200000)[t]}
\put(1.500000,0){\oval(3.000000,1.500000)[t]}
\put(1.500000,0.600000){\line(0,1){0.3}}
\put(1,0){\line(1,0){1}}
 \put(1.500000,-0.150000){\line(0,1){0.3}}
\put(2.000000,0){\oval(1.800000,0.800000)[b]}
\put(2.000000,0){\oval(2.000000,1.000000)[b]}
\put(2.000000,0){\oval(2.200000,1.200000)[b]}
\put(2,-0.1){\line(1,0){1}}
 \put(2,0){\line(1,0){1}}
 \put(2,0.1){\line(1,0){1}}
 \end{picture}
 \begin{picture}(5,3.500000)
\multiput(0,0)(1,0){4}{\circle*{0.3}}
\scriptsize\put(-1,0){4}
{\scriptsize {\put(-1,-1.250000){(5,5 ; 2,1,0,0)}}}
\put(0,0){\line(1,0){1}}
 \put(1.000000,0){\oval(1.800000,0.800000)[t]}
\put(1.000000,0){\oval(2.000000,1.000000)[t]}
\put(1.000000,0){\oval(2.200000,1.200000)[t]}
\put(1.500000,0){\oval(2.800000,1.300000)[t]}
\put(1.500000,0){\oval(3.000000,1.500000)[t]}
\put(1.500000,0){\oval(3.200000,1.700000)[t]}
\put(1.300000,0.600000){\line(0,1){0.3}}
\put(1.700000,0.600000){\line(0,1){0.3}}
\put(2.000000,0){\oval(1.800000,0.800000)[b]}
\put(2.000000,0){\oval(2.000000,1.000000)[b]}
\put(2.000000,0){\oval(2.200000,1.200000)[b]}
\put(2,0){\line(1,0){1}}
 \end{picture}

 \begin{picture}(5,3.500000)
\multiput(0,0)(1,0){4}{\circle*{0.3}}
\scriptsize\put(-1,0){5}
{\scriptsize {\put(-1,-1.250000){(6,5 ; 2,2,1,1)}}}
\put(0,-0.1){\line(1,0){1}}
 \put(0,0){\line(1,0){1}}
 \put(0,0.1){\line(1,0){1}}
 \put(0.500000,-0.150000){\line(0,1){0.3}}
\put(1.000000,0){\oval(2.000000,1.000000)[t]}
\put(1.500000,0){\oval(2.800000,1.300000)[t]}
\put(1.500000,0){\oval(3.000000,1.500000)[t]}
\put(1.500000,0){\oval(3.200000,1.700000)[t]}
\put(1.500000,0.600000){\line(0,1){0.3}}
\put(1,0){\line(1,0){1}}
 \put(2.000000,0){\oval(1.800000,0.800000)[b]}
\put(2.000000,0){\oval(2.000000,1.000000)[b]}
\put(2.000000,0){\oval(2.200000,1.200000)[b]}
\put(2.000000,-0.650000){\line(0,1){0.3}}
\put(2,0){\line(1,0){1}}
 \end{picture}
\\
 \begin{picture}(5,3.500000)
{\bf \put(-2,1){3}
}\multiput(0,0)(1,0){4}{\circle*{0.3}}
\scriptsize\put(-1,0){0}
\put(0,-0.1){\line(1,0){1}}
 \put(0,0){\line(1,0){1}}
 \put(0,0.1){\line(1,0){1}}
 \put(1.000000,0){\oval(2.000000,1.000000)[t]}
\put(1.500000,0){\oval(3.000000,1.500000)[t]}
\end{picture}
 \begin{picture}(5,3.500000)
\multiput(0,0)(1,0){4}{\circle*{0.3}}
\scriptsize\put(-1,0){1}
{\scriptsize {\put(-1,-1.250000){(2,3 ; 0,0,1,1)}}}
\put(0,0){\line(1,0){1}}
 \put(1.000000,0){\oval(2.000000,1.000000)[t]}
\put(1.500000,0){\oval(2.800000,1.300000)[t]}
\put(1.500000,0){\oval(3.000000,1.500000)[t]}
\put(1.500000,0){\oval(3.200000,1.700000)[t]}
\put(1.500000,0.600000){\line(0,1){0.3}}
\put(2.000000,0){\oval(2.000000,1.000000)[b]}
\put(2,0){\line(1,0){1}}
 \end{picture}
 \begin{picture}(5,3.500000)
\multiput(0,0)(1,0){4}{\circle*{0.3}}
\scriptsize\put(-1,0){2}
{\scriptsize {\put(-1,-1.250000){(2,3 ; 0,0,2,2)}}}
\put(0,0){\line(1,0){1}}
 \put(1.000000,0){\oval(1.800000,0.800000)[t]}
\put(1.000000,0){\oval(2.000000,1.000000)[t]}
\put(1.000000,0){\oval(2.200000,1.200000)[t]}
\put(1.500000,0){\oval(3.000000,1.500000)[t]}
\put(1.500000,0.600000){\line(0,1){0.3}}
\put(2.000000,0){\oval(2.000000,1.000000)[b]}
\put(2,-0.1){\line(1,0){1}}
 \put(2,0){\line(1,0){1}}
 \put(2,0.1){\line(1,0){1}}
 \end{picture}
 \begin{picture}(5,3.500000)
\multiput(0,0)(1,0){4}{\circle*{0.3}}
\scriptsize\put(-1,0){3}
{\scriptsize {\put(-1,-1.250000){(4,4 ; 2,2,0,0)}}}
\put(0,-0.1){\line(1,0){1}}
 \put(0,0){\line(1,0){1}}
 \put(0,0.1){\line(1,0){1}}
 \put(1.500000,0){\oval(2.800000,1.300000)[t]}
\put(1.500000,0){\oval(3.000000,1.500000)[t]}
\put(1.500000,0){\oval(3.200000,1.700000)[t]}
\put(1.500000,0.600000){\line(0,1){0.3}}
\put(1,0){\line(1,0){1}}
 \put(2,-0.1){\line(1,0){1}}
 \put(2,0){\line(1,0){1}}
 \put(2,0.1){\line(1,0){1}}
 \end{picture}
\\

 \begin{picture}(5,3.500000)
{\bf \put(-2,1){4}
}\multiput(0,0)(1,0){4}{\circle*{0.3}}
\scriptsize\put(-1,0){0}
\put(0,-0.1){\line(1,0){1}}
 \put(0,0){\line(1,0){1}}
 \put(0,0.1){\line(1,0){1}}
 \put(1,0){\line(1,0){1}}
 \put(2,-0.1){\line(1,0){1}}
 \put(2,0){\line(1,0){1}}
 \put(2,0.1){\line(1,0){1}}
 \end{picture}
 \begin{picture}(5,3.500000)
\multiput(0,0)(1,0){4}{\circle*{0.3}}
\scriptsize\put(-1,0){1}
{\scriptsize {\put(-1,-1.250000){(2,3 ; 0,0,0,0)}}}
\put(0,0){\line(1,0){1}}
 \put(1.500000,0){\oval(2.800000,1.300000)[t]}
\put(1.500000,0){\oval(3.000000,1.500000)[t]}
\put(1.500000,0){\oval(3.200000,1.700000)[t]}
\put(1.500000,0.600000){\line(0,1){0.3}}
\put(1,-0.1){\line(1,0){1}}
 \put(1,0){\line(1,0){1}}
 \put(1,0.1){\line(1,0){1}}
 \put(2.000000,0){\oval(2.000000,1.000000)[t]}
\end{picture}
 \begin{picture}(5,3.500000)
\multiput(0,0)(1,0){4}{\circle*{0.3}}
\scriptsize\put(-1,0){2}
{\scriptsize {\put(-1,-1.250000){(4,4 ; 1,1,2,2)}}}
\put(0,0){\line(1,0){1}}
 \put(1.000000,0){\oval(2.000000,1.000000)[t]}
\put(1.500000,0){\oval(2.800000,1.300000)[t]}
\put(1.500000,0){\oval(3.000000,1.500000)[t]}
\put(1.500000,0){\oval(3.200000,1.700000)[t]}
\put(1.500000,0.600000){\line(0,1){0.3}}
\put(1,-0.1){\line(1,0){1}}
 \put(1,0){\line(1,0){1}}
 \put(1,0.1){\line(1,0){1}}
 \put(1.500000,-0.150000){\line(0,1){0.3}}
\put(2.000000,0){\oval(2.000000,1.000000)[b]}
\put(2,0){\line(1,0){1}}
 \end{picture}

%% file: pic/pic_3h2.tex
\special{em:linewidth 0.010in}
\setlength{\unitlength}{0.2in}

 \begin{picture}(5,3.500000)
{\bf \put(-2,1){1}
}\multiput(0,0)(1,0){4}{\circle*{0.3}}
\scriptsize\put(-1,0){0}
\put(0,0){\line(1,0){1}}
 \put(1,0){\line(1,0){1}}
 \put(2,-0.1){\line(1,0){1}}
 \put(2,-0.033){\line(1,0){1}}
 \put(2,0.033){\line(1,0){1}}
 \put(2,0.1){\line(1,0){1}}
 \end{picture}
 \begin{picture}(5,3.500000)
\multiput(0,0)(1,0){4}{\circle*{0.3}}
\scriptsize\put(-1,0){1}
{\scriptsize {\put(-1,-1.250000){(2,3 ; 0,0,0,0)}}}
\put(0,0){\line(1,0){1}}
 \put(1.500000,0){\oval(3.000000,1.500000)[t]}
\put(1.500000,0.600000){\line(0,1){0.3}}
\put(1,-0.1){\line(1,0){1}}
 \put(1,-0.033){\line(1,0){1}}
 \put(1,0.033){\line(1,0){1}}
 \put(1,0.1){\line(1,0){1}}
 \put(2.000000,0){\oval(2.000000,1.000000)[t]}
\end{picture}
 \begin{picture}(5,3.500000)
\multiput(0,0)(1,0){4}{\circle*{0.3}}
\scriptsize\put(-1,0){2}
{\scriptsize {\put(-1,-1.250000){(2,3 ; 0,0,3,3)}}}
\put(1.000000,0){\oval(2.000000,1.000000)[t]}
\put(1.500000,0){\oval(3.300000,1.800000)[t]}
\put(1.500000,0){\oval(3.100000,1.600000)[t]}
\put(1.500000,0){\oval(2.900000,1.400000)[t]}
\put(1.500000,0){\oval(2.700000,1.200000)[t]}
\put(1.500000,0.600000){\line(0,1){0.3}}
\put(1,0){\line(1,0){1}}
 \put(2,0){\line(1,0){1}}
 \end{picture}
 \begin{picture}(5,3.500000)
\multiput(0,0)(1,0){4}{\circle*{0.3}}
\scriptsize\put(-1,0){3}
{\scriptsize {\put(-1,-1.250000){(3,4 ; 1,0,0,1)}}}
\put(1.000000,0){\oval(2.300000,1.300000)[t]}
\put(1.000000,0){\oval(2.100000,1.100000)[t]}
\put(1.000000,0){\oval(1.900000,0.900000)[t]}
\put(1.000000,0){\oval(1.700000,0.700000)[t]}
\put(1.500000,0){\oval(3.000000,1.500000)[t]}
\put(1.500000,0.600000){\line(0,1){0.3}}
\put(2.000000,0){\oval(2.000000,1.000000)[b]}
\put(2,-0.1){\line(1,0){1}}
 \put(2,-0.033){\line(1,0){1}}
 \put(2,0.033){\line(1,0){1}}
 \put(2,0.1){\line(1,0){1}}
 \end{picture}
 \begin{picture}(5,3.500000)
\multiput(0,0)(1,0){4}{\circle*{0.3}}
\scriptsize\put(-1,0){4}
{\scriptsize {\put(-1,-1.250000){(4,4 ; 1,1,2,2)}}}
\put(0,0){\line(1,0){1}}
 \put(1.000000,0){\oval(2.000000,1.000000)[t]}
\put(1.500000,0){\oval(3.300000,1.800000)[t]}
\put(1.500000,0){\oval(3.100000,1.600000)[t]}
\put(1.500000,0){\oval(2.900000,1.400000)[t]}
\put(1.500000,0){\oval(2.700000,1.200000)[t]}
\put(1.500000,0.600000){\line(0,1){0.3}}
\put(1,0){\line(1,0){1}}
 \put(1.500000,-0.150000){\line(0,1){0.3}}
\put(2.000000,0){\oval(2.000000,1.000000)[b]}
\put(2,0){\line(1,0){1}}
 \end{picture}

 \begin{picture}(5,3.500000)
\multiput(0,0)(1,0){4}{\circle*{0.3}}
\scriptsize\put(-1,0){5}
{\scriptsize {\put(-1,-1.250000){(4,5 ; 3,0,0,2)}}}
\put(0,-0.1){\line(1,0){1}}
 \put(0,-0.033){\line(1,0){1}}
 \put(0,0.033){\line(1,0){1}}
 \put(0,0.1){\line(1,0){1}}
 \put(1.500000,0){\oval(3.300000,1.800000)[t]}
\put(1.500000,0){\oval(3.100000,1.600000)[t]}
\put(1.500000,0){\oval(2.900000,1.400000)[t]}
\put(1.500000,0){\oval(2.700000,1.200000)[t]}
\put(1.500000,0.600000){\line(0,1){0.3}}
\put(1,0){\line(1,0){1}}
 \end{picture}
 \begin{picture}(5,3.500000)
\multiput(0,0)(1,0){4}{\circle*{0.3}}
\scriptsize\put(-1,0){6}
{\scriptsize {\put(-1,-1.250000){(5,5 ; 3,2,0,0)}}}
\put(1.000000,0){\oval(2.000000,1.000000)[t]}
\put(1.500000,0){\oval(3.300000,1.800000)[t]}
\put(1.500000,0){\oval(3.100000,1.600000)[t]}
\put(1.500000,0){\oval(2.900000,1.400000)[t]}
\put(1.500000,0){\oval(2.700000,1.200000)[t]}
\put(1.300000,0.600000){\line(0,1){0.3}}
\put(1.700000,0.600000){\line(0,1){0.3}}
\put(1,0){\line(1,0){1}}
 \put(2,-0.1){\line(1,0){1}}
 \put(2,-0.033){\line(1,0){1}}
 \put(2,0.033){\line(1,0){1}}
 \put(2,0.1){\line(1,0){1}}
 \end{picture}
 \begin{picture}(5,3.500000)
\multiput(0,0)(1,0){4}{\circle*{0.3}}
\scriptsize\put(-1,0){7}
{\scriptsize {\put(-1,-1.250000){(6,5 ; 3,3,1,1)}}}
\put(0,0){\line(1,0){1}}
 \put(0.500000,-0.150000){\line(0,1){0.3}}
\put(1.000000,0){\oval(2.300000,1.300000)[t]}
\put(1.000000,0){\oval(2.100000,1.100000)[t]}
\put(1.000000,0){\oval(1.900000,0.900000)[t]}
\put(1.000000,0){\oval(1.700000,0.700000)[t]}
\put(1.500000,0){\oval(3.000000,1.500000)[t]}
\put(1.500000,0.600000){\line(0,1){0.3}}
\put(1,-0.1){\line(1,0){1}}
 \put(1,-0.033){\line(1,0){1}}
 \put(1,0.033){\line(1,0){1}}
 \put(1,0.1){\line(1,0){1}}
 \put(2.000000,0){\oval(2.000000,1.000000)[b]}
\put(2.000000,-0.650000){\line(0,1){0.3}}
\put(2,-0.1){\line(1,0){1}}
 \put(2,-0.033){\line(1,0){1}}
 \put(2,0.033){\line(1,0){1}}
 \put(2,0.1){\line(1,0){1}}
 \end{picture}
 \begin{picture}(5,3.500000)
\multiput(0,0)(1,0){4}{\circle*{0.3}}
\scriptsize\put(-1,0){8}
{\scriptsize {\put(-1,-1.250000){(6,5 ; 3,3,3,3)}}}
\put(0,-0.1){\line(1,0){1}}
 \put(0,-0.033){\line(1,0){1}}
 \put(0,0.033){\line(1,0){1}}
 \put(0,0.1){\line(1,0){1}}
 \put(1.500000,0){\oval(3.300000,1.800000)[t]}
\put(1.500000,0){\oval(3.100000,1.600000)[t]}
\put(1.500000,0){\oval(2.900000,1.400000)[t]}
\put(1.500000,0){\oval(2.700000,1.200000)[t]}
\put(1.500000,0.600000){\line(0,1){0.3}}
\put(2,-0.1){\line(1,0){1}}
 \put(2,-0.033){\line(1,0){1}}
 \put(2,0.033){\line(1,0){1}}
 \put(2,0.1){\line(1,0){1}}
 \end{picture}
 \begin{picture}(5,3.500000)
\multiput(0,0)(1,0){4}{\circle*{0.3}}
\scriptsize\put(-1,0){9}
{\scriptsize {\put(-1,-1.250000){(8,6 ; 5,5,2,2)}}}
\put(0,-0.1){\line(1,0){1}}
 \put(0,-0.033){\line(1,0){1}}
 \put(0,0.033){\line(1,0){1}}
 \put(0,0.1){\line(1,0){1}}
 \put(1.500000,0){\oval(3.000000,1.500000)[t]}
\put(1.500000,0.600000){\line(0,1){0.3}}
\put(1,-0.1){\line(1,0){1}}
 \put(1,-0.033){\line(1,0){1}}
 \put(1,0.033){\line(1,0){1}}
 \put(1,0.1){\line(1,0){1}}
 \put(1.500000,-0.150000){\line(0,1){0.3}}
\put(2.000000,0){\oval(2.300000,1.300000)[t]}
\put(2.000000,0){\oval(2.100000,1.100000)[t]}
\put(2.000000,0){\oval(1.900000,0.900000)[t]}
\put(2.000000,0){\oval(1.700000,0.700000)[t]}
\end{picture}

 \begin{picture}(5,3.500000)
\multiput(0,0)(1,0){4}{\circle*{0.3}}
\scriptsize\put(-1,0){10}
{\scriptsize {\put(-1,-1.250000){(8,6 ; 5,5,3,3)}}}
\put(1.000000,0){\oval(2.300000,1.300000)[t]}
\put(1.000000,0){\oval(2.100000,1.100000)[t]}
\put(1.000000,0){\oval(1.900000,0.900000)[t]}
\put(1.000000,0){\oval(1.700000,0.700000)[t]}
\put(1.500000,0){\oval(3.300000,1.800000)[t]}
\put(1.500000,0){\oval(3.100000,1.600000)[t]}
\put(1.500000,0){\oval(2.900000,1.400000)[t]}
\put(1.500000,0){\oval(2.700000,1.200000)[t]}
\put(1.350000,0.600000){\line(0,1){0.3}}
\put(1.500000,0.600000){\line(0,1){0.3}}
\put(1.650000,0.600000){\line(0,1){0.3}}
\put(1,0){\line(1,0){1}}
 \put(2,-0.1){\line(1,0){1}}
 \put(2,-0.033){\line(1,0){1}}
 \put(2,0.033){\line(1,0){1}}
 \put(2,0.1){\line(1,0){1}}
 \end{picture}
 \begin{picture}(5,3.500000)
\multiput(0,0)(1,0){4}{\circle*{0.3}}
\scriptsize\put(-1,0){11}
{\scriptsize {\put(-1,-1.250000){(10,6 ; 6,6,0,0)}}}
\put(0,0){\line(1,0){1}}
 \put(1.000000,0){\oval(2.300000,1.300000)[t]}
\put(1.000000,0){\oval(2.100000,1.100000)[t]}
\put(1.000000,0){\oval(1.900000,0.900000)[t]}
\put(1.000000,0){\oval(1.700000,0.700000)[t]}
\put(1.500000,0){\oval(3.000000,1.500000)[t]}
\put(1.500000,0.600000){\line(0,1){0.3}}
\put(2.000000,0){\oval(2.000000,1.000000)[b]}
\put(2,-0.1){\line(1,0){1}}
 \put(2,-0.033){\line(1,0){1}}
 \put(2,0.033){\line(1,0){1}}
 \put(2,0.1){\line(1,0){1}}
 \end{picture}
 \begin{picture}(5,3.500000)
\multiput(0,0)(1,0){4}{\circle*{0.3}}
\scriptsize\put(-1,0){12}
{\scriptsize {\put(-1,-1.250000){(10,6 ; 6,6,1,1)}}}
\put(0,-0.1){\line(1,0){1}}
 \put(0,-0.033){\line(1,0){1}}
 \put(0,0.033){\line(1,0){1}}
 \put(0,0.1){\line(1,0){1}}
 \put(1.000000,0){\oval(2.000000,1.000000)[t]}
\put(1.500000,0){\oval(3.000000,1.500000)[t]}
\put(1.500000,0.600000){\line(0,1){0.3}}
\put(1,-0.1){\line(1,0){1}}
 \put(1,-0.033){\line(1,0){1}}
 \put(1,0.033){\line(1,0){1}}
 \put(1,0.1){\line(1,0){1}}
 \put(1.300000,-0.150000){\line(0,1){0.3}}
\put(1.700000,-0.150000){\line(0,1){0.3}}
\put(2.000000,0){\oval(2.300000,1.300000)[b]}
\put(2.000000,0){\oval(2.100000,1.100000)[b]}
\put(2.000000,0){\oval(1.900000,0.900000)[b]}
\put(2.000000,0){\oval(1.700000,0.700000)[b]}
\put(2,0){\line(1,0){1}}
 \end{picture}
 \begin{picture}(5,3.500000)
\multiput(0,0)(1,0){4}{\circle*{0.3}}
\scriptsize\put(-1,0){13}
{\scriptsize {\put(-1,-1.250000){(10,6 ; 6,6,3,3)}}}
\put(0,0){\line(1,0){1}}
 \put(1.000000,0){\oval(2.000000,1.000000)[t]}
\put(1.500000,0){\oval(3.300000,1.800000)[t]}
\put(1.500000,0){\oval(3.100000,1.600000)[t]}
\put(1.500000,0){\oval(2.900000,1.400000)[t]}
\put(1.500000,0){\oval(2.700000,1.200000)[t]}
\put(1.500000,0.600000){\line(0,1){0.3}}
\put(2.000000,0){\oval(2.000000,1.000000)[b]}
\put(2,0){\line(1,0){1}}
 \end{picture}
 \begin{picture}(5,3.500000)
\multiput(0,0)(1,0){4}{\circle*{0.3}}
\scriptsize\put(-1,0){14}
{\scriptsize {\put(-1,-1.250000){(12,6 ; 8,8,1,1)}}}
\put(1.000000,0){\oval(2.300000,1.300000)[t]}
\put(1.000000,0){\oval(2.100000,1.100000)[t]}
\put(1.000000,0){\oval(1.900000,0.900000)[t]}
\put(1.000000,0){\oval(1.700000,0.700000)[t]}
\put(1.000000,0.350000){\line(0,1){0.3}}
\put(1.500000,0){\oval(3.300000,1.800000)[t]}
\put(1.500000,0){\oval(3.100000,1.600000)[t]}
\put(1.500000,0){\oval(2.900000,1.400000)[t]}
\put(1.500000,0){\oval(2.700000,1.200000)[t]}
\put(1.500000,0.600000){\line(0,1){0.3}}
\put(1,-0.1){\line(1,0){1}}
 \put(1,-0.033){\line(1,0){1}}
 \put(1,0.033){\line(1,0){1}}
 \put(1,0.1){\line(1,0){1}}
 \put(2,-0.1){\line(1,0){1}}
 \put(2,-0.033){\line(1,0){1}}
 \put(2,0.033){\line(1,0){1}}
 \put(2,0.1){\line(1,0){1}}
 \put(2.500000,-0.150000){\line(0,1){0.3}}
\end{picture}

 \begin{picture}(5,3.500000)
\multiput(0,0)(1,0){4}{\circle*{0.3}}
\scriptsize\put(-1,0){15}
{\scriptsize {\put(-1,-1.250000){(16,7 ; 9,9,2,2)}}}
\put(0,-0.1){\line(1,0){1}}
 \put(0,-0.033){\line(1,0){1}}
 \put(0,0.033){\line(1,0){1}}
 \put(0,0.1){\line(1,0){1}}
 \put(1.000000,0){\oval(2.300000,1.300000)[t]}
\put(1.000000,0){\oval(2.100000,1.100000)[t]}
\put(1.000000,0){\oval(1.900000,0.900000)[t]}
\put(1.000000,0){\oval(1.700000,0.700000)[t]}
\put(1.500000,0){\oval(3.300000,1.800000)[t]}
\put(1.500000,0){\oval(3.100000,1.600000)[t]}
\put(1.500000,0){\oval(2.900000,1.400000)[t]}
\put(1.500000,0){\oval(2.700000,1.200000)[t]}
\put(1.350000,0.600000){\line(0,1){0.3}}
\put(1.500000,0.600000){\line(0,1){0.3}}
\put(1.650000,0.600000){\line(0,1){0.3}}
\put(1,0){\line(1,0){1}}
 \put(1.500000,-0.150000){\line(0,1){0.3}}
\put(2.000000,0){\oval(2.300000,1.300000)[b]}
\put(2.000000,0){\oval(2.100000,1.100000)[b]}
\put(2.000000,0){\oval(1.900000,0.900000)[b]}
\put(2.000000,0){\oval(1.700000,0.700000)[b]}
\put(2,-0.1){\line(1,0){1}}
 \put(2,-0.033){\line(1,0){1}}
 \put(2,0.033){\line(1,0){1}}
 \put(2,0.1){\line(1,0){1}}
 \end{picture}
 \begin{picture}(5,3.500000)
\multiput(0,0)(1,0){4}{\circle*{0.3}}
\scriptsize\put(-1,0){16}
{\scriptsize {\put(-1,-1.250000){(12,7 ; 10,5,0,0)}}}
\put(0,-0.1){\line(1,0){1}}
 \put(0,-0.033){\line(1,0){1}}
 \put(0,0.033){\line(1,0){1}}
 \put(0,0.1){\line(1,0){1}}
 \put(1.000000,0){\oval(2.000000,1.000000)[t]}
\put(1.500000,0){\oval(3.300000,1.800000)[t]}
\put(1.500000,0){\oval(3.100000,1.600000)[t]}
\put(1.500000,0){\oval(2.900000,1.400000)[t]}
\put(1.500000,0){\oval(2.700000,1.200000)[t]}
\put(1.500000,0.600000){\line(0,1){0.3}}
\put(2,0){\line(1,0){1}}
 \end{picture}
 \begin{picture}(5,3.500000)
\multiput(0,0)(1,0){4}{\circle*{0.3}}
\scriptsize\put(-1,0){17}
{\scriptsize {\put(-1,-1.250000){(20,7 ; 12,11,0,0)}}}
\put(1.000000,0){\oval(2.300000,1.300000)[t]}
\put(1.000000,0){\oval(2.100000,1.100000)[t]}
\put(1.000000,0){\oval(1.900000,0.900000)[t]}
\put(1.000000,0){\oval(1.700000,0.700000)[t]}
\put(1.500000,0){\oval(3.300000,1.800000)[t]}
\put(1.500000,0){\oval(3.100000,1.600000)[t]}
\put(1.500000,0){\oval(2.900000,1.400000)[t]}
\put(1.500000,0){\oval(2.700000,1.200000)[t]}
\put(1.500000,0.600000){\line(0,1){0.3}}
\put(1,0){\line(1,0){1}}
 \put(2.000000,0){\oval(2.300000,1.300000)[b]}
\put(2.000000,0){\oval(2.100000,1.100000)[b]}
\put(2.000000,0){\oval(1.900000,0.900000)[b]}
\put(2.000000,0){\oval(1.700000,0.700000)[b]}
\put(2,-0.1){\line(1,0){1}}
 \put(2,-0.033){\line(1,0){1}}
 \put(2,0.033){\line(1,0){1}}
 \put(2,0.1){\line(1,0){1}}
 \put(2.300000,-0.150000){\line(0,1){0.3}}
\put(2.700000,-0.150000){\line(0,1){0.3}}
\end{picture}
 \begin{picture}(5,3.500000)
\multiput(0,0)(1,0){4}{\circle*{0.3}}
\scriptsize\put(-1,0){18}
{\scriptsize {\put(-1,-1.250000){(24,7 ; 14,14,1,1)}}}
\put(0,-0.1){\line(1,0){1}}
 \put(0,-0.033){\line(1,0){1}}
 \put(0,0.033){\line(1,0){1}}
 \put(0,0.1){\line(1,0){1}}
 \put(0.500000,-0.150000){\line(0,1){0.3}}
\put(1.000000,0){\oval(2.300000,1.300000)[t]}
\put(1.000000,0){\oval(2.100000,1.100000)[t]}
\put(1.000000,0){\oval(1.900000,0.900000)[t]}
\put(1.000000,0){\oval(1.700000,0.700000)[t]}
\put(1.000000,0.350000){\line(0,1){0.3}}
\put(1.500000,0){\oval(3.300000,1.800000)[t]}
\put(1.500000,0){\oval(3.100000,1.600000)[t]}
\put(1.500000,0){\oval(2.900000,1.400000)[t]}
\put(1.500000,0){\oval(2.700000,1.200000)[t]}
\put(1.500000,0.600000){\line(0,1){0.3}}
\put(1,-0.1){\line(1,0){1}}
 \put(1,-0.033){\line(1,0){1}}
 \put(1,0.033){\line(1,0){1}}
 \put(1,0.1){\line(1,0){1}}
 \put(1.500000,-0.150000){\line(0,1){0.3}}
\put(2.000000,0){\oval(2.300000,1.300000)[b]}
\put(2.000000,0){\oval(2.100000,1.100000)[b]}
\put(2.000000,0){\oval(1.900000,0.900000)[b]}
\put(2.000000,0){\oval(1.700000,0.700000)[b]}
\put(2.000000,-0.650000){\line(0,1){0.3}}
\put(2,-0.1){\line(1,0){1}}
 \put(2,-0.033){\line(1,0){1}}
 \put(2,0.033){\line(1,0){1}}
 \put(2,0.1){\line(1,0){1}}
 \put(2.500000,-0.150000){\line(0,1){0.3}}
\end{picture}
 \begin{picture}(5,3.500000)
\multiput(0,0)(1,0){4}{\circle*{0.3}}
\scriptsize\put(-1,0){19}
{\scriptsize {\put(-1,-1.250000){(24,8 ; 16,16,1,1)}}}
\put(0,-0.1){\line(1,0){1}}
 \put(0,-0.033){\line(1,0){1}}
 \put(0,0.033){\line(1,0){1}}
 \put(0,0.1){\line(1,0){1}}
 \put(0.500000,-0.150000){\line(0,1){0.3}}
\put(1.000000,0){\oval(2.000000,1.000000)[t]}
\put(1.500000,0){\oval(3.300000,1.800000)[t]}
\put(1.500000,0){\oval(3.100000,1.600000)[t]}
\put(1.500000,0){\oval(2.900000,1.400000)[t]}
\put(1.500000,0){\oval(2.700000,1.200000)[t]}
\put(1.500000,0.600000){\line(0,1){0.3}}
\put(1,0){\line(1,0){1}}
 \put(2.000000,0){\oval(2.300000,1.300000)[b]}
\put(2.000000,0){\oval(2.100000,1.100000)[b]}
\put(2.000000,0){\oval(1.900000,0.900000)[b]}
\put(2.000000,0){\oval(1.700000,0.700000)[b]}
\put(2.000000,-0.650000){\line(0,1){0.3}}
\put(2,0){\line(1,0){1}}
 \end{picture}
\\

 \begin{picture}(5,3.500000)
{\bf \put(-2,1){2}
}\multiput(0,0)(1,0){4}{\circle*{0.3}}
\scriptsize\put(-1,0){0}
\put(0,0){\line(1,0){1}}
 \put(1,-0.05){\line(1,0){1}}
 \put(1,0.05){\line(1,0){1}}
 \put(2,-0.05){\line(1,0){1}}
 \put(2,0.05){\line(1,0){1}}
 \end{picture}
 \begin{picture}(5,3.500000)
\multiput(0,0)(1,0){4}{\circle*{0.3}}
\scriptsize\put(-1,0){1}
{\scriptsize {\put(-1,-1.250000){(2,3 ; 0,0,0,0)}}}
\put(0,-0.05){\line(1,0){1}}
 \put(0,0.05){\line(1,0){1}}
 \put(1.500000,0){\oval(3.000000,1.500000)[t]}
\put(1.500000,0.600000){\line(0,1){0.3}}
\put(1,-0.05){\line(1,0){1}}
 \put(1,0.05){\line(1,0){1}}
 \put(2.000000,0){\oval(1.900000,0.900000)[t]}
\put(2.000000,0){\oval(2.100000,1.100000)[t]}
\end{picture}
 \begin{picture}(5,3.500000)
\multiput(0,0)(1,0){4}{\circle*{0.3}}
\scriptsize\put(-1,0){2}
{\scriptsize {\put(-1,-1.250000){(2,3 ; 0,0,3,3)}}}
\put(1.000000,0){\oval(1.900000,0.900000)[t]}
\put(1.000000,0){\oval(2.100000,1.100000)[t]}
\put(1.500000,0){\oval(2.900000,1.400000)[t]}
\put(1.500000,0){\oval(3.100000,1.600000)[t]}
\put(1.500000,0.600000){\line(0,1){0.3}}
\put(1,0){\line(1,0){1}}
 \put(2,-0.05){\line(1,0){1}}
 \put(2,0.05){\line(1,0){1}}
 \end{picture}
 \begin{picture}(5,3.500000)
\multiput(0,0)(1,0){4}{\circle*{0.3}}
\scriptsize\put(-1,0){3}
{\scriptsize {\put(-1,-1.250000){(3,4 ; 1,0,0,1)}}}
\put(1.000000,0){\oval(1.900000,0.900000)[t]}
\put(1.000000,0){\oval(2.100000,1.100000)[t]}
\put(1.500000,0){\oval(2.900000,1.400000)[t]}
\put(1.500000,0){\oval(3.100000,1.600000)[t]}
\put(1.500000,0.600000){\line(0,1){0.3}}
\put(2.000000,0){\oval(1.900000,0.900000)[b]}
\put(2.000000,0){\oval(2.100000,1.100000)[b]}
\put(2,-0.05){\line(1,0){1}}
 \put(2,0.05){\line(1,0){1}}
 \end{picture}
 \begin{picture}(5,3.500000)
\multiput(0,0)(1,0){4}{\circle*{0.3}}
\scriptsize\put(-1,0){4}
{\scriptsize {\put(-1,-1.250000){(4,4 ; 1,1,2,2)}}}
\put(0,-0.05){\line(1,0){1}}
 \put(0,0.05){\line(1,0){1}}
 \put(1.000000,0){\oval(1.900000,0.900000)[t]}
\put(1.000000,0){\oval(2.100000,1.100000)[t]}
\put(1.500000,0){\oval(2.900000,1.400000)[t]}
\put(1.500000,0){\oval(3.100000,1.600000)[t]}
\put(1.500000,0.600000){\line(0,1){0.3}}
\put(1,0){\line(1,0){1}}
 \put(1.500000,-0.150000){\line(0,1){0.3}}
\put(2.000000,0){\oval(1.900000,0.900000)[b]}
\put(2.000000,0){\oval(2.100000,1.100000)[b]}
\put(2,-0.05){\line(1,0){1}}
 \put(2,0.05){\line(1,0){1}}
 \end{picture}

 \begin{picture}(5,3.500000)
\multiput(0,0)(1,0){4}{\circle*{0.3}}
\scriptsize\put(-1,0){5}
{\scriptsize {\put(-1,-1.250000){(4,4 ; 2,2,0,0)}}}
\put(0,-0.05){\line(1,0){1}}
 \put(0,0.05){\line(1,0){1}}
 \put(1.000000,0){\oval(2.000000,1.000000)[t]}
\put(1.500000,0){\oval(2.900000,1.400000)[t]}
\put(1.500000,0){\oval(3.100000,1.600000)[t]}
\put(1.500000,0.600000){\line(0,1){0.3}}
\put(2.000000,0){\oval(1.900000,0.900000)[b]}
\put(2.000000,0){\oval(2.100000,1.100000)[b]}
\put(2,0){\line(1,0){1}}
 \end{picture}
 \begin{picture}(5,3.500000)
\multiput(0,0)(1,0){4}{\circle*{0.3}}
\scriptsize\put(-1,0){6}
{\scriptsize {\put(-1,-1.250000){(6,5 ; 3,3,0,0)}}}
\put(1.000000,0){\oval(1.900000,0.900000)[t]}
\put(1.000000,0){\oval(2.100000,1.100000)[t]}
\put(1.500000,0){\oval(2.900000,1.400000)[t]}
\put(1.500000,0){\oval(3.100000,1.600000)[t]}
\put(1.500000,0.600000){\line(0,1){0.3}}
\put(1,-0.05){\line(1,0){1}}
 \put(1,0.05){\line(1,0){1}}
 \put(2,-0.05){\line(1,0){1}}
 \put(2,0.05){\line(1,0){1}}
 \end{picture}
 \begin{picture}(5,3.500000)
\multiput(0,0)(1,0){4}{\circle*{0.3}}
\scriptsize\put(-1,0){7}
{\scriptsize {\put(-1,-1.250000){(6,5 ; 3,3,1,1)}}}
\put(0,-0.05){\line(1,0){1}}
 \put(0,0.05){\line(1,0){1}}
 \put(1.000000,0){\oval(1.900000,0.900000)[t]}
\put(1.000000,0){\oval(2.100000,1.100000)[t]}
\put(1.000000,0.350000){\line(0,1){0.3}}
\put(1.500000,0){\oval(2.900000,1.400000)[t]}
\put(1.500000,0){\oval(3.100000,1.600000)[t]}
\put(1.500000,0.600000){\line(0,1){0.3}}
\put(1,-0.05){\line(1,0){1}}
 \put(1,0.05){\line(1,0){1}}
 \put(2.000000,0){\oval(1.900000,0.900000)[b]}
\put(2.000000,0){\oval(2.100000,1.100000)[b]}
\put(2,-0.05){\line(1,0){1}}
 \put(2,0.05){\line(1,0){1}}
 \put(2.500000,-0.150000){\line(0,1){0.3}}
\end{picture}
 \begin{picture}(5,3.500000)
\multiput(0,0)(1,0){4}{\circle*{0.3}}
\scriptsize\put(-1,0){8}
{\scriptsize {\put(-1,-1.250000){(12,6 ; 6,6,0,0)}}}
\put(0,-0.05){\line(1,0){1}}
 \put(0,0.05){\line(1,0){1}}
 \put(1.000000,0){\oval(1.900000,0.900000)[t]}
\put(1.000000,0){\oval(2.100000,1.100000)[t]}
\put(1.500000,0){\oval(2.900000,1.400000)[t]}
\put(1.500000,0){\oval(3.100000,1.600000)[t]}
\put(1.500000,0.600000){\line(0,1){0.3}}
\put(2.000000,0){\oval(1.900000,0.900000)[b]}
\put(2.000000,0){\oval(2.100000,1.100000)[b]}
\put(2,-0.05){\line(1,0){1}}
 \put(2,0.05){\line(1,0){1}}
 \end{picture}
 \begin{picture}(5,3.500000)
\multiput(0,0)(1,0){4}{\circle*{0.3}}
\scriptsize\put(-1,0){9}
{\scriptsize {\put(-1,-1.250000){(12,6 ; 6,6,1,1)}}}
\put(0,-0.05){\line(1,0){1}}
 \put(0,0.05){\line(1,0){1}}
 \put(1.000000,0){\oval(1.900000,0.900000)[t]}
\put(1.000000,0){\oval(2.100000,1.100000)[t]}
\put(1.500000,0){\oval(2.900000,1.400000)[t]}
\put(1.500000,0){\oval(3.100000,1.600000)[t]}
\put(1.500000,0.600000){\line(0,1){0.3}}
\put(1,-0.05){\line(1,0){1}}
 \put(1,0.05){\line(1,0){1}}
 \put(1.500000,-0.150000){\line(0,1){0.3}}
\put(2.000000,0){\oval(1.900000,0.900000)[b]}
\put(2.000000,0){\oval(2.100000,1.100000)[b]}
\put(2,-0.05){\line(1,0){1}}
 \put(2,0.05){\line(1,0){1}}
 \end{picture}
\\

 \begin{picture}(5,3.500000)
{\bf \put(-2,1){3}
}\multiput(0,0)(1,0){4}{\circle*{0.3}}
\scriptsize\put(-1,0){0}
\put(0,0){\line(1,0){1}}
 \put(1,0){\line(1,0){1}}
 \put(2.000000,0){\oval(2.000000,1.000000)[t]}
\put(2,0){\line(1,0){1}}
 \end{picture}
 \begin{picture}(5,3.500000)
\multiput(0,0)(1,0){4}{\circle*{0.3}}
\scriptsize\put(-1,0){1}
{\scriptsize {\put(-1,-1.250000){(2,3 ; 0,0,0,0)}}}
\put(0,0){\line(1,0){1}}
 \put(1.000000,0){\oval(2.000000,1.000000)[t]}
\put(1.500000,0){\oval(3.000000,1.500000)[t]}
\put(1.500000,0.600000){\line(0,1){0.3}}
\put(1,0){\line(1,0){1}}
 \put(2.000000,0){\oval(2.000000,1.000000)[b]}
\put(2,0){\line(1,0){1}}
 \end{picture}
\\

 \begin{picture}(5,3.500000)
{\bf \put(-2,1){4}
}\multiput(0,0)(1,0){4}{\circle*{0.3}}
\scriptsize\put(-1,0){0}
\put(0,-0.05){\line(1,0){1}}
 \put(0,0.05){\line(1,0){1}}
 \put(1,0){\line(1,0){1}}
 \put(2,-0.1){\line(1,0){1}}
 \put(2,-0.033){\line(1,0){1}}
 \put(2,0.033){\line(1,0){1}}
 \put(2,0.1){\line(1,0){1}}
 \end{picture}
 \begin{picture}(5,3.500000)
\multiput(0,0)(1,0){4}{\circle*{0.3}}
\scriptsize\put(-1,0){1}
{\scriptsize {\put(-1,-1.250000){(2,3 ; 0,0,0,0)}}}
\put(0,0){\line(1,0){1}}
 \put(1.500000,0){\oval(2.900000,1.400000)[t]}
\put(1.500000,0){\oval(3.100000,1.600000)[t]}
\put(1.500000,0.600000){\line(0,1){0.3}}
\put(1,-0.1){\line(1,0){1}}
 \put(1,-0.033){\line(1,0){1}}
 \put(1,0.033){\line(1,0){1}}
 \put(1,0.1){\line(1,0){1}}
 \put(2.000000,0){\oval(2.000000,1.000000)[t]}
\end{picture}
 \begin{picture}(5,3.500000)
\multiput(0,0)(1,0){4}{\circle*{0.3}}
\scriptsize\put(-1,0){2}
{\scriptsize {\put(-1,-1.250000){(2,3 ; 0,0,3,3)}}}
\put(1.000000,0){\oval(2.000000,1.000000)[t]}
\put(1.500000,0){\oval(3.300000,1.800000)[t]}
\put(1.500000,0){\oval(3.100000,1.600000)[t]}
\put(1.500000,0){\oval(2.900000,1.400000)[t]}
\put(1.500000,0){\oval(2.700000,1.200000)[t]}
\put(1.500000,0.600000){\line(0,1){0.3}}
\put(1,-0.05){\line(1,0){1}}
 \put(1,0.05){\line(1,0){1}}
 \put(2,0){\line(1,0){1}}
 \end{picture}
 \begin{picture}(5,3.500000)
\multiput(0,0)(1,0){4}{\circle*{0.3}}
\scriptsize\put(-1,0){3}
{\scriptsize {\put(-1,-1.250000){(4,4 ; 1,1,0,0)}}}
\put(0,0){\line(1,0){1}}
 \put(1.000000,0){\oval(2.300000,1.300000)[t]}
\put(1.000000,0){\oval(2.100000,1.100000)[t]}
\put(1.000000,0){\oval(1.900000,0.900000)[t]}
\put(1.000000,0){\oval(1.700000,0.700000)[t]}
\put(1.500000,0){\oval(3.000000,1.500000)[t]}
\put(1.500000,0.600000){\line(0,1){0.3}}
\put(2.000000,0){\oval(2.000000,1.000000)[b]}
\put(2,-0.1){\line(1,0){1}}
 \put(2,-0.033){\line(1,0){1}}
 \put(2,0.033){\line(1,0){1}}
 \put(2,0.1){\line(1,0){1}}
 \end{picture}
 \begin{picture}(5,3.500000)
\multiput(0,0)(1,0){4}{\circle*{0.3}}
\scriptsize\put(-1,0){4}
{\scriptsize {\put(-1,-1.250000){(4,4 ; 1,1,2,2)}}}
\put(0,0){\line(1,0){1}}
 \put(1.000000,0){\oval(2.000000,1.000000)[t]}
\put(1.500000,0){\oval(3.300000,1.800000)[t]}
\put(1.500000,0){\oval(3.100000,1.600000)[t]}
\put(1.500000,0){\oval(2.900000,1.400000)[t]}
\put(1.500000,0){\oval(2.700000,1.200000)[t]}
\put(1.500000,0.600000){\line(0,1){0.3}}
\put(1,-0.05){\line(1,0){1}}
 \put(1,0.05){\line(1,0){1}}
 \put(1.500000,-0.150000){\line(0,1){0.3}}
\put(2.000000,0){\oval(2.000000,1.000000)[b]}
\put(2,0){\line(1,0){1}}
 \end{picture}

 \begin{picture}(5,3.500000)
\multiput(0,0)(1,0){4}{\circle*{0.3}}
\scriptsize\put(-1,0){5}
{\scriptsize {\put(-1,-1.250000){(8,5 ; 3,3,0,0)}}}
\put(0,-0.1){\line(1,0){1}}
 \put(0,-0.033){\line(1,0){1}}
 \put(0,0.033){\line(1,0){1}}
 \put(0,0.1){\line(1,0){1}}
 \put(1.500000,0){\oval(3.300000,1.800000)[t]}
\put(1.500000,0){\oval(3.100000,1.600000)[t]}
\put(1.500000,0){\oval(2.900000,1.400000)[t]}
\put(1.500000,0){\oval(2.700000,1.200000)[t]}
\put(1.500000,0.600000){\line(0,1){0.3}}
\put(1,0){\line(1,0){1}}
 \put(2,-0.1){\line(1,0){1}}
 \put(2,-0.033){\line(1,0){1}}
 \put(2,0.033){\line(1,0){1}}
 \put(2,0.1){\line(1,0){1}}
 \end{picture}
\\

 \begin{picture}(5,3.500000)
{\bf \put(-2,1){5}
}\multiput(0,0)(1,0){4}{\circle*{0.3}}
\scriptsize\put(-1,0){0}
\put(0,0){\line(1,0){1}}
 \put(1,-0.05){\line(1,0){1}}
 \put(1,0.05){\line(1,0){1}}
 \put(2.000000,0){\oval(1.900000,0.900000)[t]}
\put(2.000000,0){\oval(2.100000,1.100000)[t]}
\end{picture}
 \begin{picture}(5,3.500000)
\multiput(0,0)(1,0){4}{\circle*{0.3}}
\scriptsize\put(-1,0){1}
{\scriptsize {\put(-1,-1.250000){(2,3 ; 0,0,0,0)}}}
\put(0,-0.05){\line(1,0){1}}
 \put(0,0.05){\line(1,0){1}}
 \put(1.000000,0){\oval(1.900000,0.900000)[t]}
\put(1.000000,0){\oval(2.100000,1.100000)[t]}
\put(1.500000,0){\oval(3.000000,1.500000)[t]}
\put(1.500000,0.600000){\line(0,1){0.3}}
\put(2.000000,0){\oval(1.900000,0.900000)[b]}
\put(2.000000,0){\oval(2.100000,1.100000)[b]}
\put(2,-0.05){\line(1,0){1}}
 \put(2,0.05){\line(1,0){1}}
 \end{picture}
 \begin{picture}(5,3.500000)
\multiput(0,0)(1,0){4}{\circle*{0.3}}
\scriptsize\put(-1,0){2}
{\scriptsize {\put(-1,-1.250000){(2,3 ; 0,0,2,2)}}}
\put(0,-0.05){\line(1,0){1}}
 \put(0,0.05){\line(1,0){1}}
 \put(1.000000,0){\oval(2.000000,1.000000)[t]}
\put(1.500000,0){\oval(2.900000,1.400000)[t]}
\put(1.500000,0){\oval(3.100000,1.600000)[t]}
\put(1.500000,0.600000){\line(0,1){0.3}}
\put(2.000000,0){\oval(1.900000,0.900000)[b]}
\put(2.000000,0){\oval(2.100000,1.100000)[b]}
\put(2,0){\line(1,0){1}}
 \end{picture}
\\

 \begin{picture}(5,3.500000)
{\bf \put(-2,1){6}
}\multiput(0,0)(1,0){4}{\circle*{0.3}}
\scriptsize\put(-1,0){0}
\put(0,0){\line(1,0){1}}
 \put(1,-0.1){\line(1,0){1}}
 \put(1,-0.033){\line(1,0){1}}
 \put(1,0.033){\line(1,0){1}}
 \put(1,0.1){\line(1,0){1}}
 \put(2,0){\line(1,0){1}}
 \end{picture}
 \begin{picture}(5,3.500000)
\multiput(0,0)(1,0){4}{\circle*{0.3}}
\scriptsize\put(-1,0){1}
{\scriptsize {\put(-1,-1.250000){(2,3 ; 0,0,0,0)}}}
\put(0,-0.1){\line(1,0){1}}
 \put(0,-0.033){\line(1,0){1}}
 \put(0,0.033){\line(1,0){1}}
 \put(0,0.1){\line(1,0){1}}
 \put(1.500000,0){\oval(3.000000,1.500000)[t]}
\put(1.500000,0.600000){\line(0,1){0.3}}
\put(1,0){\line(1,0){1}}
 \put(2.000000,0){\oval(2.300000,1.300000)[t]}
\put(2.000000,0){\oval(2.100000,1.100000)[t]}
\put(2.000000,0){\oval(1.900000,0.900000)[t]}
\put(2.000000,0){\oval(1.700000,0.700000)[t]}
\end{picture}
 \begin{picture}(5,3.500000)
\multiput(0,0)(1,0){4}{\circle*{0.3}}
\scriptsize\put(-1,0){2}
{\scriptsize {\put(-1,-1.250000){(3,4 ; 1,0,0,1)}}}
\put(1.000000,0){\oval(2.000000,1.000000)[t]}
\put(1.500000,0){\oval(3.300000,1.800000)[t]}
\put(1.500000,0){\oval(3.100000,1.600000)[t]}
\put(1.500000,0){\oval(2.900000,1.400000)[t]}
\put(1.500000,0){\oval(2.700000,1.200000)[t]}
\put(1.500000,0.600000){\line(0,1){0.3}}
\put(2.000000,0){\oval(2.300000,1.300000)[b]}
\put(2.000000,0){\oval(2.100000,1.100000)[b]}
\put(2.000000,0){\oval(1.900000,0.900000)[b]}
\put(2.000000,0){\oval(1.700000,0.700000)[b]}
\put(2,0){\line(1,0){1}}
 \end{picture}
 \begin{picture}(5,3.500000)
\multiput(0,0)(1,0){4}{\circle*{0.3}}
\scriptsize\put(-1,0){3}
{\scriptsize {\put(-1,-1.250000){(4,4 ; 1,1,2,2)}}}
\put(0,-0.1){\line(1,0){1}}
 \put(0,-0.033){\line(1,0){1}}
 \put(0,0.033){\line(1,0){1}}
 \put(0,0.1){\line(1,0){1}}
 \put(1.000000,0){\oval(2.300000,1.300000)[t]}
\put(1.000000,0){\oval(2.100000,1.100000)[t]}
\put(1.000000,0){\oval(1.900000,0.900000)[t]}
\put(1.000000,0){\oval(1.700000,0.700000)[t]}
\put(1.500000,0){\oval(3.000000,1.500000)[t]}
\put(1.500000,0.600000){\line(0,1){0.3}}
\put(1,0){\line(1,0){1}}
 \put(1.500000,-0.150000){\line(0,1){0.3}}
\put(2.000000,0){\oval(2.300000,1.300000)[b]}
\put(2.000000,0){\oval(2.100000,1.100000)[b]}
\put(2.000000,0){\oval(1.900000,0.900000)[b]}
\put(2.000000,0){\oval(1.700000,0.700000)[b]}
\put(2,-0.1){\line(1,0){1}}
 \put(2,-0.033){\line(1,0){1}}
 \put(2,0.033){\line(1,0){1}}
 \put(2,0.1){\line(1,0){1}}
 \end{picture}
 \begin{picture}(5,3.500000)
\multiput(0,0)(1,0){4}{\circle*{0.3}}
\scriptsize\put(-1,0){4}
{\scriptsize {\put(-1,-1.250000){(5,5 ; 2,1,0,0)}}}
\put(0,0){\line(1,0){1}}
 \put(1.000000,0){\oval(2.300000,1.300000)[t]}
\put(1.000000,0){\oval(2.100000,1.100000)[t]}
\put(1.000000,0){\oval(1.900000,0.900000)[t]}
\put(1.000000,0){\oval(1.700000,0.700000)[t]}
\put(1.500000,0){\oval(3.300000,1.800000)[t]}
\put(1.500000,0){\oval(3.100000,1.600000)[t]}
\put(1.500000,0){\oval(2.900000,1.400000)[t]}
\put(1.500000,0){\oval(2.700000,1.200000)[t]}
\put(1.300000,0.600000){\line(0,1){0.3}}
\put(1.700000,0.600000){\line(0,1){0.3}}
\put(2.000000,0){\oval(2.300000,1.300000)[b]}
\put(2.000000,0){\oval(2.100000,1.100000)[b]}
\put(2.000000,0){\oval(1.900000,0.900000)[b]}
\put(2.000000,0){\oval(1.700000,0.700000)[b]}
\put(2,0){\line(1,0){1}}
 \end{picture}

 \begin{picture}(5,3.500000)
\multiput(0,0)(1,0){4}{\circle*{0.3}}
\scriptsize\put(-1,0){5}
{\scriptsize {\put(-1,-1.250000){(6,5 ; 2,2,1,1)}}}
\put(0,-0.1){\line(1,0){1}}
 \put(0,-0.033){\line(1,0){1}}
 \put(0,0.033){\line(1,0){1}}
 \put(0,0.1){\line(1,0){1}}
 \put(0.500000,-0.150000){\line(0,1){0.3}}
\put(1.000000,0){\oval(2.000000,1.000000)[t]}
\put(1.500000,0){\oval(3.300000,1.800000)[t]}
\put(1.500000,0){\oval(3.100000,1.600000)[t]}
\put(1.500000,0){\oval(2.900000,1.400000)[t]}
\put(1.500000,0){\oval(2.700000,1.200000)[t]}
\put(1.500000,0.600000){\line(0,1){0.3}}
\put(1,0){\line(1,0){1}}
 \put(2.000000,0){\oval(2.300000,1.300000)[b]}
\put(2.000000,0){\oval(2.100000,1.100000)[b]}
\put(2.000000,0){\oval(1.900000,0.900000)[b]}
\put(2.000000,0){\oval(1.700000,0.700000)[b]}
\put(2.000000,-0.650000){\line(0,1){0.3}}
\put(2,0){\line(1,0){1}}
 \end{picture}
\\

 \begin{picture}(5,3.500000)
{\bf \put(-2,1){7}
}\multiput(0,0)(1,0){4}{\circle*{0.3}}
\scriptsize\put(-1,0){0}
\put(0,-0.1){\line(1,0){1}}
 \put(0,0){\line(1,0){1}}
 \put(0,0.1){\line(1,0){1}}
 \put(1,0){\line(1,0){1}}
 \put(2,-0.1){\line(1,0){1}}
 \put(2,-0.033){\line(1,0){1}}
 \put(2,0.033){\line(1,0){1}}
 \put(2,0.1){\line(1,0){1}}
 \end{picture}
 \begin{picture}(5,3.500000)
\multiput(0,0)(1,0){4}{\circle*{0.3}}
\scriptsize\put(-1,0){1}
{\scriptsize {\put(-1,-1.250000){(2,3 ; 0,0,0,0)}}}
\put(0,0){\line(1,0){1}}
 \put(1.500000,0){\oval(2.800000,1.300000)[t]}
\put(1.500000,0){\oval(3.000000,1.500000)[t]}
\put(1.500000,0){\oval(3.200000,1.700000)[t]}
\put(1.500000,0.600000){\line(0,1){0.3}}
\put(1,-0.1){\line(1,0){1}}
 \put(1,-0.033){\line(1,0){1}}
 \put(1,0.033){\line(1,0){1}}
 \put(1,0.1){\line(1,0){1}}
 \put(2.000000,0){\oval(2.000000,1.000000)[t]}
\end{picture}
 \begin{picture}(5,3.500000)
\multiput(0,0)(1,0){4}{\circle*{0.3}}
\scriptsize\put(-1,0){2}
{\scriptsize {\put(-1,-1.250000){(2,3 ; 0,0,3,3)}}}
\put(1.000000,0){\oval(2.000000,1.000000)[t]}
\put(1.500000,0){\oval(3.300000,1.800000)[t]}
\put(1.500000,0){\oval(3.100000,1.600000)[t]}
\put(1.500000,0){\oval(2.900000,1.400000)[t]}
\put(1.500000,0){\oval(2.700000,1.200000)[t]}
\put(1.500000,0.600000){\line(0,1){0.3}}
\put(1,-0.1){\line(1,0){1}}
 \put(1,0){\line(1,0){1}}
 \put(1,0.1){\line(1,0){1}}
 \put(2,0){\line(1,0){1}}
 \end{picture}
 \begin{picture}(5,3.500000)
\multiput(0,0)(1,0){4}{\circle*{0.3}}
\scriptsize\put(-1,0){3}
{\scriptsize {\put(-1,-1.250000){(4,4 ; 1,1,2,2)}}}
\put(0,0){\line(1,0){1}}
 \put(1.000000,0){\oval(2.000000,1.000000)[t]}
\put(1.500000,0){\oval(3.300000,1.800000)[t]}
\put(1.500000,0){\oval(3.100000,1.600000)[t]}
\put(1.500000,0){\oval(2.900000,1.400000)[t]}
\put(1.500000,0){\oval(2.700000,1.200000)[t]}
\put(1.500000,0.600000){\line(0,1){0.3}}
\put(1,-0.1){\line(1,0){1}}
 \put(1,0){\line(1,0){1}}
 \put(1,0.1){\line(1,0){1}}
 \put(1.500000,-0.150000){\line(0,1){0.3}}
\put(2.000000,0){\oval(2.000000,1.000000)[b]}
\put(2,0){\line(1,0){1}}
 \end{picture}
\\

 \begin{picture}(5,3.500000)
{\bf \put(-2,1){8}
}\multiput(0,0)(1,0){4}{\circle*{0.3}}
\scriptsize\put(-1,0){0}
\put(0,-0.05){\line(1,0){1}}
 \put(0,0.05){\line(1,0){1}}
 \put(1,0){\line(1,0){1}}
 \put(2.000000,0){\oval(2.000000,1.000000)[t]}
\put(2,0){\line(1,0){1}}
 \end{picture}
 \begin{picture}(5,3.500000)
\multiput(0,0)(1,0){4}{\circle*{0.3}}
\scriptsize\put(-1,0){1}
{\scriptsize {\put(-1,-1.250000){(2,3 ; 0,0,0,0)}}}
\put(0,0){\line(1,0){1}}
 \put(1.000000,0){\oval(2.000000,1.000000)[t]}
\put(1.500000,0){\oval(2.900000,1.400000)[t]}
\put(1.500000,0){\oval(3.100000,1.600000)[t]}
\put(1.500000,0.600000){\line(0,1){0.3}}
\put(1,0){\line(1,0){1}}
 \put(2.000000,0){\oval(2.000000,1.000000)[b]}
\put(2,0){\line(1,0){1}}
 \end{picture}
\\

 \begin{picture}(5,3.500000)
{\bf \put(-2,1){9}
}\multiput(0,0)(1,0){4}{\circle*{0.3}}
\scriptsize\put(-1,0){0}
\put(0,-0.1){\line(1,0){1}}
 \put(0,-0.033){\line(1,0){1}}
 \put(0,0.033){\line(1,0){1}}
 \put(0,0.1){\line(1,0){1}}
 \put(1.000000,0){\oval(2.000000,1.000000)[t]}
\put(1.500000,0){\oval(3.000000,1.500000)[t]}
\end{picture}
 \begin{picture}(5,3.500000)
\multiput(0,0)(1,0){4}{\circle*{0.3}}
\scriptsize\put(-1,0){1}
{\scriptsize {\put(-1,-1.250000){(2,3 ; 0,0,1,1)}}}
\put(0,0){\line(1,0){1}}
 \put(1.000000,0){\oval(2.000000,1.000000)[t]}
\put(1.500000,0){\oval(3.300000,1.800000)[t]}
\put(1.500000,0){\oval(3.100000,1.600000)[t]}
\put(1.500000,0){\oval(2.900000,1.400000)[t]}
\put(1.500000,0){\oval(2.700000,1.200000)[t]}
\put(1.500000,0.600000){\line(0,1){0.3}}
\put(2.000000,0){\oval(2.000000,1.000000)[b]}
\put(2,0){\line(1,0){1}}
 \end{picture}
 \begin{picture}(5,3.500000)
\multiput(0,0)(1,0){4}{\circle*{0.3}}
\scriptsize\put(-1,0){2}
{\scriptsize {\put(-1,-1.250000){(2,3 ; 0,0,2,2)}}}
\put(0,0){\line(1,0){1}}
 \put(1.000000,0){\oval(2.300000,1.300000)[t]}
\put(1.000000,0){\oval(2.100000,1.100000)[t]}
\put(1.000000,0){\oval(1.900000,0.900000)[t]}
\put(1.000000,0){\oval(1.700000,0.700000)[t]}
\put(1.500000,0){\oval(3.000000,1.500000)[t]}
\put(1.500000,0.600000){\line(0,1){0.3}}
\put(2.000000,0){\oval(2.000000,1.000000)[b]}
\put(2,-0.1){\line(1,0){1}}
 \put(2,-0.033){\line(1,0){1}}
 \put(2,0.033){\line(1,0){1}}
 \put(2,0.1){\line(1,0){1}}
 \end{picture}
 \begin{picture}(5,3.500000)
\multiput(0,0)(1,0){4}{\circle*{0.3}}
\scriptsize\put(-1,0){3}
{\scriptsize {\put(-1,-1.250000){(4,4 ; 2,2,0,0)}}}
\put(0,-0.1){\line(1,0){1}}
 \put(0,-0.033){\line(1,0){1}}
 \put(0,0.033){\line(1,0){1}}
 \put(0,0.1){\line(1,0){1}}
 \put(1.500000,0){\oval(3.300000,1.800000)[t]}
\put(1.500000,0){\oval(3.100000,1.600000)[t]}
\put(1.500000,0){\oval(2.900000,1.400000)[t]}
\put(1.500000,0){\oval(2.700000,1.200000)[t]}
\put(1.500000,0.600000){\line(0,1){0.3}}
\put(1,0){\line(1,0){1}}
 \put(2,-0.1){\line(1,0){1}}
 \put(2,-0.033){\line(1,0){1}}
 \put(2,0.033){\line(1,0){1}}
 \put(2,0.1){\line(1,0){1}}
 \end{picture}
\\

 \begin{picture}(5,3.500000)
{\bf \put(-2,1){10}
}\multiput(0,0)(1,0){4}{\circle*{0.3}}
\scriptsize\put(-1,0){0}
\put(0,-0.05){\line(1,0){1}}
 \put(0,0.05){\line(1,0){1}}
 \put(1,-0.05){\line(1,0){1}}
 \put(1,0.05){\line(1,0){1}}
 \put(2,-0.05){\line(1,0){1}}
 \put(2,0.05){\line(1,0){1}}
 \end{picture}
 \begin{picture}(5,3.500000)
\multiput(0,0)(1,0){4}{\circle*{0.3}}
\scriptsize\put(-1,0){1}
{\scriptsize {\put(-1,-1.250000){(2,3 ; 0,0,0,0)}}}
\put(0,-0.05){\line(1,0){1}}
 \put(0,0.05){\line(1,0){1}}
 \put(1.500000,0){\oval(2.900000,1.400000)[t]}
\put(1.500000,0){\oval(3.100000,1.600000)[t]}
\put(1.500000,0.600000){\line(0,1){0.3}}
\put(1,-0.05){\line(1,0){1}}
 \put(1,0.05){\line(1,0){1}}
 \put(2.000000,0){\oval(1.900000,0.900000)[t]}
\put(2.000000,0){\oval(2.100000,1.100000)[t]}
\end{picture}
 \begin{picture}(5,3.500000)
\multiput(0,0)(1,0){4}{\circle*{0.3}}
\scriptsize\put(-1,0){2}
{\scriptsize {\put(-1,-1.250000){(4,4 ; 1,1,0,0)}}}
\put(0,-0.05){\line(1,0){1}}
 \put(0,0.05){\line(1,0){1}}
 \put(1.000000,0){\oval(1.900000,0.900000)[t]}
\put(1.000000,0){\oval(2.100000,1.100000)[t]}
\put(1.500000,0){\oval(2.900000,1.400000)[t]}
\put(1.500000,0){\oval(3.100000,1.600000)[t]}
\put(1.500000,0.600000){\line(0,1){0.3}}
\put(2.000000,0){\oval(1.900000,0.900000)[b]}
\put(2.000000,0){\oval(2.100000,1.100000)[b]}
\put(2,-0.05){\line(1,0){1}}
 \put(2,0.05){\line(1,0){1}}
 \end{picture}
 \begin{picture}(5,3.500000)
\multiput(0,0)(1,0){4}{\circle*{0.3}}
\scriptsize\put(-1,0){3}
{\scriptsize {\put(-1,-1.250000){(4,4 ; 1,1,2,2)}}}
\put(0,-0.05){\line(1,0){1}}
 \put(0,0.05){\line(1,0){1}}
 \put(1.000000,0){\oval(1.900000,0.900000)[t]}
\put(1.000000,0){\oval(2.100000,1.100000)[t]}
\put(1.500000,0){\oval(2.900000,1.400000)[t]}
\put(1.500000,0){\oval(3.100000,1.600000)[t]}
\put(1.500000,0.600000){\line(0,1){0.3}}
\put(1,-0.05){\line(1,0){1}}
 \put(1,0.05){\line(1,0){1}}
 \put(1.500000,-0.150000){\line(0,1){0.3}}
\put(2.000000,0){\oval(1.900000,0.900000)[b]}
\put(2.000000,0){\oval(2.100000,1.100000)[b]}
\put(2,-0.05){\line(1,0){1}}
 \put(2,0.05){\line(1,0){1}}
 \end{picture}
\\

 \begin{picture}(5,3.500000)
{\bf \put(-2,1){11}
}\multiput(0,0)(1,0){4}{\circle*{0.3}}
\scriptsize\put(-1,0){0}
\put(0,-0.1){\line(1,0){1}}
 \put(0,-0.033){\line(1,0){1}}
 \put(0,0.033){\line(1,0){1}}
 \put(0,0.1){\line(1,0){1}}
 \put(1,0){\line(1,0){1}}
 \put(2,-0.1){\line(1,0){1}}
 \put(2,-0.033){\line(1,0){1}}
 \put(2,0.033){\line(1,0){1}}
 \put(2,0.1){\line(1,0){1}}
 \end{picture}
 \begin{picture}(5,3.500000)
\multiput(0,0)(1,0){4}{\circle*{0.3}}
\scriptsize\put(-1,0){1}
{\scriptsize {\put(-1,-1.250000){(2,3 ; 0,0,0,0)}}}
\put(0,0){\line(1,0){1}}
 \put(1.500000,0){\oval(3.300000,1.800000)[t]}
\put(1.500000,0){\oval(3.100000,1.600000)[t]}
\put(1.500000,0){\oval(2.900000,1.400000)[t]}
\put(1.500000,0){\oval(2.700000,1.200000)[t]}
\put(1.500000,0.600000){\line(0,1){0.3}}
\put(1,-0.1){\line(1,0){1}}
 \put(1,-0.033){\line(1,0){1}}
 \put(1,0.033){\line(1,0){1}}
 \put(1,0.1){\line(1,0){1}}
 \put(2.000000,0){\oval(2.000000,1.000000)[t]}
\end{picture}
 \begin{picture}(5,3.500000)
\multiput(0,0)(1,0){4}{\circle*{0.3}}
\scriptsize\put(-1,0){2}
{\scriptsize {\put(-1,-1.250000){(4,4 ; 1,1,2,2)}}}
\put(0,0){\line(1,0){1}}
 \put(1.000000,0){\oval(2.000000,1.000000)[t]}
\put(1.500000,0){\oval(3.300000,1.800000)[t]}
\put(1.500000,0){\oval(3.100000,1.600000)[t]}
\put(1.500000,0){\oval(2.900000,1.400000)[t]}
\put(1.500000,0){\oval(2.700000,1.200000)[t]}
\put(1.500000,0.600000){\line(0,1){0.3}}
\put(1,-0.1){\line(1,0){1}}
 \put(1,-0.033){\line(1,0){1}}
 \put(1,0.033){\line(1,0){1}}
 \put(1,0.1){\line(1,0){1}}
 \put(1.500000,-0.150000){\line(0,1){0.3}}
\put(2.000000,0){\oval(2.000000,1.000000)[b]}
\put(2,0){\line(1,0){1}}
 \end{picture}
\\

 \begin{picture}(5,3.500000)
{\bf \put(-2,1){12}
}\multiput(0,0)(1,0){4}{\circle*{0.3}}
\scriptsize\put(-1,0){0}
\put(0,-0.1){\line(1,0){1}}
 \put(0,0){\line(1,0){1}}
 \put(0,0.1){\line(1,0){1}}
 \put(1,0){\line(1,0){1}}
 \put(2.000000,0){\oval(2.000000,1.000000)[t]}
\put(2,0){\line(1,0){1}}
 \end{picture}
 \begin{picture}(5,3.500000)
\multiput(0,0)(1,0){4}{\circle*{0.3}}
\scriptsize\put(-1,0){1}
{\scriptsize {\put(-1,-1.250000){(2,3 ; 0,0,0,0)}}}
\put(0,0){\line(1,0){1}}
 \put(1.000000,0){\oval(2.000000,1.000000)[t]}
\put(1.500000,0){\oval(2.800000,1.300000)[t]}
\put(1.500000,0){\oval(3.000000,1.500000)[t]}
\put(1.500000,0){\oval(3.200000,1.700000)[t]}
\put(1.500000,0.600000){\line(0,1){0.3}}
\put(1,0){\line(1,0){1}}
 \put(2.000000,0){\oval(2.000000,1.000000)[b]}
\put(2,0){\line(1,0){1}}
 \end{picture}
\\

 \begin{picture}(5,3.500000)
{\bf \put(-2,1){13}
}\multiput(0,0)(1,0){4}{\circle*{0.3}}
\scriptsize\put(-1,0){0}
\put(0,-0.05){\line(1,0){1}}
 \put(0,0.05){\line(1,0){1}}
 \put(1,-0.05){\line(1,0){1}}
 \put(1,0.05){\line(1,0){1}}
 \put(2.000000,0){\oval(1.900000,0.900000)[t]}
\put(2.000000,0){\oval(2.100000,1.100000)[t]}
\end{picture}
 \begin{picture}(5,3.500000)
\multiput(0,0)(1,0){4}{\circle*{0.3}}
\scriptsize\put(-1,0){1}
{\scriptsize {\put(-1,-1.250000){(2,3 ; 0,0,0,0)}}}
\put(0,-0.05){\line(1,0){1}}
 \put(0,0.05){\line(1,0){1}}
 \put(1.000000,0){\oval(1.900000,0.900000)[t]}
\put(1.000000,0){\oval(2.100000,1.100000)[t]}
\put(1.500000,0){\oval(2.900000,1.400000)[t]}
\put(1.500000,0){\oval(3.100000,1.600000)[t]}
\put(1.500000,0.600000){\line(0,1){0.3}}
\put(2.000000,0){\oval(1.900000,0.900000)[b]}
\put(2.000000,0){\oval(2.100000,1.100000)[b]}
\put(2,-0.05){\line(1,0){1}}
 \put(2,0.05){\line(1,0){1}}
 \end{picture}
\\

 \begin{picture}(5,3.500000)
{\bf \put(-2,1){14}
}\multiput(0,0)(1,0){4}{\circle*{0.3}}
\scriptsize\put(-1,0){0}
\put(0,-0.1){\line(1,0){1}}
 \put(0,-0.033){\line(1,0){1}}
 \put(0,0.033){\line(1,0){1}}
 \put(0,0.1){\line(1,0){1}}
 \put(1,0){\line(1,0){1}}
 \put(2.000000,0){\oval(2.000000,1.000000)[t]}
\put(2,0){\line(1,0){1}}
 \end{picture}
 \begin{picture}(5,3.500000)
\multiput(0,0)(1,0){4}{\circle*{0.3}}
\scriptsize\put(-1,0){1}
{\scriptsize {\put(-1,-1.250000){(2,3 ; 0,0,0,0)}}}
\put(0,0){\line(1,0){1}}
 \put(1.000000,0){\oval(2.000000,1.000000)[t]}
\put(1.500000,0){\oval(3.300000,1.800000)[t]}
\put(1.500000,0){\oval(3.100000,1.600000)[t]}
\put(1.500000,0){\oval(2.900000,1.400000)[t]}
\put(1.500000,0){\oval(2.700000,1.200000)[t]}
\put(1.500000,0.600000){\line(0,1){0.3}}
\put(1,0){\line(1,0){1}}
 \put(2.000000,0){\oval(2.000000,1.000000)[b]}
\put(2,0){\line(1,0){1}}
 \end{picture}

%% file: pic/f3.tex
%TeX representation of Charton picture (version 1.96)
%Copyright (c) S.Trifonov, D.Tejblum, 1993-4
%FILE: F3.tex
%DATE: Sat Nov 30 15:58:50 2002
%------------
\setlength{\unitlength}{0.075in}
\begin{picture}(70.03,28.78)
\put(23.08,1.83){{\setbox0=\hbox{$\scriptstyle\bullet$}\kern-.4\wd0\lower.5\ht0\box0}}
\put(26.15,1.83){{\setbox0=\hbox{$\scriptstyle\bullet$}\kern-.4\wd0\lower.5\ht0\box0}}
\put(26.15,4.83){{\setbox0=\hbox{$\scriptstyle\bullet$}\kern-.4\wd0\lower.5\ht0\box0}}
\put(23.08,4.83){{\setbox0=\hbox{$\scriptstyle\bullet$}\kern-.4\wd0\lower.5\ht0\box0}}
\put(12.95,5.00){{\setbox0=\hbox{$\scriptstyle\bullet$}\kern-.4\wd0\lower.5\ht0\box0}}
\put(12.95,1.66){{\setbox0=\hbox{$\scriptstyle\bullet$}\kern-.4\wd0\lower.5\ht0\box0}}
\put(2.95,7.83){{\setbox0=\hbox{$N=5$}\lower\ht0\box0}}
\put(18.28,4.00){{\setbox0=\hbox{$P=$}\lower\ht0\box0}}
\put(2.68,4.00){{\setbox0=\hbox{$F=$}\lower\ht0\box0}}
\special{em:linewidth 0.014in}
\put(23.03,4.98){\special{em:moveto}}
\put(26.23,1.86){\special{em:lineto}}
\put(23.03,4.98){\special{em:moveto}}
\put(26.23,4.98){\special{em:lineto}}
\put(26.23,1.86){\special{em:lineto}}
\put(23.03,1.86){\special{em:lineto}}
\put(23.03,4.98){\special{em:lineto}}
\put(7.51,3.21){{\setbox0=\hbox{$\scriptstyle\bullet$}\kern-.4\wd0\lower.5\ht0\box0}}
\put(7.51,3.21){\special{em:moveto}}
\put(10.00,3.21){\special{em:lineto}}
\put(12.95,5.13){\special{em:lineto}}
\put(12.95,1.73){\special{em:lineto}}
\put(10.00,3.21){\special{em:lineto}}
\put(61.08,1.83){{\setbox0=\hbox{$\scriptstyle\bullet$}\kern-.4\wd0\lower.5\ht0\box0}}
\put(64.15,1.83){{\setbox0=\hbox{$\scriptstyle\bullet$}\kern-.4\wd0\lower.5\ht0\box0}}
\put(64.15,4.83){{\setbox0=\hbox{$\scriptstyle\bullet$}\kern-.4\wd0\lower.5\ht0\box0}}
\put(61.08,4.83){{\setbox0=\hbox{$\scriptstyle\bullet$}\kern-.4\wd0\lower.5\ht0\box0}}
\put(50.95,5.00){{\setbox0=\hbox{$\scriptstyle\bullet$}\kern-.4\wd0\lower.5\ht0\box0}}
\put(50.95,1.66){{\setbox0=\hbox{$\scriptstyle\bullet$}\kern-.4\wd0\lower.5\ht0\box0}}
\put(40.68,7.16){{\setbox0=\hbox{$N=12$}\lower\ht0\box0}}
\put(56.28,4.00){{\setbox0=\hbox{$P=$}\lower\ht0\box0}}
\put(40.68,4.00){{\setbox0=\hbox{$F=$}\lower\ht0\box0}}
\put(61.03,4.98){\special{em:moveto}}
\put(64.23,1.86){\special{em:lineto}}
\put(64.23,4.98){\special{em:moveto}}
\put(61.03,1.86){\special{em:lineto}}
\put(61.03,4.98){\special{em:moveto}}
\put(64.23,4.98){\special{em:lineto}}
\put(64.23,1.86){\special{em:lineto}}
\put(61.03,1.86){\special{em:lineto}}
\put(61.03,4.98){\special{em:lineto}}
\put(45.51,3.21){{\setbox0=\hbox{$\scriptstyle\bullet$}\kern-.4\wd0\lower.5\ht0\box0}}
\put(45.51,3.21){\special{em:moveto}}
\put(48.00,3.21){\special{em:lineto}}
\put(50.95,5.13){\special{em:lineto}}
\put(50.95,1.73){\special{em:lineto}}
\put(48.00,3.21){\special{em:lineto}}
\put(56.20,9.66){\special{em:moveto}}
\put(56.20,9.73){\special{em:lineto}}
\put(56.20,9.80){\special{em:lineto}}
\put(56.21,9.90){\special{em:lineto}}
\put(56.26,10.03){\special{em:lineto}}
\put(56.31,10.13){\special{em:lineto}}
\put(56.36,10.20){\special{em:lineto}}
\put(56.41,10.26){\special{em:lineto}}
\put(56.48,10.33){\special{em:lineto}}
\put(56.55,10.40){\special{em:lineto}}
\put(56.65,10.46){\special{em:lineto}}
\put(56.78,10.53){\special{em:lineto}}
\put(57.11,10.58){\special{em:lineto}}
\put(57.48,10.51){\special{em:lineto}}
\put(57.61,10.45){\special{em:lineto}}
\put(57.70,10.38){\special{em:lineto}}
\put(57.76,10.31){\special{em:lineto}}
\put(57.83,10.25){\special{em:lineto}}
\put(57.88,10.18){\special{em:lineto}}
\put(57.91,10.11){\special{em:lineto}}
\put(57.95,10.05){\special{em:lineto}}
\put(57.98,9.98){\special{em:lineto}}
\put(58.00,9.91){\special{em:lineto}}
\put(58.01,9.85){\special{em:lineto}}
\put(58.03,9.78){\special{em:lineto}}
\put(58.03,9.71){\special{em:lineto}}
\put(58.03,9.66){\special{em:lineto}}
\put(58.03,9.60){\special{em:lineto}}
\put(58.03,9.53){\special{em:lineto}}
\put(58.01,9.46){\special{em:lineto}}
\put(58.00,9.40){\special{em:lineto}}
\put(57.98,9.33){\special{em:lineto}}
\put(57.95,9.26){\special{em:lineto}}
\put(57.91,9.20){\special{em:lineto}}
\put(57.86,9.13){\special{em:lineto}}
\put(57.81,9.06){\special{em:lineto}}
\put(57.75,9.00){\special{em:lineto}}
\put(57.68,8.93){\special{em:lineto}}
\put(57.58,8.86){\special{em:lineto}}
\put(57.45,8.80){\special{em:lineto}}
\put(57.11,8.75){\special{em:lineto}}
\put(56.75,8.81){\special{em:lineto}}
\put(56.61,8.88){\special{em:lineto}}
\put(56.53,8.95){\special{em:lineto}}
\put(56.46,9.01){\special{em:lineto}}
\put(56.40,9.08){\special{em:lineto}}
\put(56.35,9.15){\special{em:lineto}}
\put(56.31,9.21){\special{em:lineto}}
\put(56.28,9.28){\special{em:lineto}}
\put(56.25,9.35){\special{em:lineto}}
\put(56.23,9.41){\special{em:lineto}}
\put(56.21,9.48){\special{em:lineto}}
\put(56.20,9.55){\special{em:lineto}}
\put(56.20,9.61){\special{em:lineto}}
\put(56.20,9.66){\special{em:lineto}}
\put(68.40,20.18){\special{em:moveto}}
\put(68.40,20.25){\special{em:lineto}}
\put(68.40,20.31){\special{em:lineto}}
\put(68.41,20.38){\special{em:lineto}}
\put(68.43,20.45){\special{em:lineto}}
\put(68.46,20.51){\special{em:lineto}}
\put(68.50,20.58){\special{em:lineto}}
\put(68.53,20.65){\special{em:lineto}}
\put(68.60,20.71){\special{em:lineto}}
\put(68.65,20.78){\special{em:lineto}}
\put(68.75,20.86){\special{em:lineto}}
\put(68.91,20.95){\special{em:lineto}}
\put(69.21,21.00){\special{em:lineto}}
\put(69.56,20.93){\special{em:lineto}}
\put(69.68,20.86){\special{em:lineto}}
\put(69.76,20.80){\special{em:lineto}}
\put(69.83,20.73){\special{em:lineto}}
\put(69.88,20.66){\special{em:lineto}}
\put(69.91,20.60){\special{em:lineto}}
\put(69.96,20.53){\special{em:lineto}}
\put(69.98,20.46){\special{em:lineto}}
\put(70.00,20.40){\special{em:lineto}}
\put(70.01,20.33){\special{em:lineto}}
\put(70.03,20.26){\special{em:lineto}}
\put(70.03,20.20){\special{em:lineto}}
\put(70.03,20.18){\special{em:lineto}}
\put(70.03,20.11){\special{em:lineto}}
\put(70.03,20.05){\special{em:lineto}}
\put(70.01,19.98){\special{em:lineto}}
\put(70.00,19.91){\special{em:lineto}}
\put(69.96,19.85){\special{em:lineto}}
\put(69.93,19.78){\special{em:lineto}}
\put(69.90,19.71){\special{em:lineto}}
\put(69.83,19.65){\special{em:lineto}}
\put(69.78,19.58){\special{em:lineto}}
\put(69.70,19.51){\special{em:lineto}}
\put(69.60,19.45){\special{em:lineto}}
\put(69.41,19.38){\special{em:lineto}}
\put(69.21,19.36){\special{em:lineto}}
\put(68.86,19.43){\special{em:lineto}}
\put(68.75,19.50){\special{em:lineto}}
\put(68.66,19.56){\special{em:lineto}}
\put(68.60,19.63){\special{em:lineto}}
\put(68.55,19.70){\special{em:lineto}}
\put(68.51,19.76){\special{em:lineto}}
\put(68.46,19.83){\special{em:lineto}}
\put(68.45,19.90){\special{em:lineto}}
\put(68.43,19.96){\special{em:lineto}}
\put(68.41,20.03){\special{em:lineto}}
\put(68.40,20.10){\special{em:lineto}}
\put(68.40,20.16){\special{em:lineto}}
\put(68.40,20.18){\special{em:lineto}}
\put(48.26,20.18){\special{em:moveto}}
\put(48.28,20.33){\special{em:lineto}}
\put(48.33,20.48){\special{em:lineto}}
\put(48.35,20.55){\special{em:lineto}}
\put(48.38,20.61){\special{em:lineto}}
\put(48.43,20.70){\special{em:lineto}}
\put(48.55,20.83){\special{em:lineto}}
\put(48.65,20.91){\special{em:lineto}}
\put(48.80,21.00){\special{em:lineto}}
\put(49.00,21.06){\special{em:lineto}}
\put(49.16,21.08){\special{em:lineto}}
\put(49.50,21.01){\special{em:lineto}}
\put(49.63,20.95){\special{em:lineto}}
\put(49.71,20.88){\special{em:lineto}}
\put(49.80,20.81){\special{em:lineto}}
\put(49.85,20.75){\special{em:lineto}}
\put(49.90,20.68){\special{em:lineto}}
\put(49.95,20.61){\special{em:lineto}}
\put(49.98,20.55){\special{em:lineto}}
\put(50.00,20.48){\special{em:lineto}}
\put(50.03,20.41){\special{em:lineto}}
\put(50.05,20.35){\special{em:lineto}}
\put(50.05,20.28){\special{em:lineto}}
\put(50.05,20.21){\special{em:lineto}}
\put(50.06,20.18){\special{em:lineto}}
\put(50.05,20.11){\special{em:lineto}}
\put(50.05,20.05){\special{em:lineto}}
\put(50.03,19.98){\special{em:lineto}}
\put(50.01,19.91){\special{em:lineto}}
\put(50.00,19.85){\special{em:lineto}}
\put(49.96,19.78){\special{em:lineto}}
\put(49.93,19.71){\special{em:lineto}}
\put(49.88,19.65){\special{em:lineto}}
\put(49.83,19.58){\special{em:lineto}}
\put(49.76,19.51){\special{em:lineto}}
\put(49.68,19.45){\special{em:lineto}}
\put(49.56,19.38){\special{em:lineto}}
\put(49.40,19.31){\special{em:lineto}}
\put(49.16,19.28){\special{em:lineto}}
\put(48.83,19.35){\special{em:lineto}}
\put(48.70,19.41){\special{em:lineto}}
\put(48.61,19.48){\special{em:lineto}}
\put(48.53,19.55){\special{em:lineto}}
\put(48.48,19.61){\special{em:lineto}}
\put(48.43,19.68){\special{em:lineto}}
\put(48.38,19.75){\special{em:lineto}}
\put(48.35,19.81){\special{em:lineto}}
\put(48.33,19.88){\special{em:lineto}}
\put(48.30,19.95){\special{em:lineto}}
\put(48.28,20.08){\special{em:lineto}}
\put(48.26,20.18){\special{em:lineto}}
\put(57.76,27.90){\special{em:moveto}}
\put(57.76,27.96){\special{em:lineto}}
\put(57.78,28.10){\special{em:lineto}}
\put(57.83,28.23){\special{em:lineto}}
\put(57.86,28.30){\special{em:lineto}}
\put(57.90,28.36){\special{em:lineto}}
\put(57.95,28.43){\special{em:lineto}}
\put(58.00,28.50){\special{em:lineto}}
\put(58.06,28.56){\special{em:lineto}}
\put(58.15,28.63){\special{em:lineto}}
\put(58.26,28.70){\special{em:lineto}}
\put(58.45,28.76){\special{em:lineto}}
\put(58.65,28.78){\special{em:lineto}}
\put(59.00,28.71){\special{em:lineto}}
\put(59.11,28.65){\special{em:lineto}}
\put(59.21,28.58){\special{em:lineto}}
\put(59.28,28.51){\special{em:lineto}}
\put(59.35,28.45){\special{em:lineto}}
\put(59.38,28.38){\special{em:lineto}}
\put(59.43,28.31){\special{em:lineto}}
\put(59.46,28.25){\special{em:lineto}}
\put(59.48,28.18){\special{em:lineto}}
\put(59.50,28.11){\special{em:lineto}}
\put(59.51,28.05){\special{em:lineto}}
\put(59.53,27.98){\special{em:lineto}}
\put(59.53,27.91){\special{em:lineto}}
\put(59.53,27.90){\special{em:lineto}}
\put(59.53,27.83){\special{em:lineto}}
\put(59.51,27.76){\special{em:lineto}}
\put(59.51,27.70){\special{em:lineto}}
\put(59.50,27.63){\special{em:lineto}}
\put(59.46,27.56){\special{em:lineto}}
\put(59.43,27.50){\special{em:lineto}}
\put(59.40,27.43){\special{em:lineto}}
\put(59.35,27.36){\special{em:lineto}}
\put(59.30,27.30){\special{em:lineto}}
\put(59.23,27.23){\special{em:lineto}}
\put(59.15,27.16){\special{em:lineto}}
\put(59.03,27.10){\special{em:lineto}}
\put(58.85,27.03){\special{em:lineto}}
\put(58.65,27.01){\special{em:lineto}}
\put(58.30,27.08){\special{em:lineto}}
\put(58.18,27.15){\special{em:lineto}}
\put(58.08,27.21){\special{em:lineto}}
\put(58.01,27.28){\special{em:lineto}}
\put(57.95,27.35){\special{em:lineto}}
\put(57.91,27.41){\special{em:lineto}}
\put(57.86,27.50){\special{em:lineto}}
\put(57.83,27.56){\special{em:lineto}}
\put(57.80,27.63){\special{em:lineto}}
\put(57.78,27.76){\special{em:lineto}}
\put(57.76,27.85){\special{em:lineto}}
\put(57.76,27.90){\special{em:lineto}}
\put(10.28,15.45){\special{em:moveto}}
\put(10.30,15.51){\special{em:lineto}}
\put(10.30,15.58){\special{em:lineto}}
\put(10.31,15.65){\special{em:lineto}}
\put(10.33,15.71){\special{em:lineto}}
\put(10.36,15.78){\special{em:lineto}}
\put(10.40,15.85){\special{em:lineto}}
\put(10.43,15.91){\special{em:lineto}}
\put(10.48,15.98){\special{em:lineto}}
\put(10.55,16.05){\special{em:lineto}}
\put(10.61,16.11){\special{em:lineto}}
\put(10.73,16.18){\special{em:lineto}}
\put(10.88,16.25){\special{em:lineto}}
\put(11.11,16.28){\special{em:lineto}}
\put(11.43,16.21){\special{em:lineto}}
\put(11.56,16.15){\special{em:lineto}}
\put(11.66,16.06){\special{em:lineto}}
\put(11.75,15.98){\special{em:lineto}}
\put(11.81,15.90){\special{em:lineto}}
\put(11.86,15.80){\special{em:lineto}}
\put(11.90,15.73){\special{em:lineto}}
\put(11.91,15.63){\special{em:lineto}}
\put(11.93,15.56){\special{em:lineto}}
\put(11.95,15.45){\special{em:lineto}}
\put(11.93,15.28){\special{em:lineto}}
\put(11.88,15.15){\special{em:lineto}}
\put(11.86,15.08){\special{em:lineto}}
\put(11.81,15.00){\special{em:lineto}}
\put(11.76,14.93){\special{em:lineto}}
\put(11.70,14.86){\special{em:lineto}}
\put(11.63,14.80){\special{em:lineto}}
\put(11.53,14.73){\special{em:lineto}}
\put(11.40,14.66){\special{em:lineto}}
\put(11.11,14.61){\special{em:lineto}}
\put(10.80,14.68){\special{em:lineto}}
\put(10.66,14.75){\special{em:lineto}}
\put(10.58,14.81){\special{em:lineto}}
\put(10.51,14.88){\special{em:lineto}}
\put(10.45,14.95){\special{em:lineto}}
\put(10.41,15.01){\special{em:lineto}}
\put(10.36,15.08){\special{em:lineto}}
\put(10.35,15.15){\special{em:lineto}}
\put(10.31,15.21){\special{em:lineto}}
\put(10.30,15.28){\special{em:lineto}}
\put(10.30,15.35){\special{em:lineto}}
\put(10.28,15.41){\special{em:lineto}}
\put(10.28,15.45){\special{em:lineto}}
\put(23.00,26.26){\special{em:moveto}}
\put(23.00,26.33){\special{em:lineto}}
\put(23.00,26.40){\special{em:lineto}}
\put(23.01,26.46){\special{em:lineto}}
\put(23.03,26.53){\special{em:lineto}}
\put(23.06,26.60){\special{em:lineto}}
\put(23.10,26.66){\special{em:lineto}}
\put(23.13,26.73){\special{em:lineto}}
\put(23.20,26.80){\special{em:lineto}}
\put(23.25,26.86){\special{em:lineto}}
\put(23.33,26.93){\special{em:lineto}}
\put(23.43,27.00){\special{em:lineto}}
\put(23.61,27.06){\special{em:lineto}}
\put(23.81,27.08){\special{em:lineto}}
\put(24.16,27.01){\special{em:lineto}}
\put(24.28,26.95){\special{em:lineto}}
\put(24.36,26.88){\special{em:lineto}}
\put(24.43,26.81){\special{em:lineto}}
\put(24.48,26.75){\special{em:lineto}}
\put(24.51,26.68){\special{em:lineto}}
\put(24.56,26.61){\special{em:lineto}}
\put(24.58,26.55){\special{em:lineto}}
\put(24.60,26.48){\special{em:lineto}}
\put(24.61,26.41){\special{em:lineto}}
\put(24.63,26.35){\special{em:lineto}}
\put(24.63,26.28){\special{em:lineto}}
\put(24.63,26.26){\special{em:lineto}}
\put(24.63,26.20){\special{em:lineto}}
\put(24.63,26.13){\special{em:lineto}}
\put(24.61,26.06){\special{em:lineto}}
\put(24.60,26.00){\special{em:lineto}}
\put(24.56,25.93){\special{em:lineto}}
\put(24.53,25.86){\special{em:lineto}}
\put(24.50,25.80){\special{em:lineto}}
\put(24.43,25.73){\special{em:lineto}}
\put(24.38,25.66){\special{em:lineto}}
\put(24.28,25.58){\special{em:lineto}}
\put(24.11,25.50){\special{em:lineto}}
\put(23.81,25.45){\special{em:lineto}}
\put(23.46,25.51){\special{em:lineto}}
\put(23.35,25.58){\special{em:lineto}}
\put(23.26,25.65){\special{em:lineto}}
\put(23.20,25.71){\special{em:lineto}}
\put(23.15,25.78){\special{em:lineto}}
\put(23.11,25.85){\special{em:lineto}}
\put(23.06,25.91){\special{em:lineto}}
\put(23.05,25.98){\special{em:lineto}}
\put(23.03,26.05){\special{em:lineto}}
\put(23.01,26.11){\special{em:lineto}}
\put(23.00,26.18){\special{em:lineto}}
\put(23.00,26.25){\special{em:lineto}}
\put(23.00,26.26){\special{em:lineto}}
\put(38.01,0.00){{\setbox0=\hbox{b)}\lower\ht0\box0}}
\put(0.01,0.00){{\setbox0=\hbox{a)}\lower\ht0\box0}}
\put(58.65,27.75){{\setbox0=\hbox{$\scriptstyle\bullet$}\kern-.4\wd0\lower.5\ht0\box0}}
\put(49.16,20.18){{\setbox0=\hbox{$\scriptstyle\bullet$}\kern-.4\wd0\lower.5\ht0\box0}}
\put(69.21,20.18){{\setbox0=\hbox{$\scriptstyle\bullet$}\kern-.4\wd0\lower.5\ht0\box0}}
\put(57.11,9.66){{\setbox0=\hbox{$\scriptstyle\bullet$}\kern-.4\wd0\lower.5\ht0\box0}}
\put(55.10,12.33){\special{em:moveto}}
\put(55.38,12.48){\special{em:lineto}}
\put(55.66,12.63){\special{em:moveto}}
\put(55.95,12.80){\special{em:lineto}}
\put(56.23,12.96){\special{em:moveto}}
\put(56.51,13.11){\special{em:lineto}}
\put(56.80,13.28){\special{em:moveto}}
\put(57.08,13.45){\special{em:lineto}}
\put(57.36,13.61){\special{em:moveto}}
\put(57.46,13.66){\special{em:lineto}}
\put(57.68,13.61){\special{em:lineto}}
\put(58.00,13.53){\special{em:moveto}}
\put(58.31,13.45){\special{em:lineto}}
\put(58.63,13.36){\special{em:moveto}}
\put(58.95,13.28){\special{em:lineto}}
\put(59.26,13.20){\special{em:moveto}}
\put(59.58,13.11){\special{em:lineto}}
\put(59.90,13.03){\special{em:moveto}}
\put(60.21,12.93){\special{em:lineto}}
\put(60.53,12.83){\special{em:moveto}}
\put(60.85,12.73){\special{em:lineto}}
\put(58.65,27.60){\special{em:moveto}}
\put(57.11,9.51){\special{em:lineto}}
\put(61.50,20.18){\special{em:moveto}}
\put(61.78,20.35){\special{em:lineto}}
\put(62.06,20.51){\special{em:moveto}}
\put(62.35,20.68){\special{em:lineto}}
\put(62.63,20.85){\special{em:moveto}}
\put(62.91,21.01){\special{em:lineto}}
\put(63.20,21.18){\special{em:moveto}}
\put(63.48,21.35){\special{em:lineto}}
\put(63.76,21.51){\special{em:moveto}}
\put(64.05,21.68){\special{em:lineto}}
\put(64.33,21.85){\special{em:moveto}}
\put(64.61,22.01){\special{em:lineto}}
\put(64.90,22.18){\special{em:moveto}}
\put(65.18,22.35){\special{em:lineto}}
\put(65.46,22.51){\special{em:moveto}}
\put(65.76,22.68){\special{em:lineto}}
\put(65.76,22.38){\special{em:moveto}}
\put(65.75,22.05){\special{em:lineto}}
\put(65.73,21.71){\special{em:moveto}}
\put(65.71,21.38){\special{em:lineto}}
\put(65.70,21.05){\special{em:moveto}}
\put(65.68,20.71){\special{em:lineto}}
\put(65.66,20.38){\special{em:moveto}}
\put(65.65,20.05){\special{em:lineto}}
\put(65.63,19.71){\special{em:moveto}}
\put(65.61,19.38){\special{em:lineto}}
\put(65.60,19.05){\special{em:moveto}}
\put(65.58,18.71){\special{em:lineto}}
\put(65.56,18.38){\special{em:moveto}}
\put(65.55,18.05){\special{em:lineto}}
\put(65.53,17.71){\special{em:moveto}}
\put(65.50,17.38){\special{em:lineto}}
\put(65.46,17.05){\special{em:moveto}}
\put(65.43,16.71){\special{em:lineto}}
\put(51.76,16.63){\special{em:moveto}}
\put(51.76,16.96){\special{em:lineto}}
\put(51.76,17.30){\special{em:moveto}}
\put(51.76,17.63){\special{em:lineto}}
\put(51.76,17.96){\special{em:moveto}}
\put(51.76,18.30){\special{em:lineto}}
\put(51.76,18.63){\special{em:moveto}}
\put(51.76,18.96){\special{em:lineto}}
\put(51.76,19.30){\special{em:moveto}}
\put(51.76,19.63){\special{em:lineto}}
\put(51.76,19.96){\special{em:moveto}}
\put(51.76,20.30){\special{em:lineto}}
\put(51.76,20.63){\special{em:moveto}}
\put(51.76,20.96){\special{em:lineto}}
\put(51.76,21.30){\special{em:moveto}}
\put(51.76,21.63){\special{em:lineto}}
\put(51.76,21.96){\special{em:moveto}}
\put(51.76,22.26){\special{em:lineto}}
\put(51.78,22.26){\special{em:lineto}}
\put(52.08,22.13){\special{em:moveto}}
\put(52.38,22.00){\special{em:lineto}}
\put(52.68,21.86){\special{em:moveto}}
\put(52.98,21.73){\special{em:lineto}}
\put(53.28,21.60){\special{em:moveto}}
\put(53.58,21.46){\special{em:lineto}}
\put(53.86,21.33){\special{em:moveto}}
\put(54.15,21.20){\special{em:lineto}}
\put(54.43,21.05){\special{em:moveto}}
\put(54.71,20.91){\special{em:lineto}}
\put(55.00,20.76){\special{em:moveto}}
\put(55.28,20.63){\special{em:lineto}}
\put(55.56,20.48){\special{em:moveto}}
\put(55.86,20.33){\special{em:lineto}}
\put(58.65,27.75){\special{em:moveto}}
\put(49.05,20.18){\special{em:lineto}}
\put(58.65,27.75){\special{em:moveto}}
\put(69.33,20.18){\special{em:lineto}}
\put(54.98,12.33){\special{em:moveto}}
\put(57.11,9.51){\special{em:lineto}}
\put(61.15,12.48){\special{em:lineto}}
\put(56.28,20.18){\special{em:moveto}}
\put(49.05,20.18){\special{em:lineto}}
\put(51.76,16.50){\special{em:lineto}}
\put(61.38,20.18){\special{em:moveto}}
\put(69.33,20.18){\special{em:lineto}}
\put(65.41,16.63){\special{em:lineto}}
\put(56.28,20.18){\special{em:moveto}}
\put(58.65,27.75){\special{em:lineto}}
\put(61.15,12.48){\special{em:lineto}}
\put(54.98,12.33){\special{em:moveto}}
\put(58.65,27.75){\special{em:lineto}}
\put(61.38,20.18){\special{em:lineto}}
\put(51.76,16.50){\special{em:moveto}}
\put(58.65,27.75){\special{em:lineto}}
\put(65.41,16.50){\special{em:lineto}}
\put(58.65,16.63){\special{em:moveto}}
\put(58.65,27.60){\special{em:lineto}}
\put(61.15,12.48){\special{em:moveto}}
\put(61.11,12.55){\special{em:lineto}}
\put(60.98,12.75){\special{em:moveto}}
\put(60.95,12.81){\special{em:lineto}}
\put(60.81,13.01){\special{em:moveto}}
\put(60.78,13.08){\special{em:lineto}}
\put(60.65,13.28){\special{em:moveto}}
\put(60.61,13.35){\special{em:lineto}}
\put(60.48,13.55){\special{em:moveto}}
\put(60.45,13.61){\special{em:lineto}}
\put(60.31,13.81){\special{em:moveto}}
\put(60.28,13.88){\special{em:lineto}}
\put(60.15,14.08){\special{em:moveto}}
\put(60.11,14.15){\special{em:lineto}}
\put(59.98,14.35){\special{em:moveto}}
\put(59.95,14.41){\special{em:lineto}}
\put(59.81,14.61){\special{em:moveto}}
\put(59.78,14.68){\special{em:lineto}}
\put(59.65,14.88){\special{em:moveto}}
\put(59.61,14.95){\special{em:lineto}}
\put(59.48,15.15){\special{em:moveto}}
\put(59.45,15.21){\special{em:lineto}}
\put(59.31,15.41){\special{em:moveto}}
\put(59.28,15.48){\special{em:lineto}}
\put(59.15,15.68){\special{em:moveto}}
\put(59.11,15.75){\special{em:lineto}}
\put(58.98,15.95){\special{em:moveto}}
\put(58.95,16.01){\special{em:lineto}}
\put(58.81,16.21){\special{em:moveto}}
\put(58.78,16.28){\special{em:lineto}}
\put(58.65,16.48){\special{em:moveto}}
\put(58.61,16.55){\special{em:lineto}}
\put(58.48,16.75){\special{em:moveto}}
\put(58.45,16.81){\special{em:lineto}}
\put(58.31,17.01){\special{em:moveto}}
\put(58.28,17.08){\special{em:lineto}}
\put(58.15,17.28){\special{em:moveto}}
\put(58.11,17.35){\special{em:lineto}}
\put(57.98,17.55){\special{em:moveto}}
\put(57.95,17.61){\special{em:lineto}}
\put(57.81,17.81){\special{em:moveto}}
\put(57.78,17.88){\special{em:lineto}}
\put(57.65,18.08){\special{em:moveto}}
\put(57.61,18.15){\special{em:lineto}}
\put(57.48,18.35){\special{em:moveto}}
\put(57.45,18.41){\special{em:lineto}}
\put(57.31,18.61){\special{em:moveto}}
\put(57.28,18.68){\special{em:lineto}}
\put(57.15,18.88){\special{em:moveto}}
\put(57.11,18.95){\special{em:lineto}}
\put(56.98,19.15){\special{em:moveto}}
\put(56.95,19.21){\special{em:lineto}}
\put(56.81,19.41){\special{em:moveto}}
\put(56.78,19.48){\special{em:lineto}}
\put(56.65,19.68){\special{em:moveto}}
\put(56.61,19.75){\special{em:lineto}}
\put(56.46,19.95){\special{em:moveto}}
\put(56.41,20.01){\special{em:lineto}}
\put(61.38,20.18){\special{em:moveto}}
\put(61.33,20.13){\special{em:lineto}}
\put(61.18,19.95){\special{em:moveto}}
\put(61.13,19.90){\special{em:lineto}}
\put(60.98,19.71){\special{em:moveto}}
\put(60.93,19.66){\special{em:lineto}}
\put(60.78,19.48){\special{em:moveto}}
\put(60.73,19.43){\special{em:lineto}}
\put(60.58,19.25){\special{em:moveto}}
\put(60.53,19.20){\special{em:lineto}}
\put(60.38,19.01){\special{em:moveto}}
\put(60.33,18.96){\special{em:lineto}}
\put(60.18,18.78){\special{em:moveto}}
\put(60.13,18.73){\special{em:lineto}}
\put(59.98,18.55){\special{em:moveto}}
\put(59.93,18.50){\special{em:lineto}}
\put(59.78,18.31){\special{em:moveto}}
\put(59.73,18.26){\special{em:lineto}}
\put(59.58,18.08){\special{em:moveto}}
\put(59.53,18.03){\special{em:lineto}}
\put(59.38,17.85){\special{em:moveto}}
\put(59.33,17.80){\special{em:lineto}}
\put(59.18,17.61){\special{em:moveto}}
\put(59.13,17.56){\special{em:lineto}}
\put(58.98,17.38){\special{em:moveto}}
\put(58.93,17.33){\special{em:lineto}}
\put(58.78,17.15){\special{em:moveto}}
\put(58.73,17.10){\special{em:lineto}}
\put(58.58,16.91){\special{em:moveto}}
\put(58.53,16.86){\special{em:lineto}}
\put(58.38,16.68){\special{em:moveto}}
\put(58.33,16.63){\special{em:lineto}}
\put(58.18,16.45){\special{em:moveto}}
\put(58.13,16.40){\special{em:lineto}}
\put(57.98,16.21){\special{em:moveto}}
\put(57.93,16.16){\special{em:lineto}}
\put(57.78,15.98){\special{em:moveto}}
\put(57.73,15.93){\special{em:lineto}}
\put(57.58,15.75){\special{em:moveto}}
\put(57.53,15.70){\special{em:lineto}}
\put(57.38,15.51){\special{em:moveto}}
\put(57.33,15.45){\special{em:lineto}}
\put(57.18,15.25){\special{em:moveto}}
\put(57.13,15.18){\special{em:lineto}}
\put(56.98,14.98){\special{em:moveto}}
\put(56.93,14.91){\special{em:lineto}}
\put(56.78,14.71){\special{em:moveto}}
\put(56.73,14.65){\special{em:lineto}}
\put(56.58,14.45){\special{em:moveto}}
\put(56.53,14.38){\special{em:lineto}}
\put(56.38,14.18){\special{em:moveto}}
\put(56.33,14.11){\special{em:lineto}}
\put(56.18,13.91){\special{em:moveto}}
\put(56.13,13.85){\special{em:lineto}}
\put(55.98,13.65){\special{em:moveto}}
\put(55.93,13.58){\special{em:lineto}}
\put(55.78,13.38){\special{em:moveto}}
\put(55.73,13.31){\special{em:lineto}}
\put(55.58,13.11){\special{em:moveto}}
\put(55.53,13.05){\special{em:lineto}}
\put(55.38,12.85){\special{em:moveto}}
\put(55.33,12.78){\special{em:lineto}}
\put(55.18,12.58){\special{em:moveto}}
\put(55.13,12.51){\special{em:lineto}}
\put(51.76,16.50){\special{em:moveto}}
\put(51.85,16.50){\special{em:lineto}}
\put(52.10,16.50){\special{em:moveto}}
\put(52.18,16.50){\special{em:lineto}}
\put(52.43,16.50){\special{em:moveto}}
\put(52.51,16.50){\special{em:lineto}}
\put(52.76,16.50){\special{em:moveto}}
\put(52.85,16.50){\special{em:lineto}}
\put(53.10,16.50){\special{em:moveto}}
\put(53.18,16.50){\special{em:lineto}}
\put(53.43,16.50){\special{em:moveto}}
\put(53.51,16.50){\special{em:lineto}}
\put(53.76,16.50){\special{em:moveto}}
\put(53.85,16.50){\special{em:lineto}}
\put(54.10,16.50){\special{em:moveto}}
\put(54.18,16.50){\special{em:lineto}}
\put(54.43,16.50){\special{em:moveto}}
\put(54.51,16.50){\special{em:lineto}}
\put(54.76,16.50){\special{em:moveto}}
\put(54.85,16.50){\special{em:lineto}}
\put(55.10,16.50){\special{em:moveto}}
\put(55.18,16.50){\special{em:lineto}}
\put(55.43,16.50){\special{em:moveto}}
\put(55.51,16.50){\special{em:lineto}}
\put(55.76,16.50){\special{em:moveto}}
\put(55.85,16.50){\special{em:lineto}}
\put(56.10,16.50){\special{em:moveto}}
\put(56.18,16.50){\special{em:lineto}}
\put(56.43,16.50){\special{em:moveto}}
\put(56.51,16.50){\special{em:lineto}}
\put(56.76,16.50){\special{em:moveto}}
\put(56.85,16.50){\special{em:lineto}}
\put(57.10,16.50){\special{em:moveto}}
\put(57.18,16.50){\special{em:lineto}}
\put(57.43,16.50){\special{em:moveto}}
\put(57.51,16.50){\special{em:lineto}}
\put(57.76,16.50){\special{em:moveto}}
\put(57.85,16.50){\special{em:lineto}}
\put(58.10,16.50){\special{em:moveto}}
\put(58.18,16.50){\special{em:lineto}}
\put(58.43,16.50){\special{em:moveto}}
\put(58.51,16.50){\special{em:lineto}}
\put(58.76,16.50){\special{em:moveto}}
\put(58.85,16.50){\special{em:lineto}}
\put(59.10,16.50){\special{em:moveto}}
\put(59.18,16.50){\special{em:lineto}}
\put(59.43,16.50){\special{em:moveto}}
\put(59.51,16.50){\special{em:lineto}}
\put(59.76,16.50){\special{em:moveto}}
\put(59.85,16.50){\special{em:lineto}}
\put(60.10,16.50){\special{em:moveto}}
\put(60.18,16.50){\special{em:lineto}}
\put(60.43,16.50){\special{em:moveto}}
\put(60.51,16.50){\special{em:lineto}}
\put(60.76,16.50){\special{em:moveto}}
\put(60.85,16.50){\special{em:lineto}}
\put(61.10,16.50){\special{em:moveto}}
\put(61.18,16.50){\special{em:lineto}}
\put(61.43,16.50){\special{em:moveto}}
\put(61.51,16.50){\special{em:lineto}}
\put(61.76,16.50){\special{em:moveto}}
\put(61.85,16.50){\special{em:lineto}}
\put(62.10,16.50){\special{em:moveto}}
\put(62.18,16.50){\special{em:lineto}}
\put(62.43,16.50){\special{em:moveto}}
\put(62.51,16.50){\special{em:lineto}}
\put(62.76,16.50){\special{em:moveto}}
\put(62.85,16.50){\special{em:lineto}}
\put(63.10,16.50){\special{em:moveto}}
\put(63.18,16.50){\special{em:lineto}}
\put(63.43,16.50){\special{em:moveto}}
\put(63.51,16.50){\special{em:lineto}}
\put(63.76,16.50){\special{em:moveto}}
\put(63.85,16.50){\special{em:lineto}}
\put(64.10,16.50){\special{em:moveto}}
\put(64.18,16.50){\special{em:lineto}}
\put(64.43,16.50){\special{em:moveto}}
\put(64.51,16.50){\special{em:lineto}}
\put(64.76,16.50){\special{em:moveto}}
\put(64.85,16.50){\special{em:lineto}}
\put(65.10,16.50){\special{em:moveto}}
\put(65.18,16.50){\special{em:lineto}}
\put(56.28,20.18){\special{em:moveto}}
\put(61.38,20.18){\special{em:lineto}}
\put(65.41,16.50){\special{em:lineto}}
\put(61.15,12.48){\special{em:lineto}}
\put(54.98,12.48){\special{em:lineto}}
\put(51.76,16.50){\special{em:lineto}}
\put(56.28,20.18){\special{em:lineto}}
\put(11.11,15.45){{\setbox0=\hbox{$\scriptstyle\bullet$}\kern-.4\wd0\lower.5\ht0\box0}}
\put(23.81,26.26){{\setbox0=\hbox{$\scriptstyle\bullet$}\kern-.4\wd0\lower.5\ht0\box0}}
\put(24.53,14.41){\special{em:moveto}}
\put(20.50,17.08){\special{em:lineto}}
\put(15.86,19.46){\special{em:lineto}}
\put(16.70,12.93){\special{em:lineto}}
\put(24.41,14.41){\special{em:lineto}}
\put(23.81,26.26){\special{em:lineto}}
\put(18.71,16.78){\special{em:lineto}}
\put(16.81,12.93){\special{em:lineto}}
\put(20.50,17.08){\special{em:lineto}}
\put(23.81,26.26){\special{em:lineto}}
\put(23.93,26.41){\special{em:moveto}}
\put(26.43,18.26){\special{em:lineto}}
\put(11.00,15.30){\special{em:moveto}}
\put(23.93,26.41){\special{em:lineto}}
\put(22.51,10.56){\special{em:lineto}}
\put(26.43,18.26){\special{em:lineto}}
\put(11.00,15.30){\special{em:lineto}}
\put(22.51,10.56){\special{em:lineto}}
\end{picture}